\documentclass[11pt,reqno,tbtags]{amsart}
\usepackage{amssymb}
\usepackage{longtable}
\usepackage{morefloats}
\usepackage{url}
\usepackage{natbib}
\bibpunct[, ]{[}{]}{;}{n}{,}{,}

\title{Tibetan calendar mathematics}

\date{31 December, 2007 (Sunday 23, month 11, Fire--Pig year);\\
revised 8 January, 2014 (Wednesday 8, month 11,  Water--Snake year)}

\author{Svante Janson}
\address{Department of Mathematics, Uppsala University, PO Box 480,
SE-751~06 Uppsala, Sweden}
\email{svante.janson@math.uu.se}
\newcommand\urladdrx[1]{{\urladdr{\def~{{\tiny$\sim$}}#1}}}
\urladdrx{http://www.math.uu.se/~svante/}

\numberwithin{equation}{section}

\newtheorem*{property*}{Property \csname @currentlabel\endcsname}

\theoremstyle{definition}

\newtheorem*{acks}{Acknowledgements}

\theoremstyle{remark}
\newtheorem{remark}{Remark}

\newenvironment{romenumerate}[1][0pt]{
\addtolength{\leftmargini}{#1}\begin{enumerate}
 }{\end{enumerate}}

\newcounter{oldenumi}
{\setcounter{oldenumi}{\value{enumi}}
\begin{romenumerate} \setcounter{enumi}{\value{oldenumi}}}
{\end{romenumerate}}

\newcounter{thmenumerate}

\newcounter{xenumerate}   

\newcommand{\refR}[1]{Remark~\ref{#1}}
\newcommand{\refS}[1]{Section~\ref{#1}}

\newcommand{\refF}[1]{Footnote~\ref{#1}}
\newcommand{\refApp}[1]{Appendix~\ref{#1}}
\newcommand{\refand}[2]{\ref{#1} and~\ref{#2}}
\newcommand{\refT}[1]{Table~\ref{#1}}

\begingroup
  \divide\count255 by 60
  \multiply\count255 by -60
  \advance\count255 by \time
  \ifnum \count255 < 10 \xdef\klockan{\the\count1.0\the\count255}
  \else\xdef\klockan{\the\count1.\the\count255}\fi
\endgroup

\newcommand\set[1]{\ensuremath{\{#1\}}}

\newcommand\xpar[1]{(#1)}
\newcommand\bigpar[1]{\bigl(#1\bigr)}
\newcommand\Bigpar[1]{\Bigl(#1\Bigr)}

\newcommand\lrpar[1]{\left(#1\right)}

\def\rompar(#1){\textup(#1\textup)}    
\newcommand\xfrac[2]{#1/#2}

\def\xexp(#1){e^{#1}}
\newcommand\ceil[1]{\lceil#1\rceil}

\newcommand\Ceil[1]{\left\lceil#1\right\rceil}
\newcommand\floor[1]{\lfloor#1\rfloor}
\newcommand\Floor[1]{\left\lfloor#1\right\rfloor}

\newcommand\ie{i.e.\spacefactor=1000}
\newcommand\eg{e.g.\spacefactor=1000}
\newcommand\viz{viz.\spacefactor=1000}
\newcommand\cf{cf.\spacefactor=1000}

\newcommand\bbZ{\mathbb Z}

\newcommand\FRAC{\operatorname{frac}}

\newcommand\ga{\alpha}
\newcommand\gb{\beta}

\newcommand\gD{\Delta}

\newcommand\gl{\lambda}

\newcommand\gamx{\gamma^*}

\newcommand\cP{\mathcal P}
\newcommand\Bool[1]{\left\{#1\right\}} 
\newcommand\bigbool[1]{\bigl\{#1\bigr\}} 
\newcommand\boolx[1]{[#1]}

\newcommand\xfootnote[1]{\unskip\footnote{#1}}

\newcommand\rr[1]{\enspace (#1)}
\newcommand\grad{\ensuremath{^\circ}}  
\newcommand\true{\ensuremath{\mathit{true}}}  
\newcommand\false{\ensuremath{\mathit{false}}}  
\newcommand\tmc{true month count}
\newcommand\tm{true month}
\newcommand\ixx{intercalation index}
\newcommand\MM{M^*}
\newcommand\MSL{\mathit{mean\_sun}}
\newcommand\dpp{definition point}  

\newcommand\tibx[1]{\emph{#1}}
\newcommand\tibxx[1]{#1}
\newcommand\tib[1]{(\tibx{#1})}
\newcommand\Caitra{\tibx{Caitra}}
\newcommand\nag{\tibx{nag pa}}
\newcommand\Kc{K\=alacakra}
\newcommand\KT{\Kc{} Tantra}
\newcommand\PH{Phugpa}
\newcommand\TS{Tsurphu} 
\newcommand\YP{Yeshe Paljor}
\newcommand\eS{\text{E806}}
\newcommand\eH{\text{E1927}}
\newcommand\eX{\text{E1987}}
\newcommand\gbx{\beta^*}
\newcommand\gax{\alpha'}
\newcommand\JD{\mathrm{JD}} 
\newcommand\RD{\mathrm{RD}} 

\newcommand\amod{\allowbreak \mkern 8mu\operatorname{amod}\,\,}
\newcommand\bxmod{\; \bmod}

\newcommand\karana{kara{\d n}a}
\newcommand\Karana{Kara{\d n}a}
\newcommand\siddhanta{siddh\=anta}
\newcommand\Hp{[Henning, personal communication]}

\newenvironment{Rule}{\begin{quote}\em}{\end{quote}}

\begin{document}

\begin{abstract}
The calculations of the Tibetan calendar are described, using modern
mathematical notations instead of the traditional methods.
\end{abstract}

\maketitle

\section{Introduction}\label{Sintro}

The Tibetan calendar is derived from the Indian calendar tradition; it
has the same general structure as Indian calendars,
but the details differ significantly. 
The basis for the Tibetan calendar is the \KT, which was translated
from Sanskrit into Tibetan in the 11th century. 
(Traditional date of the translation is 1027
when the first 60 year cycle starts.) 
It is based on Indian astronomy, but much modified.
The calendar became the standard in Tibet in the second half of the
thirteenth century. 

As in Indian calendars 
(see \citet{CC}),
months are lunar (from new moon to new moon) but numbered according to the
corresponding solar months,   
\ie{} the position of the sun,
\xfootnote{Solar months are used only in the astronomical theory behind the
  calendar, and not in the calendar itself.} 
while days are numbered by the corresponding
lunar days. Since these correspondences are not perfect, there are occasionally
two months with the same number, in which case the first of them is
regarded as a leap month, and occasionally a skipped date or two days
with the same date
(then the first of them is regarded as a
leap day).
Unlike modern Indian calendars, there are no skipped months.

Various improvements of the calendar calculations have been suggested over the
centuries, see \cite{Schuh}, \cite[Chapter VI]{Henning},
\cite{tibetenc}, \cite{kalacakra},
and different traditions follow different rules for the details of the
calculation.
There are two main versions (\emph{\PH} and \emph{\TS})
of the Tibetan calendar in use today by
different groups inside and outside Tibet;
moreover, Bhutan and Mongolia also use versions of the Tibetan calendar.
See further \refApp{Aversions}.
The different versions frequently differ by a day or a month, see
\refApp{Acomp}. 

The description below refers to the \PH{} version,
introduced 1447, 
which is the most common version;
it is the version followed by \eg{} the Dalai Lama and 
it can be regarded as the standard version of the Tibetan calendar.
The differences 
in the \TS{} version and other versions are discussed in \refApp{Aversions}.

\begin{table}[t]
\begin{tabular}{r r l | r r l }
4/3 1927 & 28/2 1987 & Fire--Rabbit & 2/3 1957 & 27/2 2017 & Fire--Bird\\ 
22/2 1928 & 18/2 1988 & Earth--Dragon & 19/2 1958 & 16/2 2018 & Earth--Dog\\ 
10/2 1929 & 7/2 1989 & Earth--Snake & 8/2 1959 & 5/2 2019 & Earth--Pig\\ 
1/3 1930 & 26/2 1990 & Iron--Horse & 27/2 1960 & 24/2 2020 & Iron--Mouse\\ 
18/2 1931 & 15/2 1991 & Iron--Sheep & 16/2 1961 & 12/2 2021 & Iron--Ox\\ 
7/2 1932 & 5/3 1992 & Water--Monkey & 5/2 1962 & 3/3 2022 & Water--Tiger\\ 
25/2 1933 & 22/2 1993 & Water--Bird & 24/2 1963 & 21/2 2023 & Water--Rabbit\\ 
14/2 1934 & 11/2 1994 & Wood--Dog & 14/2 1964 & 10/2 2024 & Wood--Dragon\\ 
4/2 1935 & 2/3 1995 & Wood--Pig & 4/3 1965 & 28/2 2025 & Wood--Snake\\ 
23/2 1936 & 19/2 1996 & Fire--Mouse & 21/2 1966 & 18/2 2026 & Fire--Horse\\ 
12/2 1937 & 8/2 1997 & Fire--Ox & 10/2 1967 & 7/2 2027 & Fire--Sheep\\ 
3/3 1938 & 27/2 1998 & Earth--Tiger & 29/2 1968 & 26/2 2028 & Earth--Monkey\\ 
20/2 1939 & 17/2 1999 & Earth--Rabbit & 17/2 1969 & 14/2 2029 & Earth--Bird\\ 
9/2 1940 & 6/2 2000 & Iron--Dragon & 7/2 1970 & 5/3 2030 & Iron--Dog\\ 
26/2 1941 & 24/2 2001 & Iron--Snake & 26/2 1971 & 22/2 2031 & Iron--Pig\\ 
16/2 1942 & 13/2 2002 & Water--Horse & 15/2 1972 & 12/2 2032 & Water--Mouse\\ 
5/2 1943 & 3/3 2003 & Water--Sheep & 5/3 1973 & 2/3 2033 & Water--Ox\\ 
24/2 1944 & 21/2 2004 & Wood--Monkey & 22/2 1974 & 19/2 2034 & Wood--Tiger\\ 
13/2 1945 & 9/2 2005 & Wood--Bird & 11/2 1975 & 9/2 2035 & Wood--Rabbit\\ 
4/3 1946 & 28/2 2006 & Fire--Dog & 1/3 1976 & 27/2 2036 & Fire--Dragon\\ 
21/2 1947 & 18/2 2007 & Fire--Pig & 19/2 1977 & 15/2 2037 & Fire--Snake\\ 
10/2 1948 & 7/2 2008 & Earth--Mouse & 8/2 1978 & 6/3 2038 & Earth--Horse\\ 
28/2 1949 & 25/2 2009 & Earth--Ox & 27/2 1979 & 23/2 2039 & Earth--Sheep\\ 
17/2 1950 & 14/2 2010 & Iron--Tiger & 17/2 1980 & 13/2 2040 & Iron--Monkey\\ 
7/2 1951 & 5/3 2011 & Iron--Rabbit & 5/2 1981 & 3/3 2041 & Iron--Bird\\ 
26/2 1952 & 22/2 2012 & Water--Dragon & 24/2 1982 & 21/2 2042 & Water--Dog\\ 
14/2 1953 & 11/2 2013 & Water--Snake & 13/2 1983 & 10/2 2043 & Water--Pig\\ 
4/2 1954 & 2/3 2014 & Wood--Horse & 3/3 1984 & 29/2 2044 & Wood--Mouse\\ 
23/2 1955 & 19/2 2015 & Wood--Sheep & 20/2 1985 & 17/2 2045 & Wood--Ox\\ 
12/2 1956 & 9/2 2016 & Fire--Monkey & 9/2 1986 & 6/2 2046 & Fire--Tiger\\ 
\end{tabular}
\caption{Tibetan New Year, \tibx{Losar}, 
(\PH{} version) and year names for the last and
  current 60 year cycles.}
\label{Tlosar}  
\end{table}

\refT{Tlosar} gives the Gregorian dates of the Tibetan New Year  
\tib{Losar}
for 
the years in the last and current 60 year cycles (see \refS{Syear}).
See also \refT{T4Losar} in \refApp{Aversions}.
The table shows that New Year at present
occurs in February or the first week of March; the extreme dates during the
20th century are 4 February (\eg{} 1954) and 5 March (\eg{} 1992), and
during the 21st century 5 February (2019) and 7 March (2095). (The dates get
slowly later in the Gregorian calendar, see \refS{Smean}.) 
\xfootnote{
The Tibetan New Year thus frequently coincides with the Chinese New Year, 
which also is at new moon and always in the range 21 January -- 21 February
\cite[Section 17.6]{CC}, although often the New Years differ by a day
because of the different calculation methods; however, in about each second
year (at present), the Tibetan New Year is a month after the Chinese.
}

The purpose of this paper is to describe the mathematics of the Tibetan
calendar, using modern 
mathematical notations instead of the traditional methods.
After some preliminaries (Sections \ref{Snot}--\ref{Sdef})
and a description of the naming of years (\refS{Syear} and \refApp{A60}),
the calculations of  the Tibetan calendar are presented in two main
parts. In the first part (\refS{Smonths} and \refApp{Aleap}), 
the months are regarded as units and I 
discuss how they are numbered, which implies the partitioning of them
into years and also shows which months that are leap months.
In the second part (Sections \ref{Sdays}--\ref{Sweek}),
I discuss the coupling between months
and days, including finding the actual days when a month begins and
ends and the numbering of the days.
Finally, some further calculations are described 
(Sections \ref{Sfurther}--\ref{Sholidays}) and
some mathematical consequences are given (Sections
\ref{Smean}--\ref{Speriod}).
Calculations for the planets and some other astrological 
calculations are
described in Appendices \ref{Aplanet}--\ref{Aastro}.

The description is based mainly on the books by \citet{Schuh} and 
\citet{Henning}, but the
analysis and mathematical formulations are often my own.
(Unfortunately I do not read Tibetan, so I have to use secondary
sources instead of Tibetan texts. This is of course is a serious drawback,
although I have been able to check the calculations against published
almanacs such as \cite{tib2013}.)
For further study I recommend in particular the detailed recent book by
\citet{Henning}, 
which contains much more material than this paper; 
see also his web site \cite{kalacakra} with further information. 
I (mostly) describe only the contemporary versions and ignore  the historical
development; for the history of the calendar, see  \citet{Schuh} and
\cite{tibetenc}.
See also \citet{CC} for a related description 
and computer implementation of the Tibetan calendar.
Open source computer programs can be obtained from
\citet{kalacakra}.
Computer generated tables and calendars covering almost 1000 years are 
given in \citet{Schuh} and by
\citet[Traditional Tibetan calendar archive]{kalacakra}. 
Some further related references are
\citet[pp.~403--409]{Ginzel} (a classic, but dated and incomplete),
\citet{Petri} (on Tibetan astronomy),
\citet{Tseng1,Tseng2} (on astrology
\xfootnote{\label{f:astro}%
 \citet{Tseng1,Tseng2} and \citet{Schuh-review} make a distinction between
 astrology 
(using planets etc.)\
and divination (\tibx{nag-rtsis} or \tibx{'byung-rtsis}, using elements etc.).
I use ``astrology'' in a wider sense, encompassing both.
}%
),
and further articles in the collection 
\cite{Schuh-contributions}.
See also the review by \citet{Schuh-review} of the first (2007)
version of the present paper (together with the books \cite{CC} and
\cite{Henning}), and the reply by \citet{Henning-comments}.
\xfootnote{
The review by \citet{Schuh-review} is quite critical of many details.
I am grateful for some of his comments, and I have tried
to make some improvements in the present version. 
(In other cases I agree that improvements of formulations might be
desirable, but I find the present version acceptable and leave it for
various reasons.
I also agree that several of the internet references are unreliable and not
scholarly; however, I have nevertheless included them, either to illustrate 
actual usage among Buddhist groups or to show interesting claims that I have
been unable to verify or disprove; the reader is advised to regard these
with caution.) 
I disagree with some other comments, of which I will mention a few here.
First, the present paper is primarily written for people who, like me, do not
know Tibetan, and I therefore find it useful to use more or less standard
Anglicized 
forms of names rather than the scientific transliterations (Wylie) used by
Tibetologists (who hopefully will understand what I mean in any case),
although I sometimes give the latter as well (in italics).
(For example, I write 
``\TS'',  
used \eg{} by the \TS{} monastery itself on its English web pages 
and on the almanacs from Rumtek shown in 
\cite[Open source Tsurphu calendar software]{kalacakra}.)
Similarly, I have chosen to follow \cite{Henning} and call the main version of
the Tibetan 
calendar the ``Phugpa version'', for convenience and without  
going into the detailed historical background (see \cite{Schuh}).
Furthermore, I usually use the English terminology of \cite{Henning} for
Tibetan terms. Schuh does not like these translations; I cannot argue with
his expert linguistic remarks, but I am not convinced that
they are relevant for a description of the mathematical content rather than
a linguistic and cultural study of the calendar tradition
(see also \cite{Henning-comments});
moreover, 
it seems better to use an existing English terminology \cite{Henning} rather
than making my own second translation from a German 
translation \cite{Schuh} of the terms.
}

\begin{acks}
This study was made possible by the help of
Nachum Dershowitz,   Edward Henning, Edward Reingold,
Olli Salmi
and Dieter Schuh,
who have provided me with essential references and given me helpful
comments;
in particular, Edward Henning has been very helpful by answering numerous
questions on many details.
\end{acks}

\section{Notation}\label{Snot}

\subsection*{Mixed radix notation}
Traditional Tibetan calculations are made 
expressing the various quantities as sequences (written as columns) of
integers, which 
should be interpreted as rational numbers 
in a positional system with mixed radices;
\xfootnote{
In the same way as we denote time in days, hours, minutes and
  seconds, with radices 24, 60, 60.
}
the radices are fixed but 
different sequences of radices are used for different quantities.
I usually use standard notation for 
rational numbers, but when quoting the traditional
expressions, I use notations of the type (of varying length)
\xfootnote{
This notation is taken from \citet{Henning}, although he usually omits all
or most of the radices since they are given by the context. \citet{Schuh},
\cite{tibetenc}
uses the similar notation 
$[a_1,a_2,\dots,a_n]/(b_1,b_2,\dots,b_n)$ meaning either
$0;a_1,a_2,\dots,a_n \rr{b_1,b_2,\dots,b_n}$
or
$a_1;a_2,\dots,a_n \rr{b_2,\dots,b_n}$ depending on the context.
}
\begin{equation*}
  a_0;a_1,a_2,a_3 \rr{b_1,b_2,b_3}
\qquad\text{meaning}\qquad
a_0+\frac{a_1+(a_2+(a_3/b_3))/b_2}{b_1}
.
\end{equation*}
Formally, we have the inductive definition
\begin{equation*}
  a_0;a_1,a_2,\dots,a_n \rr{b_1,b_2,\dots,b_n}
=a_0+\frac{ a_1;a_2,\dots,a_n \rr{b_2,\dots,b_n}}{b_1}
\end{equation*}
for $n\ge2$, and $a_0;a_1\rr{b_1}=a_0+a_1/b_1$. 
We will usually omit a leading 0 (and the semicolon) and write
just \eg{}  $a_1,a_2,a_3 \rr{b_1,b_2,b_3}$ for numbers between 0 and 1.

For explanations and examples of the  way the traditional hand calculations
are performed 
\xfootnote{Traditionally in sand rather than on paper.}, 
see \citet{Henning}, \citet{Schuh} and
\cite[Kalenderrechnung, Sandabakus]{tibetenc}.

\subsection*{Angular units}
It will be convenient, although somewhat unconventional, to express
longitudes and other angular measurements in units of full circles
in our formulas.
To obtain values in degrees, the numbers should thus be multiplied by
360. (A Tibetan would probably prefer multiplying by 27 to obtain
lunar mansions (= lunar station = \tibx{naksatra}) and fractions thereof;
this is the unit usually used for longitudes.  
A Western astrologer might prefer multiplying by 12 to obtain values
in signs. A mathematician might prefer multiplying by $2\pi$ (radians).)

For angular measurements, full circles are often to be ignored (but
see \refApp{Aleap}); this means with our convention that the numbers are
considered modulo 1, \ie{}, that only the fractional part matters.

\subsection*{Boolean variables}

For a Boolean variable $\ell$, \ie{} a variable taking one of the
values \true{} and \false, we use $\ell=\Bool{\cP}$ to denote that
$\ell=\true$ if and only if the condition $\cP$ holds; we further 
let $\boolx{\ell}$  be the number defined by
$\boolx{\ell}=1$  when $\ell=\true$ and 
$\boolx{\ell}=0$  when $\ell=\false$.

\subsection*{Julian day number}
The \emph{Julian day number} (which we abbreviate by JD) for a given day is
the number of days that have elapsed since the epoch 1 January 4713 BC
(Julian); for days before the epoch (which hardly concern the Tibetan
calendar), negative numbers are used. The Julian day numbers thus form
a continuous numbering of all days by $\dots,-1,0,1,2,\dots$. Such a
numbering is very convenient for many purposes, including conversions
between calendars. The choice of epoch for the day numbers is
arbitrary and for most purposes unimportant. 
The conventional date 1 January 4713 BC ($-4712$ with astronomical numbering
of years) was 
originally chosen by Scalinger in 1583 as the origin of the Julian
period, a (cyclic) numbering of years; 
this was  developed by 19th-century astronomers into a
numbering of days. See further 
\cite[Section 12.7]{AA}, \cite[Section 15.1.10]{AA3}.
\xfootnote{
Dershowitz and Reingold \cite{CC} use another day number, denoted by RD,
with another epoch: RD 1 is 1 January 1 (Gregorian) which is JD 1721426.
Consequently, the two day numbers
are related by $\JD=\RD+1721425$.
}

A closely related version of this idea is the \emph{Julian date},
which is a real number that defines the time of a particular instant,
meaured (in days and fractions of days) from the same epoch.
The fractional part of the Julian date shows the fraction of a day
that has elapsed since noon GMT (UT); thus, if $n$ is an integer, then
the Julian date is $n.0$ (\ie, exactly $n$) at noon GMT  on the
day with Julian day number $n$.

It is important to distinguish between the Julian
day number and the Julian date, even if they are closely related.
Both are extremely useful, but for different purposes, and much unnecessary
confusion has been caused by confusing and mixing them.
(We follow \cite{AA} in using  different names for them, but
that is not always done by other authors.
\xfootnote{
Moreover, 
astronomers simply define Julian day number as the integer part
 of the Julian date \cite[Section 3.7]{AA3}, 
which is  \emph{not} the version 
described here 
\cite[Section 12.7]{AA}, \cite[Section 15.1.10]{AA3},
used in historical chronology.
})
For this study, and most other work on calendars, the Julian day
number is the important concept.
(The Julian date is essential for exact astronomical calculations, but 
no such calculations are used in the traditional Tibetan calendar.)
Note that the Julian day number is an integer, while the Julian date
is a real number. (A computer scientist would say that they have
different types.) Moreover, the Julian day number numbers the days
regardless of when they begin and end, while the Julian date depends
on the time of day, at Greenwich.
Hence, to convert a Julian date to a Julian day number, we need in practice
to know both the local time the day begins and the time zone, while these
are irrelevant for calculations with the Julian day number.
For example, 1 January 2007 has JD 2454102, everywhere. Thus the Julian date
2454102.0 is 1 January 2007, noon GMT (UT), and the new year began at
Julian date 2454101.5 in Britain, but at other Julian dates in other
time zones. The Tibetan day begins at dawn, about 5 am local mean solar time
(see \refR{Rdate} below),
but we do not have to find the Julian date of that instant.

\subsection*{Other notations}
We let $\mod$ denote the binary operation defined by $m \mod n =x$ if
$x\equiv m \pmod n$ and $0\le x<n$ (we only consider  $n>0$, but $m$
may be of any sign; care has to be taken with this in computer
implementations). 

Similarly (following \cite{CC}), we let 
$\amod$ denote the binary operation defined by 
$m\amod n =x$ if $x\equiv m \pmod n$ and $0< x\le n$. 
This means that $m\amod n= m\mod n$ except when $m$ is a multiple of
$n$; then $m\mod n=0$ and $m\amod n=n$.
For integers $m$ and $n$ (the usual case, and the only one used in this paper),
$m\amod n=1+(m-1)\mod n$.

We use the standard notations $\floor{x}$ and $\ceil{x}$ for the
integers satisfying $x-1<\floor{x}\le x$  and 
$x\le\ceil{x}< x+1$, \ie{} $x$ rounded down and up to an integer, 
and $\FRAC(x)=x-\floor{x}=x\mod 1$ for the fractional part.
(Again, care has to be taken for $x<0$ in computer implementations.)

\section{Some concepts}\label{Sdef}
The calendar is based on considering several different types of days, months
and years, and we give a list of the most important ones here.
The calendar is in principle astronomical and based on the positions of the moon
and sun in the sky, more precisely their \emph{longitudes} 
\xfootnote{\label{fequinox}%
The longitude of a celestial object
is measured eastward along the ecliptic, with 0 at the first point of Aries, 
both in Tibetan and Western astronomy.
(This is the vernal equinox, where the sun crosses the celestial equator
northwards every year; however, in Tibetan astronomy, the vernal equinox is
erroneously believed to be much earlier, see 
\refF{fPHpoints} and \cite[p.~328]{Henning}.)
The sun is always on the ecliptic (with small
perturbations); the moon is up to $5.15\grad$ away from the ecliptic due to
the inclination of its orbit, but 
only its longitude is relevant for the calendar.
\cite[\S1.3, Table 15.4, p.~731]{AA}
}
and in particular
the \emph{elongation} of the moon, \ie, the difference between the
longitudes of the moon and the sun.
The Tibetan calendar uses two different formulas to calculate each of these:
a simple (linear) formula giving the \emph{mean longitude} (assuming uniform
motions of the sun and moon) and a more complicated formula for the
\emph{true longitude}; thus the ``lunar'' and ``solar'' concepts defined below 
have two versions,
one ``mean'' and one ``true''.
Note that the calendar always uses these theoretical values and not the
actual astronomical positions; do not confuse the ``true'' longitude with
the exact astronomical one. 
(Nowadays, the ``true longitudes'' actually have large errors, see
\refS{Smean}.) 

\begin{description}

\item [calendar day \rm(\tibx{gza'} or \tibx{nyin-zhag})]  
Solar day or natural day; in Tibet the calendar day is from
dawn to dawn. The length is thus constant, 24 hours.
The calendar days are numbered by the number of the corresponding lunar day,
but since the correspondence is not perfect, sometimes a number is skipped
or repeated, see \refS{Sdays}.
Each calendar day has also a day of week, in the same way as in Western
calendars, see \refS{Sweek}.

\item [lunar day \rm(\tibx{tshes-zhag, \textup{Sanskrit} tithi})]
$1/30$ of a lunar month; more precisely the time during which the
elongation of the moon increases by $1/30$ ($=12\grad$).
(Since the moon does not travel at uniform speed; the lunar days have
different lengths, varying between about 21.5 and 25.7 hours, see
\refR{Rlunarday}.)  
The lunar days are numbered by 1--30 in each lunar month, 
with day 1 beginning at new moon. Thus lunar day $i$ is when the elongation
is between $(i-1)/30$ and $i/30$.

\item[calendar month \rm\tib{zla-ba}] 
A period of 29 or 30 calendar days, approximating a lunar month. 
(The calendar month gets the same number as the lunar month, and they are
often regarded as the same, but strictly speaking  
the calendar month begins at  the beginning of a calendar day, at dawn,
while the lunar month begins at the instant of new moon.)

\item[lunar month \rm\tib{tshes-zla}] 
The period from the instant of a new moon to the next new moon.
(Synodic month in astronomic terminology.)
The lunar months are numbered 1--12, but sometimes there is a leap month and
a number is repeated, see \refS{Smonths}.

\item[solar month \rm\tib{khyim-zla}] 
The period during which the sun travels through one sign.
(Each sign is 1/12 of the ecliptic, \ie, $30\grad$.)
This is thus 1/12 of a solar year, although the lengths of the true solar
months vary somewhat.

\item[calendar year]
12 or 13 calendar months, according to the rules for inserting leap months.
The length of the calendar year is 354, 355, 383, 384 or 385 days,
see \refR{Rlunarday}. The average length is close to the length of the solar
year, see \refS{Smean}; 
in principle, the calendar year is on the average fixed
with respect to the seasons, although this is in reality not exact.

\item[solar year]
The period during which the sun travels a full revolution around the ecliptic.
(Tropical year in astronomical terminology.)

\end{description}

\section{Numbering and naming of years}\label{Syear}

Several ways of numbering or naming Tibetan years are used,
see \eg{} \citet[pp.~142--145]{Schuh}, \cite[Kalender]{tibetenc}.
One common method, especially among Westeners, is to simply
number a Tibetan year by the Gregorian (or
Julian) year in which it starts (and where most of it falls).
For convenience, we will use this numbering below.

Another method is to number the years from an epoch 127 BC
(the traditional ascent of the first Tibetan king); the 
Tibetan year starting in Gregorian year $Y$ will then be numbered
$Y+127$. 
\xfootnote{This method has been used from the second half of the 20th  century.
Epochs 255 and 1027 have also been used. See \cite[Kalender]{tibetenc}.
}
Both methods are used by Tibetans; 
for example, the Tibetan calendar
\cite{tib2003} has titles in both Tibetan and English; with 2003 in
the English title and 2130 in the Tibetan.

Moreover, and more importantly, each year is named according
to a 60 year cycle. Actually, there are two different 60 year cycles
of names, one Indian and one Chinese. Of course, since the cycles have
the same length, there is a 1--1 correspondence between the names.
When naming years, the Chinese cycle is almost always used, and
sometimes the Indian and Chinese cycle names are used together.
(For example, the Chinese cycle names are used in the titles of the calendars
\cite{tib2003} and \cite{nitarthacal}.)

The Indian cycle is a list of 60 different names, in Sanskrit or
Tibetan, see \refApp{A60}. The cycle is named after its first year as
Prabhava \tib{rab byung}. The cycles are numbered, with the first
cycle beginning in AD 1027, which means that each year can be
unambiguously identified by its name in the cycle and the number of
the cycle; this method of naming years has sometimes been used
\cite[Sechzigjahreszyklen]{tibetenc}.
(See \refF{fmongol} and \cite{Terbish-GenghisKhan}
for  modern Mongolian examples.)

Year $n$ in cycle $m$ (with $1\le n\le 60$ and, presumably, $m\ge1$) thus
corresponds to Gregorian or Julian year $Y$ given by
\begin{equation}\label{indian1}
  Y = 1027 + (m-1)60+(n-1)
= 60m+n+966.
\end{equation}
Conversely,
\begin{align}
  n&=(Y-1026)\amod60
=(Y-6)\amod60,
\label{indian60}
\\
  m&=\Ceil{\frac{Y-1026}{60}}. \label{indian3}
\end{align}
For example, AD 2007 is the 21st year in the 17th Prabhava cycle,
which began in 1987.

The Chinese cycle is identical to the one used in the Chinese
calendar \cite{CC}. 
The cycles start 3 years before the Indian ones, so the
first year (Prabhava, \tibx{rab byung}) in the Indian cycle is the
fourth year  (Fire--female--Rabbit) in the Chinese cycle.
The full correspondence is given in \refApp{A60}.
The last cycle thus started 1984.
Hence, or by \eqref{indian60}, year $Y$ is number 
\begin{equation}\label{chinese60}
(Y-3)\amod60  
\end{equation}
in the Chinese cycle.
The Chinese cycles are not numbered in the Tibetan calendar.

The Chinese 60 year cycle is a combination of two cycles, of 10 and 12
years respectively. (Note that the least common multiple of 10 and 12
is 60; since 10 and 12 have greatest common divisor 2, only half of
the $10\times12$ combinations are possible.)
\xfootnote{
These cycles have been used in the Chinese calendar since at least 1400 BC,
first for naming days and from the Zhou dynasty (c.~1000 BC) for naming years.
\cite[p.~163]{Richards}
}

The 10 year cycle consists in China of 10 different names (proper
names with no other English translation) called celestial stems.
Each celestial stem is associated with an element (wood, fire, earth,
iron, water) and a gender (female or male), see \refT{T10}.
(In Chinese, the two genders are the well-known \emph{yin} and
\emph{yang}.)
Note that the 2 genders are alternating and thus are given by the year
mod 2, while the 5 elements are repeated 2 years each in the 10
year cycle. As a consequence, each celestial stem corresponds to a
unique (element, gender) pair, and in the Tibetan calendar, the
element and gender are used to name the year; the Chinese celestial stems are
usually not used \cite[p.~145]{Henning}.

It follows from \eqref{chinese60} that,
using the numberings in Tables \ref{T10} and \ref{T5}, year $Y$ is
$z=(Y-3)\amod 10$ in the Chinese 10 year cycle, and has element $\ceil{z/2}$.

The 12 year cycle is the well-known cycle of animals found in the Chinese
and many other Asian calendars, see \refT{T12}.
(The English translations are taken from \cite{Henning}; several easily
recognized variants exist.
The Tibetan names are from \cite{CC}.) 

It follows from \eqref{chinese60} that
year $Y$ is
$z=(Y-3)\amod 12$ in the 12 year animal cycle.

The Tibetan name for a year according to the Chinese cycle is thus
given as Element--gender--Animal. 
Note that the gender, being the year mod 2, also is determined by
the animal (since 2 divides 12), as shown in \refT{T12}. Indeed, the
gender is often omitted and only Element--Animal is used as the name
of the year.
For example, AD 2007 is the 24th year in the Chinese cycle (21st in
the Indian) and is thus Fire--female--Pig or simply Fire--Pig.

\begin{table}[htpb]
\begin{tabular}{r r l l l l}
\rlap{year}\phantom{0}& & element & gender &Tibetan
& celestial stem (Chinese) \\
\hline 
1 &8& wood & male &\tibx{shing-pho} & ji\v{a}  \\
2 &9& wood & female &\tibx{shing-mo} & y\v{\i} \\
3 &10& fire & male &\tibx{me-pho} & b\v{\i}ng  \\
4 &1& fire & female &\tibx{me-mo} & d\={\i}ng  \\
5 &2& earth & male &\tibx{sa-pho} & w\`u  \\
6 &3& earth & female &\tibx{sa-mo} & j\v{\i}  \\
7 &4& iron & male &\tibx{lcags-pho} & g\={e}ng  \\
8 &5& iron & female &\tibx{lcags-mo} & x\={\i}n  \\
9 &6& water & male &\tibx{chu-pho} & r\'{e}n  \\
10 &7& water & female &\tibx{chu-mo} & gu\v{\i} \\
\end{tabular}
\caption{The 10 year cycle. 
The first number on each line shows the year mod 10 counted from the
  start of a Chinese cycle; the second shows the year mod 10 counted from the
start of a Prabhava cycle.}
\label{T10}  
\end{table}

\begin{table}[htpb]
\begin{tabular}{r r l l l}
\rlap{year}\phantom{0} && animal & gender & Tibetan\\
\hline
1 &10&  Mouse  & male &\tibx{byi ba}\\ 
2 &11&  Ox  & female &\tibx{glang}\\ 
3 &12&  Tiger & male  &\tibx{stag} \\
4 &1&  Rabbit & female & \tibx{yos}\\ 
5 &2&  Dragon & male & \tibx{'brug}\\
6 &3& Snake & female &\tibx{sbrul} \\
7 &4& Horse & male &\tibx{rta} \\
8 &5& Sheep & female &\tibx{lug} \\
9 &6& Monkey & male  &\tibx{spre'u} \\
10 &7& Bird  & female &\tibx{bya}\\ 
11 &8& Dog & male &\tibx{khyi}\\
12 &9& Pig & female &\tibx{phag} \\
\end{tabular}  
\caption{The 12 year cycle. 
The first number on each line shows the year mod 12 counted from the
  start of a Chinese cycle; the second shows the year mod 12 counted from the
start of a Prabhava cycle.}
\label{T12}
\end{table}







The civil year starts, unsurprisingly, with month 1; note that in 
case there is a leap month 1, the year begins with the leap month,
which precedes the regular month 1. 
(This happened in 2000. See also \refR{Rlosar}.)
The first day of the year is also celebrated as a major holiday 
(\tibx{Losar}; the Tibetan New Year).

\begin{remark}\label{Rcaitra}
Traditional Tibetan almanacs (such as \cite{tib2013})
start, however, with month 3; this month is
identified with the \Kc{} month 
\tibx{nag pa} (Tibetan)
or \Caitra{} (Sanskrit),
which was
considered the beginning of the \Kc{} year \cite[p.~194]{Henning}
and is the starting point in Tibetan astronomy.
(\Caitra{} 
is the first month of the year 
in most Indian calendars \cite{CC}.)
The standard numbering system is known as 
\emph{Mongolian month} \tib{hor-zla}
(or, in recent almanacs, \emph{Tibetan month} \tib{bod-zla}), 
and was introduced in the 13th century when the Tibetan
New Year was moved to be the same as the Mongolian 
\xfootnote{
By the Tibetan spiritual and political leader 
Chogyel Phagspa Lodro Gyeltsen 
\tib{chos-rgyal 'Phags-pa Blo-gros rgyal-mtshan}
(1235--1280), 
who was advisor of Kublai Khan and  
Preceptor of Tibet, then part of the Mongol Empire.
\cite[Pakpa Lodro Gyeltsen]{treasury}
}
\cite[p.~145]{Henning},
\cite[Kalendar,  Zeitma{\ss}e]{tibetenc}.

The year is nevertheless considered to consist of months 1--12 in the usual
way, 
as said above,
so a traditional almanac contains the 10 last
months of the year and 2 months of the next (perhaps plus a leap month).
\end{remark}

\section{Numbering of months and leap months}\label{Smonths}
The Tibetan calendar months are, as said in \refS{Sdef}, lunar months, that 
begin and end at new moon.
There are usually 12 months in a year, 
but there is also
sometimes an extra \emph{leap} (or \emph{intercalary}) month, so that the
year then has 13 months (a \emph{leap year}), 
in order to keep the calendar year roughly aligned with the (tropical) solar
year. 
Non-leap months are called \emph{regular}. There are thus always 12 regular
months in a year. 

The regular months in each year are numbered 
1--12. 
A leap month takes the same number (and belongs to the same year) 
as the following month (as in Indian calendars \cite{CC}). 
\xfootnote{
The original system in the \KT{} seems to have been the opposite
(although this is not stated explicitly),
with a leap month taking the number of the preceding month
(as in the Chinese calendar \cite{CC}),
see \cite[Published calendar explanation]{kalacakra}.
That system is also used in the Bhutanese version of the
calendar, see  \refS{ABhutan}.
}
In a leap year there are thus two months having the same number and, by the
rule just given, the first of these is the leap month.

\subsection{Month names}\label{SSmonths}
Usually, the Tibetan months are just numbered, but month names exist too;
they are printed in almanacs 
and are sometimes used, in particular in literature about the calendar.
Several different naming systems exist as follows, see \refT{Tmonths}.
(In all systems, a leap month gets the same name as the corresponding regular
month.) 
See further 
\citet[p.~145]{Schuh}, 
\cite{Schuh-review} and \cite[Zeitma{\ss}e]{tibetenc},
and \citet[pp.~147--149 and 194--196]{Henning} and
\cite[Early epochs]{kalacakra}.

\begin{romenumerate}[-15pt]
\item \label{month-mansion}
The Indian system of naming months by 
twelve of the names of the 
lunar mansions, chosen to be approximatively where the full moon occurs that
month (and thus opposite to the sign that the sun visits).
(Cf.\ \cite[Section 9.3]{CC}. See also \citet[pp.~358--359]{Henning} and
\citet[p.~92]{Petri}.) 
See \refT{Tmonths}, where both the Tibetan names and the  Sanskrit
names are given.
This system was introduced with the \Kc{} calendar in the 11th
century. 

\item \label{month-animal}
Animal names, by the same 12 animal cycle as for years, see \refT{T12}.
This can be extended to Element + Animal, or Element + Gender + Animal,
as for years, see \refApp{ASattributes-months} for details.
(This is an old Chinese system, used in Tibet from the 12th century.) 
Note that \PH{} almanacs set month 11 = Tiger (shown in \refT{Tmonths}), 
while  \TS{} almanacs (\refApp{ATS}) set  month 1 = Tiger (as in
the Chinese calendar), see \refT{T12months}.
This  method is  still important for astrological purposes.

\item 
Seasonal names, naming the 12 months
as beginning, middle and end of each of the four
seasons spring, summer, autumn, winter.
This is the oldest system, and was used at least from the 7th century.
\cite[Zeitma{\ss}e]{tibetenc}

Unfortunately, the seasonal months have been identified with the \Kc{}
months in several different ways, so there are several different versions.
As far as I know, none of them is really used today, but
they are mentioned \eg{} in almanacs, which often give two or several different
versions of them, see \eg{} \cite[pp.~195--196]{Henning}.
At least the following versions exist today.
(Including versions used by the \TS{} tradition, see \refApp{ATS}.)
\begin{enumerate}
\item \label{month-season-KT}
Month 1 = \tibx{mchu} = late-winter, etc.;
thus Month 3 = \nag{} = mid-spring.
This is the original \Kc{} identification, now used by the \TS{} tradition.

\item \label{month-season-PH}
Month 1 = \tibx{mchu} = early-spring, etc.;
thus Month 3 = \nag{} = late-spring.
See \refT{Tmonths}.
This is the \PH{} version of the \Kc{} system. 
(Introduced c.~1200 by Drakpa Gyaltsen, see \cite[Early epochs]{kalacakra}.)

\item \label{month-season-animal}
Identifying the animal names in \ref{month-animal} by Tiger = early-spring,
etc.\ (as in the Chinese calendar).
Note that this gives different seasonal names for the \PH{} and \TS{}
versions.
The \PH{} version is given in \refT{Tmonths}.

\item
A different system with 6 seasons with two months each.
(This is an Indian system, see \eg{} 
\cite[\S 75 and pp.~345--346]{Ginzel} and \cite[p.~195]{Richards}.)
Month 1 = \tibx{mchu} = early-winter;
Month 3 = \nag{} = early-spring;
see \cite[pp.~358--359]{Henning} for the complete list.
\end{enumerate}
 
\end{romenumerate}

The present numbering of months was, as said above, introduced in the 13th
century. 
\xfootnote{In the \KT, \nag{} (\Caitra) is evidently the first month, as in
  Indian  calendars, although the months are not explicitly numbered.
}

\begin{table}[htpb]
\begin{tabular}{r l l l l l}
 & Tibetan & Sanskrit & seasonal I & animal &seasonal II\\
\hline
1 & mchu & M\=agha & early-spring & dragon & late-spring\\
2 & dbo & Ph\=alguna & mid-spring & snake & early-summer\\
3 &  \tibxx{nag pa} & \tibxx{Caitra} & late-spring & horse & mid-summer\\
4 &  sa ga  & Vai\'s\=akha & early-summer & sheep & late-summer\\
5 & snron & Jye{\d s}{\d t}ha & mid-summer & monkey &early-autumn\\
6 & chu stod & \=A\d{s}\=a\d{d}ha & late-summer & bird & mid-autumn\\
7 & gro bzhin & \'Sr\=ava\d{n}a & early-autumn & dog & late-autumn\\
8 & khrums & Bh\=adrapada & mid-autumn & pig & early-winter\\
9 & tha skar & \=A\'svina & late-autumn & mouse & mid-winter \\
10 & smin drug & K\=artikka & early-winter & ox & late-winter\\
11 & mgo & M\=arga\'s\=ir\d{s}a & mid-winter & tiger& early-spring \\
12 & rgyal & Pau\d{s}a & late-winter & rabbit &mid-spring
\end{tabular}  
\caption{Month names  (\PH{} system), see \refS{SSmonths}: 
number; 
Tibetan and Sanskrit lunar mansion names,  \ref{month-mansion};
seasonal names, \ref{month-season-PH};
animal names, \ref{month-animal};
seasonal names, \ref{month-season-animal}.
}
\label{Tmonths}
\end{table}

\subsection{Leap months}
The Tibetan calendar is based on the relation
\begin{equation}\label{6765}
67 \text{ lunar months} = 65 \text{ solar months},
\end{equation}
which is regarded as exact. 
(See \refS{SSratio} for the astronomical reality.)
This is a fundamental relation in the Tibetan calendar,
which distinguishes it from other calendars such as the Indian ones.
\xfootnote{
This relation derives, as all basic features of the calendar, from the \KT, 
but there it is actually given as an approximation and not as an exact
relation, see \refApp{AKT}.
}
The leap months are regularly spaced in accordance with this relation, \ie,
with 2 leap months for each 65 solar months. In other words, there are 2
leap months for 65 regular months, and these are regularly spaced with
alternating 32 and 33 regular months between the leap months. 
(Thus the distance between successive leap months alternates between 33 and
34 months.)
\xfootnote{The average length between intercalations is thus $32\frac12$ regular
  months, or $33\frac12$ months including the leap month. Tibetan authors have
  often interpreted this as meaning that each second intercalation really 
should come in the middle of a month, although it is moved to the beginning of
  the month, see \cite{Yamaguchi}, \cite[On intercalary months]{kalacakra}.
As far as I know, this view has only been used in theoretical discussions and
no Tibetan calendar has ever been produced with the leap month in the middle
of a month. 
(However, leap months are inserted in this way in some Indian calendars
\cite[p.~277]{CC}.)
}

\begin{remark}\label{R6567}
It follows that the leap months repeat in a cycle of 65 years.
 In 65 years, there are $65\cdot12=780$ regular months and
therefore 
$2\cdot 12=24$ leap months, for a total of $804$ months.
Since 12 and 65 are relatively prime, 
the leap month can occur at any place in the year; in 65 years, each
leap month 1--12 occurs exactly twice. 
\end{remark}

We describe here the basic algorithmic calculations to determine leap months
and thus the numbering of all months. The astronomical theory behind the
rules is explained in \refApp{Aleap}.

The months can be described by year, number (from 1 to 12) and,
possibly, the label ``leap''. We can thus formally think of each month as
labelled by a triple $(Y,M,\ell)$, where $Y\in\bbZ$,
$m\in\set{1,\dots,12}$ and $\ell\in\set{\true,\false}$ is a Boolean variable.
(As said in \refS{Syear}, we  number the Tibetan years by the Gregorian (or
Julian) year in which they start.)

The Tibetan calendar calculations also use a 
consecutive numbering of all months
(regular or leap)
starting with 0  at some epoch. 
(This is thus a linear numbering, ignoring the division into years, unlike
the standard cyclic numbering that starts again with each new year.)
The epoch could in principle be the beginning of
an arbitrarily chosen month; we assume in our formulas that the
epoch is the beginning of  month
$M_0$ year $Y_0$.

Any epoch will give the same calendar (provided the correct initial data are
used), and in calculations only a single epoch is chosen.
Nevertheless, to illustrate the calendar mathematics, 
we give (and compare) data for three different epochs.  
The three epochs we use are 
month 3  year 806
(the traditional epoch from \KT{} although it is
several centuries before the calendar came to Tibet,
used  \eg{} by \citet{Schuh}),
month 3 year 1927
(used \eg{} by \citet{Henning})
and 
month 3 year
1987
(used \eg{} by \citet{LaiDolma} and the almanac 
\cite{tib2013}).
We denote these epochs by E806, E1927, E1987.
\xfootnote{
I have chosen these three epochs just as examples.
Another, historically important, traditional epoch is 1687
\cite[p.~331]{Henning}, \cite[Epoch data]{kalacakra}.
}

See further \refR{Repoch} below.

\begin{remark}\label{Repoch0}
By tradition, the epoch is always at the beginning of 
month 3 (\tibx{nag pa}, \cf{} \refR{Rcaitra}), 
so $M_0=3$. 
(But see \refR{Repoch} for the correct interpretation.)
The year $Y_0$ is by tradition usually (but not always) chosen to be 
the first year  \tib{rab byung} of a Prabhava cycle; it is common to use the
first year of the present cycle (so for example in the almanac
\cite{tib2013} where the epoch is 1987). 
(This is convenient for hand calculations 
because it gives smaller numbers than older epochs.)
\end{remark}

The number of a month in the linear numbering from a given epoch
is called the
\emph{true month} \tib{zla-dag}, and is calculated as follows
for month $M$ year $Y$:

First, solar months are counted starting after the epoch (month $M_0$, year
$Y_0$);  for month $M$ year $Y$, 
the ``number of solar months''  is 
\begin{equation}\label{MM}
\MM =12(Y-Y_0)+M-M_0.  
\end{equation}

Next, one calculates a preliminary version of
the {true month} 
as the rational number
\begin{align}\label{tm0}
&  \frac{67}{65}\MM+\frac{\gbx}{65}
=  \frac{67\bigpar{12(Y-Y_0)+(M-M_0)}+\gbx}{65}, 
\end{align}
where $\gbx$ is a constant depending on the epoch.
For our three example epochs we have
\begin{align}
  \label{gbxS}
\gbx&=61
\qquad(\eS);
\\
\gbx&=55
\qquad(\eH);
  \label{gbxH}
\\
\gbx&=0\phantom0
\qquad(\eX).
  \label{gbxX}
\end{align}
We write the
fractional part of the true month \eqref{tm0}
as $ix/65$, where the integer $ix$ is called 
the \emph{intercalation index}.
\xfootnote{
We use this name
from  \cite{Henning};
the Tibetan name \tibx{zla-bshol rtsis-'phro} simply means
``remainder for the calculation of leap-month'' 
\cite{Schuh-review}, \cite[Kalenderrechnung]{tibetenc}. 
}
Thus, the (preliminary version of the) intercalation index is
\begin{align}
  \label{ix}
ix&=(67 \MM+\gbx) \bmod{65}
=(2\MM+\gbx) \bmod{65},
\end{align}
with $\MM$ given by \eqref{MM}.
Note that $\gbx$ is the initial value of $ix$ at the epoch.

The traditional \PH{} leap month rule is:
\begin{equation}\label{leaprule-P}
\hskip-1em
\vbox {\narrower\narrower\narrower\noindent\em
A leap month is inserted when the intercalation index\\
$ix = 48\text{ or\/ }49$.
}
\hskip-3em
\end{equation}

For the rest of the calendar calculations, the true month is rounded to an
integer by the following rule.
(We follow \cite{Henning} and use the same name ``true month'' for both the
rational version \eqref{tm0}
and the rounded integer version,
but in order to avoid confusion, we will often call the latter
``\tmc''.)
\begin{equation}\label{mcrule-P}
\hskip-1em
\vbox {\narrower\narrower\narrower\noindent\em
The true month \eqref{tm0} is rounded down to the nearest integer if
$ix<48$ and rounded up if $ix\ge48$, except that for a leap month (when $ix=48$
or $49$), always round down. 
}
\hskip-3em
\raisetag\baselineskip
\end{equation}
When there is a leap month, there are two months with the same number $M$
the same year. 
The rule just given means that the true month is rounded down for the first
of these (the leap month) and rounded up for the second (the regular month);
thus the \tmc{} will increase by 1 for each new month also immediately
before and after the leap month. Similarly, it is easily checked that the
\tmc{}
increases by 1 also when the intercalation index \eqref{ix} passes
65 and drops to 0 or 1 again. Hence the \tmc{}
calculated by \eqref{mcrule-P} is really a
continuous count of months.

\begin{remark}\label{Rix}
The formulation in \eqref{mcrule-P} 
differs somewhat from the traditional formulation,
described \eg{} by \citet{Henning}.
Traditionally, the true month is calculated by \eqref{tm0} as above, but 
if this yields an 
intercalation index $>49$, one notes that there has been an earlier leap
month and therefore the numbering of months has to be corrected, so the
true month just calculated for month $M$ really applies to month $M-1$;
hence the true month for such a
month $M$ 
(and for a regular month immediately following a leap month)
is obtained by doing the calculation for
$M+1$ (which means adding 
$1;2\rr{1,65}=1\frac{2}{65}$ to the true month calculated for
$M$).
\xfootnote{\label{f6566}%
This leaves a gap when the intercalation index reaches 0 or 1:
If month $M$ has true month $n;0$ or $n;1\rr{1,65}$, then the
calculation for $M-1$ yields true month $(n-2);63$ or $(n-2);64$, which thus
applies to $M-2$. 
Hence no starting month yields a calculation that applies to $M-1$;
the missing month $M-1$ is usually given
true month $(n-1);0$ or $(n-1);1$, respectively. Note that the integer part
is correct, and that the intercalation index 0 or 1 is repeated.
(Sometimes 65 and 66 are used instead [Henning, personal communication].)
These complications do not appear in my version \eqref{mcrule-P}.
}
The \tmc{} then is always the integer part of the true month
calculated in this way;
\ie, the true month is always rounded down. 
This evidently yields the same \tmc{} as
\eqref{mcrule-P}.
However, the fractional part 
will differ by $\frac{2}{65}$ for these months,
\ie{} the intercalation index will differ by 2.
(The intercalation index is not used further for the main calendar
calculations, but
almanacs publish for each month the true month with its fractional part,
\ie, the intercalation index,
see \refS{Sfurther}, and it has at least one use, see \refS{Stent}.)

Hence, in order to get the correct intercalation index,
\eqref{mcrule-P} should be extended by the rule that when the true month is
rounded up, the intercalation index is increased by 2.
(When the intercalation index from \eqref{ix} is 63 or 64,
the new value can be given as either 65 or 66, or as 0 or 1, see \refF{f6566}.)
\end{remark}

Since
$48 =65-17$,
\eqref{mcrule-P} can be written, using \eqref{tm0} and \eqref{MM}, and denoting
the \tmc{} by $n$,
\begin{equation}\label{c2}
  n 
=\Floor{\frac{67\MM+\gbx+17}{65}}-\boolx{\ell}.
\end{equation}

\begin{remark}\label{Repoch}
As said in \refR{Repoch0} above,
the epoch is by tradition
always at the beginning of 
month 3 \tib{nag pa}.
However, because of the leap year rules, this has to be interpreted as
follows:

The calculation of the true month uses the number of solar months after
the epoch, calculated by \eqref{MM}.
For month 3 at the epoch, this is 0, so the true month is by \eqref{tm0}
$\gbx/65$; thus the
intercalation index is $\gbx$. According to the rounding rule \eqref{mcrule-P},
this month  thus has \tmc{} 0 if $\gbx<48$ but \tmc{} 1 if
$\gbx\ge48$.
(If $\gbx=48$ or 49, there is a leap month 3 and the two months 3 have
\tmc{s} 0 and 1.)
Hence, it is only when $\gbx\le49$ that the epoch really is at the beginning
of calendar month 3; otherwise it is at the beginning of
the preceding month (calendar month 2),
and the nominal epoch month 3 really has \tmc{} 1.
(This means that for an epoch with $\gbx<48$, the \tmc{} is the number of
elapsed (lunar) months
\emph{after} month 3 in the epoch year, while if $\gbx\ge48$, the \tmc{} is
the number of elapsed 
months \emph{starting with} month 3.)

For our three example epochs, the values of $\gbx$ in
\eqref{gbxS}--\eqref{gbxX}   show that
for \eS{} and \eH{}, the epoch month with \tmc{}  0 is month 2,
while for \eX, it is month 3. 
(See also the discussion in \cite{Schuh-review} and \cite{Henning-comments}
on the epoch of \eH.)

Note further that the epoch is not only a year and month; it is a specific
day, see 
\refR{Repoch2} in \refS{Sastro}.
\end{remark}

\subsection{The inverse problem}

Although not needed for the traditional construction of a yearly calendar, 
let us also consider the inverse problem:
to find the number $M$ and the year $Y$
of the month with \tmc{} $n$,
together with the indicator $\ell$ telling whether
the month is a leap month. 
We begin by writing \eqref{c2} as
\begin{equation}\label{c3}
  n 
=\frac{67\MM+\gbx+17-r}{65},
\end{equation} 
where $r$ is chosen such that the result is an integer;
for a regular month, $0 \le r \le 64$ 
($r$ is the remainder $(67\MM+\gbx+17)\bmod{65}$), 
and for a leap month, $r=65$ or $66$
(the remainder is 0 or 1, but $[\ell]=1$ in \eqref{c2}).
Rearranging, we get
\begin{equation}
  \MM = \frac{65n-\gbx-17+r}{67}
\end{equation}
with $0\le r \le 66$, and thus
\begin{equation}
  \MM = \Ceil{\frac{65n-\gbx-17}{67}}.
\end{equation}

Recalling \eqref{MM}, and assuming $M_0=3$, this can be written
\begin{equation}
  \label{jb1}
12(Y-Y_0)+M=\MM+3
=\Ceil{\frac{65n+184-\gbx}{67}}
=\Ceil{\frac{65n+\gb}{67}},
\end{equation}
where we define
\begin{equation}\label{gbgbx}
  \gb=184-\gbx.
\end{equation}
For our three example epochs, this yields,
by \eqref{gbxS}--\eqref{gbxX},
\begin{align}
  \label{gb1S0}
\gb&=123 
\qquad(\eS),
\\
\gb&=129 
\qquad(\eH),
  \label{gb1H0}
\\
\gb&=184 
\qquad(\eX).
  \label{gb1X0}
\end{align}

Consequently, we can calculate $(Y,M)$ from $n$ by
\begin{align}
  x&=\Ceil{\frac{65n+\gb}{67}},
\label{bx3}\\
M&= x \amod 12,
\label{bx4}\\
Y&=\frac{x-M}{12}+Y_0
= \Ceil{\frac{x}{12}}-1+Y_0.
\label{bx5}
\end{align}
To complete the calculations of $(Y,M,\ell)$ from $n$, we 
note that a month is leap if and only if it gets
the same number as the following one. Thus,
\begin{equation}
  \label{bx2}
\ell=\Bool{\Ceil{\frac{65n+\gb}{67}}=\Ceil{\frac{65(n+1)+\gb}{67}}}
=\bigbool{(65n+\gb)\bmod67 = 1\text{ or $2$}}.
\end{equation}
Note that $n=0$ in \eqref{bx3}--\eqref{bx4} 
yields $M=\ceil{\gb/67}$, so $\ceil{\gb/67}$ is
the number of the epoch month. (This is  2 or 3,
see \refR{Repoch} above.)

\subsection{A general rule}

By \eqref{ix}, \eqref{MM}, \eqref{leaprule-P} and $M_0=3$, 
there is a leap month $M$
year $Y$ if and only if
\begin{equation}\label{cx0}
  {24(Y-Y_0)+2M-6+\gbx} 
\equiv 48\text{ or }{49}\pmod{65}.
\end{equation}
Using \eqref{gbgbx}, $48+6-\gbx=\gb-130\equiv\gb\pmod{65}$ and thus  
\eqref{cx0} can be written
\begin{Rule}
There is a leap month $M$ in year $Y$ if and only if
\begin{equation}\label{leaprule-gb0}
24(Y-Y_0)+2M
\equiv \gb\text{ or\/ }\gb+1\pmod{65}.
\end{equation}
\end{Rule}
We will see later, in Appendices \ref{Aversions} and \ref{Aleap},
that the rule \eqref{leaprule-gb0}  holds
also for other versions of the Tibetan calendar, with appropriate $\gb$. 

 By combining \eqref{ix}, \eqref{MM} and $M_0=3$,
\begin{equation}\label{ixq}
  ix = \bigpar{24(Y-Y_0)+2M+\gbx-6}\bxmod{65}
\end{equation}
(if necessary modified as in \refR{Rix})
and hence the rule \eqref{leaprule-gb0} says 
that a
leap month has intercalation index
\begin{equation}\label{jeppe}
  (\gb+\gbx-6)\bmod65 
\qquad \text{or}\qquad
  (\gb+\gbx-5)\bmod65 ,
\end{equation}
which is the general form of the relation between the \PH{} formulas
\eqref{gbgbx} and \eqref{leaprule-P}.

\subsection{Leap years}\label{SSleapyear}%
We give some further formulas for
leap years and leap months, inspired by \cite{Salmi}, that follow from the
formulas above.
(These formulas are not traditional, and are not needed to construct
the calendar.)
We begin by stating them in a general form, valid also for other versions of
the calendar.

Since $M$ may be any number $1,\dots,12$, it follows 
from \eqref{leaprule-gb0}
that $Y$ is a leap year
if and only if
\begin{equation}\label{cu}
  24(Y-Y_0)-\gb \equiv 41,\,42,\dots,\text{ or }64\pmod{65},
\end{equation}
which also can be written as
\begin{equation}\label{cv}
\bigpar{ 24(Y-Y_0)-\gb} \bxmod 65 \ge 41.
\end{equation}
Since $24\cdot19=456\equiv 1\pmod{65}$, this can be further rewritten as
\begin{equation}\label{cw}
  24(Y-Y_0-19\gb) \bxmod 65 \ge 41.
\end{equation}
Hence, if we define
\begin{equation}\label{gam}
  \gamma= (-Y_0-19\gb)\bxmod{65},
\end{equation}
we can state the rule as:
\begin{Rule}
$Y$ is a leap year if and only if
\begin{equation}\label{lygam}
  24(Y+\gamma) \bxmod 65 \ge 41.
\end{equation}
\end{Rule}
Equivalently, if we define
\begin{equation}\label{gamx}
  \gamx=24\gamma \bxmod 65
 =(-24 Y_0-\gb)\bxmod 65,
\end{equation}
then:
\begin{Rule}
$Y$ is a leap year if and only if
\begin{equation}\label{lygamx}
  (24\,Y+\gamx) \bxmod 65 \ge 41.
\end{equation}
\end{Rule}

Furthermore, it easily seen from \eqref{leaprule-gb0}, \cf{} \eqref{cu}, 
that if $Y$ is a leap year, then the leap month has number
\begin{equation}\label{mk}
M = 1+\Floor{\frac{64-\bigpar{24(Y-Y_0)-\gb} \bxmod 65}2 }
\end{equation}
which, using \eqref{gamx}, easily is transformed to
\begin{equation}\label{mq}
M
 = \Floor{33-\frac{(24\,Y+\gamx) \mod 65}2 }
 = \Floor{33-\frac{24(Y+\gamma) \mod 65}2 };
\end{equation}
moreover,
if $Y$ is not a leap year, then \eqref{mk}--\eqref{mq} yield an impossible
value $M\ge13$. 

Finally, if $Y_1 \le Y_2$, the number of leap years (and thus the number of
leap months) in the period from $Y_1$ to $Y_2$ (inclusive) is
\begin{equation}\label{lys}
  \Floor{\frac{24(Y_2+1)+\gamx}{65}} -   \Floor{\frac{24Y_1+\gamx}{65}}.
\end{equation}
To see this it suffices to consider the case $Y_2 =Y_1$, which easily
follows from \eqref{lygamx}.

For the standard \PH{} version, 
\eqref{gam} can be written, using \eqref{gbgbx},
\begin{equation}\label{gam-P}
  \gamma= (-Y_0-19\gb)\bxmod{65}
=(-Y_0+19\gbx+14)\bxmod{65},
\end{equation}
and
the values of $\gbx$ in
\eqref{gbxS}--\eqref{gbxX} yield 
\begin{align}
\gamma &=42
\intertext{and thus}
  \gamx &=24\cdot 42\mod 65 = 33.  
\end{align}
(The values of $\gamma$ and $\gamx$ are the same for any epoch, since
the leap years are the same. However, for other versions of the calendar, the
formulas hold with other values of $\gamma$ and $\gamx$, 
see \refApp{Aversions}.)
Thus the rule is \cite{Salmi}:
\begin{Rule}
$Y$ is a leap year if and only if
\begin{equation}
  24(Y+42) \mod 65 \ge 41,
\end{equation}
\end{Rule}
or, equivalently,
\begin{Rule}
$Y$ is a leap year if and only if
\begin{equation}
  (24\,Y+33) \mod 65 \ge 41.
\end{equation}
\end{Rule}

For the \PH{} version, \eqref{mq} is
\begin{equation}
M = \Floor{33-\frac{(24\,Y+33) \mod 65}2 }.
\end{equation}

\section{Days}\label{Sdays}

As said in \refS{Sdef}, 
each lunar month (from the instant of new moon to the next new moon)
is (as in the Indian calendars \cite{CC}) divided into 30
lunar days; 
these have varying length of between 21.5 and 25.7 hours, and
do not correspond exactly to the
calendar days of 24 hours each. 
During each of these lunar days, 
the elongation of the moon
(\ie, the difference between lunar and solar longitude)
increases by 1/30 ($=12\grad$).

The calendar computations, unlike the Indian ones, do not include a
function calculating the elongation at a given time; instead the
computations use the inverse function, giving directly the time when the
elongation has a given value. We denote this function,
described in detail in \refS{Sastro}, by 
$true\_date(d,n)$, giving the date at the end of the lunar day $d$ in
\tmc{} $n$. The value of this function is a real (rational) number;
traditionally it is counted modulo 7, and the integer part yields the
day of week, but we will treat it as a real number so that the integer
part directly gives the Julian day number JD.
(A different constant $m_0$ below will give RD \cite{CC} instead.)
The fractional part shows the time the lunar day ends; 
it is used to calculate some further astronomical (and astrological)
information, see \refS{Sfurther}, but can be ignored for the present purpose.

The basic rule is:
\begin{Rule}
A calendar day is labelled by the lunar day
that is current at the beginning of the calendar day.   
\end{Rule}
In other words,
a lunar day gives its name (number and month) to the calendar day
where the lunar day 
ends. (Thus the JD of the calendar day is the 
integer part
of $true\_date$ at the end of the corresponding lunar day.) 
There are two special
cases covered by the rule above: if two lunar days end the same
calendar day, then the  calendar day gets the name of the first of them;
if no lunar day ends a given calendar day, then that day gets the same
name as the following day.
The first special case occurs when a lunar day is completely contained in one
calendar day; in that case no calendar day gets the number of this lunar
day, so this date is skipped. (In the sense that the number is skipped in the
numbering of days in the
calendar. The lunar day itself exists and can be used for astrological
purposes. The calendar days exist in real life and, of course, as such
cannot  be skipped or repeated.) 
The second case occurs when a lunar day completely contains a calendar
day; in that case this calendar day gets the same number of the next day, so
the date is repeated. (In the sense that the number is repeated.)

When a date is repeated, the first of the two days with the same number is
regarded as a leap day, and denoted ``Extra'' in the almanacs. 
\cite{Schuh-review} 
(But see \refS{Sholidays}.)

Recall also that each calendar day has a day of week, in the same way as in
Western calendars, see \refS{Sweek}. In particular, 
when a day is repeated, the
two days with the same number are distinguished by different days of week.


\begin{remark}\label{Rdate}
We do not have to worry about when the day starts; the calendar day is
from dawn 
to dawn, but the formulas take this into account (at least
theoretically), 
and no further modification is done.
In particular, no calculation of sunrise is required (as it is in Indian
calendar calculations \cite{CC}). 
Thus, the $true\_date$ should be regarded as a
kind of local Julian date that is offset from the standard
astronomical one which assumes integer values at noon UT (= GMT), so
that it instead assume integer values at local (mean) dawn.
(Henning \cite[pp.~10--11]{Henning} and \cite{Henning-comments} 
specifies the start of the day as
mean daybreak = 5 am local mean solar time.
Since the time difference is about 6 hours (Lhasa has longitude $91\grad$
which corresponds to $6^h4^m$), this is about
$-1$ UT, i.e.\ 11 p.m.\ UT the preceding day.)
\end{remark}

\begin{remark}
  By definition, new moon is (exactly) at the end of lunar day 30 and full
  moon is at the end of lunar day 15; the rules above imply that unless the
  day is skipped, (true) new moon falls in calendar day 30 and (true) 
full moon in calendar day 15 in every month. The true elongation differs
from the correct astronomical value by only about $2\grad$, corresponding to 4
hours, see \refS{SSdrift},  
so usually the same holds for the astronomcal new moon and full moon as well. 
\end{remark}

\section{Astronomical functions}\label{Sastro}

The $true\_date$ is calculated by first calculating a simpler
$mean\_date$, corresponding to the linear mean motion of the moon,
and then adjusting it by the equations of the moon and
sun, which are determined by the anomalies of the moon and sun
together with tables. (The tables are really approximation to sine,
suitably scaled. A similar table is used for each planet; no general
sine table is used, another difference from Indian calendars \cite{CC}.) 
The anomalies, in turn, are also calculated by linear functions. 

Traditional hand calculations calculate the mean quantities first for the
beginning of the month, \ie{} the end of the preceding month. (This
corresponds to taking day $d=0$ below. Note that this usually gives a time
during the last calendar day of the preceding calendar month.) 
Then the quantities are adjusted to give the values for a given day.
We will combine the two steps into one, giving the values for \tm{} $n$
and day $d$ directly.  
\xfootnote{
\citet{LaiDolma} refer to  tables rather than doing multiplications to obtain
the adjustments for $d$ for $mean\_date$ and (without explanation)
for $mean\_sun$, but the results are the same.
}

The \PH{} tradition uses the following functions. As said above,
these give the values at the end of lunar day $d$ in \tmc{} $n$.
I use below the epoch \eS, but give also the corresponding constants for \eH{}
and \eX. (Recall that the different epochs yield the same calendar.)

The mean date \tib{gza' bar pa}  
is (for \eS)
\begin{equation}\label{meandate}
mean\_date(d,n)=n\cdot m_1+d\cdot m_2+m_0,
\end{equation}
where
\begin{align}
m_1&=
29;31,50,0,480\rr{60,60,6,707}=\frac{167025}{5656} \quad(\approx 29.530587),
\label{m1}
\\
m_2&=0;59,3,4,16 \rr{60,60,6,707}=\frac{11135}{11312}
\quad\Bigpar{=\frac{m_1}{30}},
\\
m_0&=0;50,44,2,38 \rr{60,60,6,707}+2015501 =2015501+\frac{4783}{5656}.
\label{m0}
\end{align}

\begin{remark}\label{R7}
The traditional reckoning counts days modulo 7 only, \ie{} day of week
(see \refS{Sweek});
  this is of course enough to construct a calendar month by month. 
My version gives the JD directly. To be precise, the traditional
  result is obtained by adding 2 to the value in \eqref{meandate} before
  taking the  remainder modulo 7, 
\cf{} \eqref{weekday}.
(Hence, the integer added to the traditional value of $m_0$ has to be 
congruent to   $-2$ modulo 7.
Indeed, 
  $2015501\equiv -2 \pmod7$, and similarly for the other epochs below and in 
\refApp{Aversions}.)

The traditional value is therefore 
$m_1=1;31,50,0,480 \rr{60,60,6,707}$, 
subtracting 28 from the value used in this paper.
(It can also, and perhaps better, be regarded as
$1,31,50,0,480 \rr{7,60,60,6,707}$, with the denominator 7 meaning that
we regard the numbers as fractions of weeks, thus ignoring integer parts). 
Similarly, the constant 2015501 in $m_0$ (to get the result in JD) is my
  addition and not traditional.
(To get RD, use 294076 instead.)
\end{remark}

\begin{remark}
For \eH{} \cite{Henning}, 
a simple calculation using \eqref{tm0} shows that the \tmc{} $n$ differs by
13866 from the value for \eS; thus
$m_0$ is instead
\begin{align*}
m_0+13866\cdot m_1
&= 2424972+\frac{5457}{5656}
\\
&= 2015501+409471+\frac{5457}{5656}
\\
&
=2015501+409465+6;57,53,2,20 \rr{60,60,6,707},
\end{align*}
traditionally written as $6;57,53,2,20 \rr{60,60,6,707}$.
(Note that $409465=7\cdot58495$ a multiple of 7, \ie{} a
whole number of weeks; hence this gives the correct day of week.)

Similarly,
for \eX{} \cite{LaiDolma}, 
with the value of $n$ differing by 14609, $m_1$ and $m_2$ are the same
but $m_0$ is instead
\begin{align*}
m_0+14609\cdot m_1
&= 2446914+\frac{135}{707}
\\
&= 2015501+431413+\frac{135}{707}
\\
&
=2015501+431410+3;11,27,2,332 \rr{60,60,6,707},
\end{align*}
traditionally written as $3;11,27,2,332 \rr{60,60,6,707}$.
\end{remark}

Similarly, the mean longitude of the sun \tib{nyi ma bar pa} 
is
\begin{equation}\label{meansun}
mean\_sun(d,n)=n\cdot s_1+d\cdot s_2+s_0,
\end{equation}
where, 
\xfootnote{\label{f67-707}%
Note that the radices used for longitude are not the same as the
  radices used for time, see \eg{} \eqref{m1} and \eqref{s1}; 
the last radix is 67
  instead of 707. This ought to cause some problems when converting between
  times and longitudes, but the difference is only in the last term 
and is usually ignored.
}
\begin{align}
s_1&=2,10,58,1,17 \rr{27,60,60,6,67}=\frac{65}{804} 
\quad\Bigpar{=\frac{65}{12\cdot67}},   \label{s1}
\\
s_2&=0,4,21,5,43 \rr{27,60,60,6,67}=\frac{13}{4824}
\quad\Bigpar{=\frac{s_1}{30}}, \label{s2}
\\
s_0&=24,57,5,2,16 \rr{27,60,60,6,67} =\frac{743}{804}. \label{s0}
\end{align}

\begin{remark}
For \eH{} \cite{Henning}, 
with the base for $n$ differing by 13866, $s_1$ and $s_2$ are the same
but 
$s_0$ is instead (modulo 1, since only the fractional part matters)
\begin{align}\label{s0E1927}
s_0+13866\cdot s_1
\equiv
25,9,10,4,32 \rr{27,60,60,6,67}
=\frac{749}{804}
.
\end{align}
For \eX{} \cite{LaiDolma}, 
with the base for $n$ differing by 14609, 
$s_0$ is instead
\xfootnote{
This vanishing of the initial value, which recurs every 65th year
(804th month), 
is called \tibx{nyi ma stong bzhugs} 
``sun empty enter'' \cite{LaiDolma}.} 
\begin{align}\label{s0E1987}
s_0+14609\cdot s_1 \equiv 0.
\end{align}
\end{remark}

Thirdly, 
the anomaly of the moon 
\tib{ril-po dang cha-shas}  
is 
\begin{equation}\label{anomalymoon}
anomaly\_moon(d,n)=n\cdot a_1+d\cdot a_2+a_0,  
\end{equation}
where (but see \refR{R30} below for an alternative for $a_2$) 
\begin{align}
a_1&=2,1 \rr{28,126}=\frac{253}{3528},
\\
a_2&=1,0 \rr{28,126}=\frac1{28},   \label{a2}
\\
a_0&=3,97 \rr{28,126}=\frac{475}{3528}.
\end{align}

\begin{remark}
For \eH{} \cite{Henning}, 
with the base for $n$ differing by 13866, 
$a_0$ is instead (modulo 1, since only the fractional part matters)
\begin{align}
a_0+13866\cdot a_1
\equiv 13,103 \rr{28,126}=\frac{1741}{3528}.
\end{align}
For \eX{} \cite{LaiDolma}, 
with the base for $n$ differing by 14609, 
$a_0$ is instead
\begin{align}
a_0+14609\cdot a_1 \equiv 21,90 \rr{28,126}=\frac{38}{49}.
\end{align}
\end{remark}

\begin{remark}
For comparisons with modern astronomical calculations, note that
the anomaly is measured from the Moon's apogee, while Western
  astronomy measures it from the perigee, which makes the values
  differ by a half-circle.
\end{remark}

The equation of the moon \tib{zla rkang} is calculated by
\begin{equation}\label{moonequ}
moon\_equ=moon\_tab(28\cdot anomaly\_moon)  
\end{equation}
where $moon\_tab(i)$ is listed in the following table for $i=0,\dots,7$,
which extends by the symmetry rules
$moon\_tab(14-i)=moon\_tab(i)$,
$moon\_tab(14+i)=-moon\_tab(i)$, and thus
$moon\_tab(28+i)=moon\_tab(i)$;
linear interpolation is used beween integer arguments.
\begin{equation}\label{moontab}
\begin{tabular}{l r r r r r r r r}
$i$ &0 & 1 & 2 & 3 & 4 & 5 & 6 & 7\\
$moon\_tab(i)$ & 0 & 5 & 10 & 15 & 19 & 22 & 24 & 25
\end{tabular}  
\end{equation}

To find the equation of the sun \tib{nyi rkang}, first calculate the anomaly by
\begin{equation}
  \label{anosun}
anomaly\_sun=mean\_sun-1/4
\end{equation}
and then take
\begin{equation}\label{sunequ}
sun\_equ=sun\_tab(12\cdot anomaly\_sun)  
\end{equation}
where $sun\_tab(i)$ is listed in the following table for $i=0,\dots,3$,
which extends by the symmetry rules
$sun\_tab(6-i)=sun\_tab(i)$,
$sun\_tab(6+i)=-sun\_tab(i)$, and thus
$sun\_tab(12+i)=sun\_tab(i)$;
linear interpolation is used beween integer arguments.
\begin{equation}\label{suntab}
\begin{tabular}{l r r r r r r r r}
$i$ &0 & 1 & 2 & 3\\
$sun\_tab(i)$ & 0 & 6 & 10 & 11 
\end{tabular}  
\end{equation}

The date at the end of the lunar day \tib{gza' dag} 
is finally calculated as 
\begin{equation}\label{truedate}
true\_date =mean\_date +moon\_equ/60-sun\_equ/60.  
\end{equation}
(The half-corrected $mean\_date +moon\_equ/60$ is called 
semi-true date \tib{gza' phyed dag pa}.)   

Similarly, although not needed to calculate the calendar date, the
true solar longitude  
\tib{nyi dag} is  
\begin{equation}\label{truesun}
true\_sun=mean\_sun-sun\_equ/(27\cdot 60).
\end{equation}

\begin{remark}
You will not find the factors 1/60 and $ 1/(27\cdot 60)$ explicit in the
references; they are consequences of the positional system with mixed
radices.
Furthermore, $1/4$ in \eqref{anosun} is traditionally expressed as
$6,45\rr{27,60}$. 
\end{remark}

\begin{remark}\label{R30}
  We have $m_2=m_1/30$ and $s_2=s_1/30$, which is very natural, since
  it means that the functions $mean\_date$ and $mean\_sun$ are linear
  functions of the lunar day count $d+30\cdot n$, and thus increase by the
  same amount every day without any jumps at the
  beginning of a new month.
\xfootnote{\label{fm'}%
Schuh \cite{Schuh}, \cite[Kalenderrechnung]{tibetenc}
gives several examples ($m=3,4,5,6,7$ in his notation) of earlier versions
of the calendar (with $m_1$ and $s_1$ slightly different from the values
above), where simpler, rounded, values of $m_2$ and $s_2$ were
used. Moreover, most of these versions used $m_1$ and $s_1$ only for
calculating for the first month each year, and simplified value $m_1'$ and
$s_1'$ for the increments for the successive months within the year.
Some of them ($m=3,6,7$) also used the simplified $a_1'=2,0 \rr{28,126}\equiv
30 a_2$ for
the monthly increments of the anomaly.
See \refApp{AKT} for an example.
I do not know any currently used versions of the calendar that use such
simplifications. 
}

For the anomaly, however, the standard value $a_2=1/28$ does not
conform to this. Note that $a_1$ could be replaced by
$1+a_1=30,1 \rr{28,126}$ since we count modulo 1 here;
in fact, this is the ``real'' value, since the 
astronomical anomaly increases by 1 full circle 
in a little less that one month, see also
\eqref{meanano}.
Moreover, $a_2=1,0\rr{28,126}=1/28$ is a close approximation to
$(1+a_1)/30$; the 
conclusion is that one usually for convenience uses the rounded value
$a_2=1/28$. \citet{Henning}, however, uses instead the exact value
\begin{equation}\label{a2lochen}
  a_2=\frac{1+a_1}{30}=\frac{3781}{105840}=\frac{1}{28}+\frac1{105840}
=1,0,1 \rr{28,126,30};
\end{equation}
this is also used in his computed calendars 
\cite[Traditional Tibetan calendar archive]{kalacakra}. 
\xfootnote{
The value \eqref{a2lochen} was proposed by 
Minling Lochen Dharmashri (1654--1717).
\Hp.
}
The value \eqref{a2lochen} is mathematically more natural than
\eqref{a2}, since the latter value yields a (small) jump in the anomaly
at the end of each month while \eqref{a2lochen} yields same increase every
day.
Nevertheless, the simpler \eqref{a2} is usually used when calculating
calendars, for example in the almanac \cite{tib2013}.
\xfootnote{\label{fa2}%
This is witnessed by, for example, the values given for the true day of week
(including fractional part),  calculated by \eqref{truedate}.
This is given in the almanac with 3 terms (1,60,60), and 
the two calculations often differ (more often towards the end of the month), 
although
typically only by 1 unit in the third term.
A computer calculation for the 360 lunar days in 2013 yields
134 days with no difference in the true date (truncated to three terms),
149 days with a difference $\pm1$ in the third term, 75 days with $\pm2$ and
2 days with $\pm3$.
}

The difference $1/105840$ between \eqref{a2} and \eqref{a2lochen}
is small, and the resulting difference in the anomaly
is at most  $30/105840=1/3528$; the difference in the argument to
$moon\_tab$ is thus at most $28/3528=1/126$; since the increments in
$moon\_tab$ are at most 5, the difference in $moon\_equ$ is at most
$5/126$; finally, by \eqref{truedate}, the difference in the true
date is at most $(5/126)/60 =1/1512$
(see also \refF{fa2}),
so when rounding to an integer
(see \eqref{jd} below), we would expect to obtain the same result
except, on the average, at most once in 1512 days.
We would thus expect that the two different values of $a_2$ would give
calendars that differ for at most one day in 1512 on the average.
Moreover, since this was a maximum value, and the average ought to be
less by a factor of about 1/2, or more precisely $15.5/30=31/60$, since the
difference is proportional to the day $d$ which on the average is
$15.5$, and by another factor of 
$5/7$ since the average derivative of $moon\_tab$ is $25/7$ while we
just used the maximum derivative 5.
(For a sine function, the average absolute value of the derivative is
$2/\pi$ times the maximum value, but the ratio is closer to one for
the approximation in $moon\_tab$.)
Hence we would expect the average (absolute) difference in
$true\_date$ to be about 
\begin{equation*}
  \frac{31}{60}\cdot\frac{5}{7}\cdot\frac1{1512}
\approx\frac1{4100}.
\end{equation*}
Consequently, we expect
that the two versions of $a_2$ will lead to different
Tibetan dates for, on the average, one day in 4100, 
or a little less than one day in 10 years. 
This is confirmed by a computer search finding 9 examples in 1900--1999 and 8
examples in 2000--2099.
\xfootnote{On the other hand, there are 16 examples in the first 65 year
  cycle 1027--1091, so the occurences are irregular.}
Some recent examples are
10 February 2001 (JD 2451951) and
10 May 2006 (JD 2453866); the next example is
19 November 2025 (JD 2460999).
(It would be interesting to check these dates in published calendars.)
\end{remark}

\begin{remark}\label{Rlunarday}
  From one day to the next, the anomaly of the moon increases by 1/28 (or
  slightly more if \eqref{a2lochen} is used, and by the slightly larger
$1,1\rr{28,126}= 127/3528$ at each new month if \eqref{a2lochen} is not
used; we can ignore these differences); 
this means that the arguments used in \eqref{moonequ} differ by 1,
and thus the resulting values of $moon\_equ$ differ by at most 5 (see
\eqref{moontab}). Similarly, the anomaly of the sun differs by slightly less
than $1/360$, so by \eqref{sunequ} and \eqref{suntab}, the values of
$sun\_equ$ differ by at most about $1/5=0.2$. 
The mean date increases by $m_2 =0.98435$ each lunar day;
this is thus the length of the mean lunar day, measured in calendar days,
which equals $24m_2=23.6245$ hours.
By \eqref{truedate} and the calculations just made 
we see that the increase of the true date from one lunar
day to the next, \ie, the length of the true lunar day,
differs from this by at most $\pm (5+0.2)/60=0.087$ days, or $2.1$ hours. 
Consequently, the length of the (true) lunar day varies between, roughly,
21.5 and 25.7 hours.

A similar calculation shows that the length of a (true) lunar month varies
between 29.263 and 29.798 days, or about $29^d6^h$ and $29^d19^h$.
Similarly, a lunar year of 12 lunar months has length
between 354.00
and 354.74 days,
and a lunar year of 13 lunar months has length between
383.67 and 384.13 days.


The length of the calendar year is one of these lengths for a lunar year
rounded up or down to one of the nearest integers; the possibilities are
thus 354, 355, 383, 384, 385.
A computer calculation (for 10000 years) gave the frequencies:
$$
\begin{tabular}{rrrrr}
354&355&383&384&385\\\hline
42\%&21\%&3\%&33\%&1\%  
\end{tabular}.
$$
\end{remark}

\begin{remark}\label{Repoch2}
  The epoch is a specific day (or instant), \viz{} the mean new moon at
the beginning of the month with \tmc{} 0 (see \refR{Repoch} for this month).
By \eqref{meandate}, the mean date of the epoch is $m_0$ (since $n=d=0$);
hence the JD of the epoch is $\floor{m_0}$. 

Our version thus encodes the epoch in $m_0$, and the epoch is not needed
explicitly. Traditionally, with $m_0$ given modulo 7, see \refR{R7}, $m_0$
gives only the day of week of the epoch (with $\floor{m_0}$ giving the
number according to \refT{T7}, see \eqref{weekday}), and the epoch date is
needed to construct the calendar.

The epoch is always close to the beginning of month 2 or 3, depending on
$\gbx$ (the initial value of the intercalation index), see \refR{Repoch}.
Since the calculation is for the beginning of lunar day 1 
(which equals the end of the last lunar day of the preceding lunar month),
the epoch is usually day 30 of calendar month 1 or 2, respectively,
\cf{} the rule for labelling calendar days in \refS{Sdays};
however, since the true date differs from the mean date, see
\eqref{truedate}, it may be one of the adjacent calendar days 
(29/1 or 1/2, or 29/2 or 1/3).
For our three example epochs, the Tibetan dates are 
29/1 806, 
29/1 1927
and 
1/3 1987.
\end{remark}

\section{Calendrical functions}\label{Scal}

As explained in \refS{Sdays},
the Julian day number JD of a given Tibetan date can be calculated as the
JD of the 
calendar day containing the end of the corresponding lunar day, i.e.\ 
\begin{equation}\label{jd}
\JD=\floor{true\_date},  
\end{equation}
with $true\_date$ given by \eqref{truedate},
except that if the Tibetan date is repeated, this gives the JD
of the second day; for the first we thus have to subtract 1.
If we do the calculations for a day that is skipped, formula
\eqref{jd} will still give the JD of the calendar day that lunar day
ends, which by the rule in \refS{Sdays} is the day with name of the preceding
lunar day (since that lunar day ends the same calendar day).

The calendar is really given by the inverse of this mapping; a day is
given the number of the corresponding lunar day (\ie, the lunar day ending
during the calendar day).
To find the Tibetan date for a given JD, we thus compute approximate
\tmc{} and day (using the mean motion in \eqref{meandate}) and then
search the neighbouring lunar days for an exact match, if any, taking care
of the special cases when there are 0 or 2 such lunar days. For a
detailed implementation, see 
Dershowitz and Reingold \cite{CC}.

\subsection*{Beginning and end of months}
To find the last day of a month, we can just compute the JD for lunar
day 30 of the month; this gives the correct result also when day 30
is repeated or skipped. To find the first day, however, requires a little
care since lunar day 1 may end during the second day of the month
(when day 1 is repeated) or during the last day of the preceding month
(when day 1 is skipped); the simplest way to find the JD of first day of a
month is to add 1 to the JD of the last day of the preceding month.
By the comment on skipped days above, this can be computed as 1 + the
JD of day 30 the preceding month, regardless of whether day 30 is
skipped or not.
(There seems to be errors in the tables in \cite{Schuh} due to this.)

\subsection*{Tibetan New Year}
The Tibetan New Year \tib{Losar}
is celebrated starting the first day of the year. Since the
first month may be a leap month 1
(which happens twice every 65 year cycle,
the last time in 2000) and the first day may be day 2 (when day 1 is
skipped), some care has to be taken to calculate the 
date. The simplest is to add 1 to the JD of the last day the preceding
year, which thus is 1 + JD of the last day of (regular) month 12 the
preceding year (and can be calculated as 1 + JD of day 30 of (regular)
month 12  the preceding year).

\begin{remark}\label{Rlosar}
Holidays are otherwise usually not celebrated in leap months. 
However, the Tibetan New Year is really the first day of the year
even in a year that begins with a leap month 1. 
\xfootnote{
This is verified by
\cite{Salden} 
(written by a representative of the personal monastery of the Dalai Lama), 
according to which
\tibx{Losar} was
celebrated on Sunday, February 6th (2000), which was the first day of
leap month 1.
}
\end{remark}

\section{Day of week}\label{Sweek}
Each calendar day is, as explained in Sections \refand{Sdays}{Scal},
given a number in the range 1,\dots,30.
It also has a day of week, which as in the Gregorian and many other
calendars simply repeats with a period of 7.
The day of week thus corresponds uniquely to the Western day of week.

The days of week are numbered 0,\dots,6, and they also have names.
Each day of week corresponds to one of the seven ``planets'' (including sun and
moon), and the name of the day is more or less identical to the name of the
planet.
\xfootnote{\label{fplanets}%
The correspondence between days of week and planets in Tibetan 
is the same as in Latin and (incompletely) in many modern European languages
(e.g., Sunday, Monday and Saturday in English).
This correspondence goes back to Roman astrology, see 
\cite[pp.~268--273 and 391--397]{Richards} and \cite[\S12.12]{AA}, 
and came to Tibet through India.
}
The correspondence with
the English names is given in \refT{T7}.
(For the last column, see \refApp{Aastro}.)

\begin{table}[!htpb]
\begin{tabular}{r l l l l}
 & English & Tibetan & planet & element\\
\hline
0 & Saturday & spen ma & Saturn & earth \\
1 & Sunday & nyi ma & Sun & fire \\
2 & Monday & zla ba & Moon & water \\
3 & Tuesday & mig dmar & Mars & fire \\
4 & Wednesday & lhag pa & Mercury & water \\
5 & Thursday & phur bu & Jupiter & wind \\
6 & Friday & pa sangs & Venus & earth 
\end{tabular}  
\caption{Days of week.}
\label{T7}
\end{table}

The day of week is simply calculated from the Julian day number found
in \eqref{jd} by
\begin{equation}\label{weekday}
  day\_of\_week=(\JD+2) \mod 7.
\end{equation}

The date (\ie, the number of the day) and the day of week are often
given together. This resolves the ambiguity when days are repeated.
(It also helps to resolve most ambiguities when different rules of
calculation may have been used.)

\section{Further calculations}\label{Sfurther}
Tibetan almanacs traditionally also contain further information, 
mainly for astrological purposes, see \citet[Chapter IV]{Henning}.
(See also his extensive computer calculated
examples \cite[Traditional Tibetan calendar archive]{kalacakra}
and programs with explanations
\cite[Open source Tibetan calendar software]{kalacakra},
\cite[Open source Tsurphu calendar software]{kalacakra}.)
In fact, the calendar is known as ``the five components''
\tib{lnga-bsdus},
\xfootnote{
Just as in Indian calendars \tib{panchang}, 
see \cite[Chapter 18, p.~312]{CC}.
}
where the five components are: the day of week, 
the lunar day, the lunar mansion, the yoga and the \karana{} 
(for these, see below). 
\xfootnote{
\cite[Kalenderrechnung]{tibetenc} has a slightly different interpretation,
including the longitude of the sun as one of the five
and omitting the day of week (which is seen as given).
}
Nevertheless, 
a complete Tibetan almanac
typically contains not only these five but also some further data.
(The following description is based on \cite{Henning} and the almanac
\cite{tib2013}, see also the examples in \cite[pp.~202--203]{Henning},
\cite{Schuh-review} and \cite[Kalenderrechnung]{tibetenc}.
There are certainly minor variations between different almanacs.)

The daily data in an almanac include (typically) the following.  
See further \refApp{Aastro}.
The numbers are usually truncated to 3 significant terms, so not all radices
given below are used.

\begin{romenumerate}[-10pt]

 \item \label{q-gza-dag}
The true day of week \tib{gza' dag}, 
\ie{}
the day of week and fractional part of the day when the lunar day
ends. This is
$(true\_date+2)\mod 7$, 
calculated by \eqref{truedate} and
written with the radices (60,60,6,707).
The integer part (\ie, the first digit) of the true day of week
is thus the day of week;
its name (given by  \refT{T7})
is also given in letters, together with the (lunar) date.

  \item\label{q-tshes-khyud}
The (true) longitude of the moon at the end of the lunar day 
\tib{tshes 'khyud zla skar}, 
calculated for lunar day $d$, \tm{} $n$ by
\begin{equation}\label{tseskhyud}
moon\_lunar\_day(d,n)=true\_sun(d,n)+d/30
\end{equation}
and written in lunar mansions with the radices (27,60,60,6,67).
(The rationale for \eqref{tseskhyud} is that the moon's elongation, \ie{}
the difference between lunar and solar longitude, by definition is
$d/30$ at the end of lunar day $d$, see \refS{Sdef}.)

  \item\label{q-zla-skar}
The (true) longitude of the moon at the beginning of the calendar day
\tib{res 'grogs zla skar} 
calculated from the values in \eqref{tseskhyud} and \eqref{truedate}
by (recalling that $\FRAC(x)$ denotes the fractional part of $x$)
\begin{equation}\label{zla-skar}
moon\_calendar\_day=moon\_lunar\_day-\FRAC(true\_date)/27,  
\end{equation}
and written in lunar mansions with the radices (27,60,60,6,67).
(The idea here is that $\FRAC(true\_date)$ is the time from the beginning of
the calendar day to the end of the lunar day, so this formula is really an
approximation assuming that the moon moves with constant speed and making a
full circle in 27 days.
\xfootnote{\label{f27}%
This is of course not exact, but it is a rather good approximation
since the tropical
\xpar{or sidereal; the difference is negligible} 
month is 27.322 days 
\cite[Table 15.3]{AA}. 
In the calculation, one may also ignore that 
in the traditional notations, the last radix differs betwen the true date
and the longitudes, \cf{} \refF{f67-707}.
}
The division by 27 is convenient since it is
simply a shift of the terms in the mixed radix notation.)

\item\label{q-naksatra}
The name of the lunar mansion (Sanskrit \emph{naksatra}), which is
determined by the Moon's longitude \eqref{zla-skar} 
with fractions of mansions ignored.
More precisely, numbering the lunar mansions from 0 to 26,
the number of the lunar mansion is
\begin{equation}\label{mansion}
\floor{27moon\_calendar\_day}.
\end{equation}
(In the traditional notation, this is the first term of the lunar longitude
\eqref{zla-skar}.) 
The names of the mansions (in Tibetan and Sanskrit)
are listed in \cite[Appendix I]{Henning}.

  \item \label{q-nyi-dag}
The (true) longitude of the sun \tib{nyi dag},
given by 
$true\_sun$ 
and written in mansions with the radices (27,60,60,6,67).
(Note that $true\_sun$ really is computed for the end of the lunar
day, but is regarded as valid also for the calendar day; 
the motion of the sun during the day is thus ignored
and no correction as in \eqref{zla-skar} is made.  
\xfootnote{
The motion of the sun is much slower than the motion of the moon; it is
about $1/365\approx1\grad$ each day.
}) 

\item \label{q-yoga-long}
The \emph{yoga} ``longitude'' \tib{sbyor ba} is the sum of the
longitudes of the sun and the moon, calculated by
\xfootnote{
As noted by \cite{Henning}, 
this is inconsistent, since we add one longitude at the beginning of the
calendar day and one longitude at the end of the lunar day. On the other
hand, this addition has in any case no physical or astronomical meaning.
}
\begin{equation}
\mathit{yoga\_longitude}=moon\_calendar\_day+true\_sun
\pmod1,
\end{equation}
and written in mansions with the radices (27,60,60,6,67).

\item \label{q-yoga}
The name of the yoga, which is
determined by the yoga longitude with fractions of mansions ignored.
More precisely, numbering the yogas from 0 to 26,
the number of the yoga is
\begin{equation}
\floor{27\mathit{yoga\_longitude}}.
\end{equation}
The names of the yogas (in Tibetan and Sanskrit)
are listed in \cite[Appendix I]{Henning}.
(The names differ from the names of the lunar mansions.)

\item \label{q-karana}
The \emph{\karana} \tib{byed pa}
in effect at the start of the calendar day. Each
lunar day is divided into two halves, and each half-day is assigned
one of 11 different \karana{s}. There are 4 ``fixed'' \karana{s} that occur
once each every month: the first half of the first lunar day, the
second half of the 29th lunar day, and the two halves of the 30th day;
the remaining 7 \karana{s}, called ``changing'', repeat cyclically for
the remaining 56 half-days. In other words, lunar day $D$ consists of
half-days $2D-1$ and $2D$, and half-day $H$ has one of the fixed
\karana{s} if $H=1,58,59,60$, and otherwise it has the changing \karana{}
number $(H-1)\amod 7$.
The names of the \karana{s} (in Tibetan and Sanskrit)
are listed in \cite[Appendix I]{Henning}.
(The exact rule for determining the time in the middle of the lunar
day that divides it into two halves is not completely clear to
me. \citet{Henning} divides each lunar day into two halves of equal lengths, but
there might be other versions.)

\item  \label{q-nyi-bar}
The mean longitude 
of the sun 
\tib{nyi bar}
in signs and degrees (and minutes),
\ie{} $mean\_sun$ \eqref{meansun} written with the radices (12,30,60).

\item 
\label{q-zla-skar-K}
The (true) longitude of the moon at the beginning of the calendar day
\tib{zla skar} as in \ref{q-zla-skar}, but calculated according to the
\karana{} 
\xfootnote{The (Sanskrit) word \karana{} has a different meaning here
than in \ref{q-karana}. In contrast, the standard calculations described
above are called \siddhanta. See further \refApp{AKT}.
}
calculations, see \refApp{AKT}. 
\xfootnote{
This means that all calendar
calculations have to be done also for the \karana{} version, in order to
find this value.}
This is, usually at least, calculated for the correct calendar day,
regardless of whether its 
date is the same in the \PH{} and \karana{} versions or not; so for example
in the almanac \cite{tib2013}. 
(There are exceptions, possibly mistakes, for example in the page 
from a 2003 calendar shown in \cite[p.~202]{Henning}.)
\item \label{q-greg}
The Gregorian date. In modern almanacs written in Western (European) numerals.
\end{romenumerate}

As an example, the almanac \cite{tib2013}
lists for each day,
in addition to the Tibetan (lunar) date 
and the Gregorian date \ref{q-greg}, 
the six numbers
\ref{q-gza-dag},
\ref{q-zla-skar},
\ref{q-nyi-dag},
\ref{q-yoga-long},
\ref{q-zla-skar-K},
\ref{q-nyi-bar}
above 
\xfootnote{These numbers are truncated (not rounded) to 3 terms,
with radices $(1,60,60)$ or (27,60,60).
%
};
there is further (astrological) information as text.
\xfootnote{The same daily values are given  in the 2003 almanac
shown in \cite[p.~202]{Henning}, see further the discussion there.
} 

Note that \ref{q-gza-dag}, \ref{q-tshes-khyud}, \ref{q-nyi-dag}, 
\ref{q-nyi-bar}, 
refer to (the end of) the lunar day
while the others
refer to (the beginning of) the calendar day.
When a day is skipped (\ie, there is a lunar day without corresponding
calendar day), the almanac usually still gives
\ref{q-gza-dag} and \ref{q-nyi-dag} for the skipped lunar day;
when a day is repeated, so there are two calendar days corresponding to the
same lunar day, the data are given for both days, with the data referring to
the lunar day thus repeated, but
with modifications for the first day: In \cite{tib2013},
\ref{q-gza-dag} (otherwise the end of the lunar day) for the first day
is given as the end of the calendar day, written as
$x;60,0$ where $x$ is the day of week; 
\ref{q-zla-skar} is obtained as
$moon\_lunar\_day-1/27$,
\xfootnote{
This is a reasonable approximation, since the end of the lunar day is just a
little more than one calendar day later than the beginning of the first
calendar day, \cf{} \refF{f27}. However, the resulting longitude is often
\emph{smaller} than the longitude at the end of the preceding lunar day, a
short time earlier. (The latter longitude is not printed in the almanac, so
the contradiction is not visible without further calculations.)
}
\cf{} \eqref{zla-skar};
\ref{q-yoga-long} is calculated separately for both days;
\ref{q-nyi-bar} is given only for the second day.
Thus only \ref{q-nyi-dag} is identical for the two days.

There is also further information at the beginning of each month, including
the mean date \tib{gza' bar} \eqref{meandate} and the mean solar longitude 
\tib{nyi bar} \eqref{meansun} (both given with all 5 terms),
and also the lunar anomaly \tib{ril cha} \eqref{anomalymoon}, all calculated 
for the beginning of the lunar month (day $d=0$ in the formulas above);
furthermore, the \tm{} \tib{zla dag} is given, with its fractional part
(the intercalation index), see \refR{Rix}. 
These monthly values (mean date, mean
solar longitude, lunar anomaly and \tm) are also given for the \karana{}
calculation (see \refApp{AKT}). 
There is also data on
the position of the planets (see \refApp{Aplanet}).

\subsection{Some special days}\label{Stent}
Furthermore, the almanac gives extra information on some special days, for
example 
the days when the mean solar longitude passes certain values, including
the three series (with $30\grad$ intervals)
\begingroup
\addtolength{\leftmargini}{-10pt}
\begin{itemize}
\item 
$0\grad, 30\grad, 60\grad\dots$ (when the mean sun enters a new sign);
\item $8\grad, 38\grad, 68\grad\dots$ 
(the \dpp{s} \tib{sgang}, see \refApp{SSdp});
\item $23\grad, 53\grad, 83\grad\dots$ 
(the midpoints between the \dpp{s} \tib{dbugs}).
\end{itemize}
(There are also 4 further such days, with longitudes
$66\grad$, 
$132\grad$, 
$147\grad$, 
$235\grad$ 
in the 2013 almanac \cite{tib2013}.)
In all cases, the almanac gives the longitude and the
mean date for the instant $mean\_sun$ reaches this value.
This is easily found from \eqref{meansun} and \eqref{meandate}, now
letting $d$ be an arbitrary rational number denoting the lunar day including
a fractional part to show the exact instant.
\xfootnote{This is the only case that I know when the calculations involve
  fractions of a lunar day.}
If we assume that $m_2=m_1/30$ and $s_2=s_1/30$, as is the case in the \PH{}
version and in all other modern versions that I know of, see \refR{R30},
then \eqref{meansun} shows that the mean solar longitude equals a given
value $\gl$ in the year $Y$ at lunar date
\begin{equation}
  d = \frac{\gl+Y-Y_0-s_0}{s_2} = 30\frac{\gl+Y-Y_0-s_0}{s_1} 
\end{equation}
after the epoch, where $Y_0$ is the epoch year (so $Y-Y_0+\gl$ can be
regarded as the desired longitude of the sun measured on a linear scale from
the epoch); however, in order for this to give the correct year,
$s_0$ has to be chosen such that $s_0$ is close to 0;
for \eX{} we thus take $s_0=0$ by \eqref{s0E1987}, but the values 
\eqref{s0} and \eqref{s0E1927} for \eS{} and \eH{}
have to be decreased by 1 (recall that the
integer part of $s_0$ was irrelevant earlier); moreover, the integer part of
$\gl$ should similarly be adjusted so that $0\le\gl<1$, except at the
beginning of the year (before longitude $0\grad$) 
when one should subtract 1 so that $-1<\gl<0$.
(Cf.\ \refApp{Aleap} for similar considerations regarding the mean solar
longitude as a real number on a linear scale.)
The mean date is then obtained from \eqref{meandate} as 
\begin{equation}\label{jw}
  d\cdot m_2+m_0 = 
\frac{m_2}{s_2} 
\xpar{\gl+Y-Y_0-s_0}+m_0
=\frac{m_1}{s_1} 
\xpar{\gl+Y-Y_0-s_0}+m_0.
\end{equation}
\endgroup 

The traditional calculation is somewhat different
\Hp.
First the month is determined (or guessed, for possible later correction)
-- this is easy since the mean sun always passes the definition point
$((M-3)\cdot 30+8)\grad$ during month $M$, see \refApp{Aleap}.
The lunar date of a new sign ($k\cdot30\grad$, for integer $k$) 
is given by $6ix/13$, where $ix$ is the
intercalation index of the month, see \eqref{ix} and \refR{Rix}.
\xfootnote{
Note that $0\le 6ix/13\le 6\cdot64/13=384/13<30$.
}
The \dpp{} $8\grad$ later comes $8;3,1\rr{13,5}=8\frac{16}{65}$ lunar days
later and 
the midpoint $7\grad$ earlier comes $7;2,4\rr{13,5}=7\frac{14}{65}$ lunar
days earlier. 
Finally, the mean dates for these lunar dates are calculated by a
multiplication equivalent to \eqref{meandate}; tables exist to assist with the
multiplication by fractional parts of a lunar day.

To verify that this method works,
we first note that $8\frac{16}{65}=8\cdot 1\frac{2}{65}$
and $7\frac{14}{65}=7\cdot 1\frac{2}{65}$, and that during 
$1\frac{2}{65}$ lunar day, the mean sun moves, see \eqref{s1}--\eqref{s2},
\begin{equation}
 1\frac{2}{65}\cdot s_2 = 
 \frac{67}{65}\cdot s_2 = 
 \frac{67}{65}\cdot \frac{s_1}{30}= 
 \frac{67}{65}\cdot \frac{65}{67\cdot12}\cdot\frac{1}{30}= 
\frac{1}{360} = 1\grad;
\end{equation}
hence the mean sun indeed moves $8\grad$ and $7\grad$ during these periods.

To verify the entry into signs, it is convenient to use the epoch \eX, since
then 
both $\gbx=0$ \eqref{gbxX} and $s_0=0$ \eqref{s0E1987}.
Consider  month $M$, year $Y$, and let as in \eqref{MM}
$\MM=12(Y-Y_0)+(M-M_0)$ be the number of solar months since the epoch 
(year $Y_0=1987$, month $M_0=3$).
If this month has \tmc{} $n$ and \ixx{} $ix$, then the rule above yields a
lunar date (recalling that each month has 30 lunar days)
\begin{equation}\label{mi}
x= 30n+\frac{6ix}{13} = 30\Bigpar{n+\frac{ix}{65}}
\end{equation}
lunar days after the epoch. Note that $n+ix/65$ is the true month.
If the \ixx{} $ix<48$, then $n+ix/65$ is thus given by \eqref{tm0}, and is
thus $\frac{65}{67}\MM$ (since $\gbx=0$), so the lunar date \eqref{mi}
is
\begin{equation}
  x = 30\cdot\frac{65}{67}\MM
\end{equation}
and consequently 
the mean solar longitude then is,
by \eqref{meansun} (with $s_0=0$) and \eqref{s1}--\eqref{s2}, 
  \begin{equation}
x\cdot s_2+s_0 = 30 \cdot\frac{65}{67} \cdot\MM \cdot \frac{s_1}{30}+0	
=s_1 \cdot\frac{65}{67} \cdot\MM 
=\frac{1}{12} \cdot\MM 
=\MM\cdot30\grad,
  \end{equation}
showing that the mean sun enters the sign with longitude $\MM\cdot 30\grad
\equiv (M-3)\cdot 30\grad$ (since longitudes are measured modulo $360\grad$).
If $ix\ge48$, then the
\tm{} is increased by $1\frac{2}{65}$, see \refR{Rix}, which corresponds to
an increase in the mean solar longitude of
$\frac{65}{67}s_1=\frac{1}{12}=30\grad$, 
so the mean sun instead enters the next sign, with longitude
$(M-2)\cdot30\grad$ 
and the entry into the sought sign is found (by the same rule) in the preceding
month. 
\xfootnote{It is easily verified that
$\frac{6\cdot47}{13} + 8\frac{16}{65} < 30 < \frac{6\cdot48}{13} +
  8\frac{16}{65}$, which implies that in any case, the \dpp{} 
$((M-3)\cdot 30+8)\grad$ is reached during month $M$ as said above,  
see further \refApp{Aleap}.
}
(If $ix=48$ or $49$, this means leap month $M$.)

\begin{remark}
  A complication for the almanac maker is that the lunar day found in this
  way sometimes corresponds to a calendar day that is one day earlier or
  later than the calendar day given by (the integer part of the) mean date
  found above, since this correspondence uses the true date and not the mean
  date. In this case, the data are entered in the almanac straddling both
  days, see \cite{Schuh-review} for details and examples.
\end{remark}

The days discussed here have been defined by the mean solar longitude,
including day when it is 0, \ie, the mean sun passes the first point of Aries.
The day when the true solar longitude is $0$ 
is also marked in the almanac; in this case the true date is
given. 
\xfootnote{The latter (``true'') day is about 2 days before the first
  (``mean''), because of the 
  equation of the sun, see \eqref{truesun}; both are over a
month later then the astronomical vernal equinox when the real sun has
longitude $0$, see \refS{SSdrift}.
}
I do not 
know exactly how this is calculated, but I assume that it is by some more
complicated version of the rule above, including adjustments for the
equations of sun and moon.

\section{Holidays}\label{Sholidays}
A list of holidays, each occuring on a fixed Tibetan date every year, is
given in  \cite[Appendix II]{Henning}.
If a holiday is fixed to a given date, and that date is skipped, the
holiday is on the preceding day.
If the date appears twice, the holiday is on the first of these 
(\ie, on the leap day, see \refS{Sdays}). 
(These rules are given by \cite{Berzin}, but I have not checked them
against published calendars.)

Holidays are usually not celebrated in leap months, but see \refR{Rlosar}.

\section{Mean lengths and astronomical accuracy}\label{Smean}
The mean length of the month is
\begin{equation}\label{meanmonth}
  m_1=29;31,50,0,480\rr{60,60,6,707}
=\frac{167025}{5656} \approx 29.530587
\text{ days},
\end{equation}
which is essentially identical to the modern astronomical value of the
synodic month $29.5305889$ (increasing by about $2\cdot10^{-7}$ each century)
\cite[(12.11-2)]{AA}.

Consequently, the mean length of the lunar day is
\begin{equation}
\label{meanlunarday}
\frac{m_1}{30}
=\frac{167025}{30\cdot5656} 
=\frac{11135}{11312} 
\approx 0.98435
\text{ days}.
\end{equation}

The mean length of the year is, \cf{} \eqref{6765},
\begin{equation}
\label{meanyear}
  \frac{1}{s_1}
\text{ months}
=  \frac{m_1}{s_1}
=\frac{804}{65}m_1
=\frac{6714405}{18382}
\approx 365.270645
\text{ days},
\end{equation}
which is $0.02846$ days more than 
the modern astronomical value of the
tropical year $365.24219$
\cite[(12.11-1) and Table 15.3]{AA},
and $0.02815$ days longer than the mean Gregorian year $365.2425$ days.
(It is also longer than the
sidereal year $365.25636$ days
\cite[Table 15.3]{AA}.)
Hence the Tibetan year lags behind and starts on the average
later, compared to the
seasons or the Gregorian year, by almost 3 days per century or almost a month
(more precisely, 28 days) per millennium.

The year has thus drifted considerably since it was introduced, and
even more since the \KT{} was written; the drift during the 1200
years since the epoch 806 is 34 days, see further \refS{SSdrift}.

The mean length of the anomalistic month is
\begin{equation}\label{meanano}
\frac{1}{1+a_1} \text{ months}
=\frac{3528}{3781} \cdot\frac{167025}{5656} 
=\frac{10522575}{381881} 
\approx  27.55459
\text{ days}.
\end{equation}
This agrees well with the modern astronomical value
$27.55455$ days \cite[Table 15.3]{AA}; the drift is about 1 day in 2000
years.  

See also \citet{Petri}.

\subsection{Accuracy of longitudes}\label{SSdrift}

As noted above, the mean length of the month is essentially equal to the
exact astronomical value. Indeed, the times (true date)
for new moons as calculated by
the Tibetan calendar are close to the exact astronomical times.
(Recall that the new moons signify the end of lunar day 30 in each month and the
beginning of a new month; in the almanac the time is thus given as the true
date of day 30 in the month.)
If we (following \cite{Henning}) regard the Tibetan day as starting at
mean daybreak, ca.~5 a.m.\ local mean solar time, 
and set this equal to $-1$ UT (GMT), see \refR{Rdate}, 
then
the true dates given by Tibetan calculations are about 4 hours earlier
than the exact astronomical values.
\xfootnote{A calculation for the 12 new moons in the Gregorian year 2013 
and comparison with a Swedish almanac gave differences 
between $3^h16^m$ and $4^h21^m$.
}
Since the elongation increases by about $12\grad$ each day, 
see \refS{Sdays},
this corresponds
to an error in the elongation of about $2\grad$ too large.

As also noted above, the mean length of the year is less accurate, and the
year has drifted 34 days since the epoch 806.
Indeed, the solar longitude as
computed by $true\_sun$ above passes $0\grad$ during
the 16th day in the third Tibetan month 2013, which equals April 26, 
37 days after the astronomical vernal equinox on March 20, 
which agrees well with this drift.
Another way to see this is to note that the true sun calculated at the
(astronomical) vernal equinox is only about $324\grad$, which means that the
solar longitude is about $36\grad$ too small. 
\xfootnote{This varies slightly
over the year, since the equation of the sun given by \eqref{sunequ} is not
astronomically exact, both because it is intrinsically an approximation
only, and because the anomaly is by \eqref{anosun} also about $36\grad$
wrong.
However, the error does not vary by more than 1--$2\grad$ over the year 
and is always about 36--$37\grad$.
}
Since the elongation was seen above to be about $2\grad$ too large,
the Tibetan longitude of the moon is about $34\grad$ 
(= 2.5 mansion) too small.

\subsection{The ratio between solar and lunar months}\label{SSratio}

The ratio of the exact astronomical lengths of solar and lunar months is
(with values correctly rounded for both today and 1000 years ago)
\begin{equation}
  \frac{365.2422}{12\cdot 29.53059} \approx 1.030689.
\end{equation}
The standard way to find good rational approximations of a real number is
to expand it as a continued fraction, see \eg{} \cite{HardyWright}:
\begin{equation}
  1.030689=
1+\frac{1}{32+\frac1{1+\frac1{1+\frac1{2+\frac1{2+\dots}}}}}
\end{equation}
and then calculate the partial quotients obtained by truncating the
continued fraction.
In this case the first partial quotients are (after 1)
\begin{align}
1+\frac{1}{32}&=\frac{33}{32}= 1.03125
\\
1+\frac{1}{32+\frac1{1}}&=\frac{34}{33}\approx1.030303
\\
1+\frac{1}{32+\frac1{1+\frac1{1}}}&=\frac{67}{65}\approx1.030769
\\
1+\frac{1}{32+\frac1{1+\frac1{1+\frac1{2}}}}&=\frac{168}{163}\approx1.030675
\\
1+\frac{1}{32+\frac1{1+\frac1{1+\frac1{2+\frac1{2}}}}}&=\frac{403}{391}
\approx 1.030690.
\end{align}
The first two approximations 33/32 and 34/33 have errors of about 
$5\cdot 10^{-4}$ and $4\cdot 10^{-4}$, or about 1 month in 200 years.
The next approximation is 67/65, used in the Tibetan calendar, which has an
error of about $8\cdot 10^{-5}$, or as said above about 1 month in 1000 years.
The next approximation is 168/163, with an error of only 
$1\cdot 10^{-5}$, about 1 month in 6\,000 years, and 403/391 has an error 
$2\cdot 10^{-6}$, about 1 month in 50\,000 years
\xfootnote{ 
The latter accuracy is only fictional; the error  50\,000 years
from now is much larger, $10^{-4}$, because of changes in the lengths of the
month and the year.}. 

The Tibetan approximation $67/65$, see \eqref{6765}, is thus a good choice
from a mathematical point of view, if we want a good approximation with
small numbers. (168/163 or 403/391 would yield more accurate but more
complicated calendars. As far as I know, they are not used in any calendar.)

\begin{remark}
An approximation used in several other calendars (\eg{} the Jewish,
and for calculation of the Christian Easter),
but not in the Tibetan,
is the relation 
235 lunar months $\approx$ 228 solar months (19 solar years), known as
Meton's cycle \cite{CC}. We have $235/228= 1.030702$, with an error of 
$1\cdot 10^{-5}$, about 1 month in 6000 years.
This does not appear in the list of partial quotients above, and indeed,
168/163 with smaller numbers is about as accurate. 
However, the Metonic  cycle has the advantage that it comprises a whole
number of years. (It appears as one of the partial quotients for the continued
fraction expansion of the number 
of lunar months in a solar year.) 
\end{remark}

\section{Period}\label{Speriod}

The $mean\_date$ given by \eqref{meandate} repeats after 5656 months,
but $true\_date$ depends also on $anomaly\_moon$ and $true\_sun$, which
repeat after 3528 and 804 months, respectively. The least common
multiple of 5656, 3528 and 804 is 
$p=23873976$; consequently, all astronomical functions above repeat
after $p$ months, 
which equals $m_1p=705012525$ days or 
$s_1p=1930110$ years (these are necessarily
integers). 
Note further that this period contains an integral number of leap year
cycles (804 months = 65 years), see \refS{Smonths},  so the numbering of
the months repeats too.
In other words:
\begin{Rule}
The calendar (days and months) repeats after
\begin{equation*}
 705\,012\,525\text{ days}
=
   23\,873\,976\text{ months}
=
1\,930\,110 \text{ years}.
\end{equation*}  
\end{Rule}
Moreover, the number of days is divisible by 7, so also
the day of week repeats after this period.

However, $1930110\equiv 6 \pmod{12}$, so to repeat the animal names of
the years in
the 12-year cycle, or the names in the 60-year cycle, we need two of
these periods, \ie{} 3\,860\,220 years.

Of course, the period is so long that the calendar will be completely
out of phase with the tropical year and thus the seasons long before,
and it will move through the seasons hundreds of times during one
period.

If we also consider the planets, see \refApp{Aplanet},
the period becomes a whooping
2\,796\,235\,115\,048\,502\,090\,600 years, see \citet[pp.~332--333]{Henning}
and \cite[Weltzeitalter]{tibetenc}. 

\clearpage

\appendix

\section{Different versions}\label{Aversions}

Different traditions follow different rules for the details of the
calculation of the Tibetan calendar, and as said in \refS{Sintro},
there are  two main versions in use today.
(Several other versions survive according to \cite[p.~9]{Henning},
but no details are given.)
Schuh \cite{Schuh} gives historical information, 
including many versions of the constants in \refS{Sastro} above 
used or proposed during the centuries,
but says nothing about the
different versions today. 
See also  \cite[Kalenderrechnung]{tibetenc}.
\citet{Henning} discusses the \PH{} and \TS{} versions in detail, and also
one recent attempt at reform. 
He gives epoch data for several versions in
\cite[Epoch data]{kalacakra}.
(Another source on different versions today is \citet{Berzin};
he gives the
impression that several versions are in use. 
However, since  \cite{Berzin} discusses the calendar and astrology
together, it is possible that some of these versions  actually use
the same calendar but differ in other, astrological, calculations or
interpretations.) 
See also the web pages of Nitartha \cite{nitartha}.

\subsection{\PH} (Phukluk, \tibx{phug-lugs}.) \label{APH} 
This is the most widespread version, and is regarded as the official
Tibetan calendar. 
(Also by at least some followers of the rival \TS{} tradition
\cite{nitartha}.)
It was started in 1447 
(the first year,  \tibx{rab byung}, of the 8th Prabhava  cycle)
by 
Phugpa Lhundrub Gyatso   
\tib{phug-pa lhun-grub rgya-mtsho}, 
and was used by the Tibetan government from at least 1696 to 1959
\cite[Phugpa-Schule]{tibetenc}, \cite[pp.~8 and 321--337]{Henning},
\cite[p.~139]{Schuh}, \cite{Schuh-review}.
The \PH{} version is used by the Gelug, Sakya, Nyingma and Shangpa Kagyu
traditions of Tibetan buddhism, including the Dalai Lama,
and is used \eg{} in the calendars published in Dharamsala in India, where
the Tibetan exile government resides, see \cite{tib2013,men-tsee-khang}.
The Bon calendar is the same as the \PH{} \cite{Berzin}.
The \PH{} calendar is described in detail in the main body of the present paper.

\subsection{\TS}(Tsurluk.)\label{ATS} 
This version was also introduced in 1447, by
Jamyang Dondrub Wozer 
\tib{mtshur-phu 'jam-dbyangs chen-po don-grub 'od-zer}
and
derives from 14th century commentaries to \KT{} by the
3rd Karmapa Rangjung Dorje of \TS{} monastery
(the main seat of the Karma Kagyu tradition of Tibetan
buddhism) \cite{Berzin},
\cite[pp.~9 and 337--342]{Henning},
\cite[Tshurphu-Schule]{tibetenc}.
This version is used by the Karma Kagyu tradition, and it is used 
\eg{} in calendars published 
by the Rumtek monastery in India, the main exile seat of 
the Karmapa (the head of Karma Kaygyu) 
\cite[Open source Tsurphu calendar software]{kalacakra}.
Also the calendar
published by
Nitartha in  USA \cite{nitarthacal,nitartha} gives the \TS{}
version (from 2004 the \PH{} version too is given).

The \TS{} tradition uses the astronomical functions in \refS{Sastro}
with the same values as given there for the \PH{} version of the
constants $m_1$, $s_1$ and $a_1$ (and $m_2$, $s_2$ and $a_2$, see \refR{R30})
for mean motions, while the epoch values $m_0$, $s_0$ and $a_0$ are different.
See further \cite[pp.~337--342]{Henning}.
\begin{remark}
  Traditional \TS{} calculations use, however, somewhat different sets
  of radices that the \PH{} versions, with one more term 
(and thus potentially higher numerical accuracy in the calculations). 
The same constants are thus written differently:
\begin{align}
m_1&=
29;31,50,0,8,584\rr{60,60,6,13,707}
=\frac{167025}{5656},
\\
s_1&
= 2,10,58,1,3,20 \rr{27,60,60,6,13,67}
=\frac{65}{804},
\end{align}
although  $m_1$ is traditionally given as
$1;31,50,0,8,584\rr{60,60,6,13,707}$, 
as always calculating modulo 7, and I have added 28, see \refR{R7}.
\end{remark}

Two epochs given in classical text \cite[p.~340]{Henning} are,
with 
$m_0$ modified to yield JD, see \refR{R7}, 
\begin{align}
\JD&=2353745 \quad \text{(Wednesday, 26 March 1732 (Greg.))} \notag
\\
m_0&=
4;14,6,2,2,666 \rr{60,60,6,13,707}+2353741
=2353745+\frac{1795153}{7635600},
\\
  s_0&=
- \bigpar{1,29,17,5,6,1 \rr{27,60,60,6,13,67}} 
=-\frac{5983}{108540},
\label{s0T1732}
\\
a_0&= 14,99 \rr{28,126}
=\frac{207}{392}.
\end{align}
and 1485 months later (equivalent and giving the same calendar)
\begin{align}
\JD&=2397598 \quad \text{(Monday, 19 April 1852)} \notag
\\
m_0&=
2; 9,24,2,5,417 \rr{60,60,6,13,707}+2397596
=2397598+\frac{1197103}{7635600},
\\
  s_0&=
0,1,22,2,4,18 \rr{27,60,60,6,13,67}
=\frac{23}{27135}, \label{s0T1852}
\\
a_0&=0,72 \rr{28,126}
=\frac{1}{49}.
\end{align}
We denote these by E1732 and E1852.
(See also \cite[Epoch data]{kalacakra}, where also a third epoch 1824 is
given.) 

To find the leap months, the true month is calculated by 
\eqref{tm0} with a constant $\gbx =59$ (E1732) or $14$ (E1852)
(the epoch value of the intercalation index);
the reason being  that counting backwards to the \Kc{}
epoch \nag{} (\Caitra) 806 yields the intercalation index 0,
as in the \KT. (But unlike the \PH{} version, see \eqref{gbxS}.)
Moreover, instead of the Phugpa rule \eqref{leaprule-P}, 
the \TS{} version uses the simpler rule:
\begin{equation}\label{leaprule-T}
\hskip-1em
\vbox {\narrower\narrower\narrower\noindent\em  
A leap month is inserted when the intercalation index\\
$ix = 0\text{ or\/ }1$.
}
\hskip-3em
\end{equation} 
This implies that for a regular month, the \tmc{} $n$ is obtained from
the true month \eqref{tm0} by rounding down to the nearest integer; for
a leap month we further subtract 1.
(The \TS{} rules are thus simpler and more natural than the \PH{} rules 
\eqref{leaprule-P} and \eqref{mcrule-P}
in \refS{Smonths}. They also follow the original \KT, \cf{} \refApp{AKT}.
In particular, if we extend the calendar backwards, there was a leap
month at 
the \Kc{} epoch \nag{} (\Caitra) 806, in agreement with the \KT.)
\begin{remark}
This means that for the \TS\ version,
the intercalation index calculated in this paper 
always agree with the traditional definition without further correction; 
cf.\ \refR{Rix} for the \PH{} version.  
\end{remark}

By \eqref{jeppe}, 
\eqref{leaprule-T} agrees with
the general rule \eqref{leaprule-gb0} 
if $\gb$ is defined such that
\begin{equation}\label{gbx-T}
  \gb+  \gbx \equiv 6 \bmod{65}.
\end{equation}
As we will see in \refApp{ASTSdp},
the values are
$\gb=142$ (E1732) and 187 (E1852), in accordance with
\eqref{gbx-T} and the values for $\gbx$ given above.

By \eqref{gam} and \eqref{gamx}, 
the rules
\eqref{lygam}, \eqref{lygamx} and the formula \eqref{lys} hold with
$\gamma =55$ and $\gamx =20$.

\begin{remark}\label{RTSkarana}
\TS{} almanacs give essentially the same information as almanacs for the
\PH{} version, see \refS{Sfurther}, but one difference is that 
they by tradition give the true solar longitude for each day calculated by the
\karana{}  calculation, see \refApp{AKT}, instead of the calculation above.
(Cf.\ \ref{q-zla-skar-K} in \refS{Sfurther}.)

In some \TS{} almanacs, the solar equation
from the \karana{} solar longitude calculation has also been used to
calculate the (\siddhanta) $true\_date$ in \eqref{truedate}. 
\xfootnote{
This seems only to have been a time saving device, as
usually the \siddhanta{} and \karana{} calculations are kept seperate.
}
In other words, the values \eqref{s1KT} and \eqref{s0K} are used in the
calculations instead of the values for $s_1$ and $s_0$ given above.
See \cite[Open source Tsurphu calendar software]{kalacakra}.
This will lead to a slightly different $true\_date$, and occasionally
a different repeated or skipped day.
\xfootnote{
Calculations similar to the ones in \refR{R30} suggest that this will happen
about 5 times per year, and that the New Year will differ by a day about 2
times per century.
One example is day 13, month 6, 2013, for which the two versions yield
20 and 21 July.
}
However, the traditional version seems to be to use the values above for the
calculation of the true date, but to also calculate separately
the \karana{} version of
$true\_sun$ and publish it
[Henning, personal communication].
\end{remark}


\subsection{Mongolia}\label{AMongo}
The Mongolian calendar is a version of the Tibetan. 
It became the official calendar in Mongolia in 1911 when 
Mongolia declared independence from China 
after the Chinese revolution that overthrew the last Chinese emperor;
however, it was replaced in the 1920s (under Communist rule) by
the Gregorian calendar (officially in 1948)
and the authorities even tried to abolish the traditional celebrations of
the Mongolian New Year. 
The calendar has had a revival together with other traditions
after the end of Communist rule in the 1990s, although
the Gregorian calendar remains the official calendar and is used for
everyday civil use;
\xfootnote{\label{fmongol}%
The first democratic constitution came into force ``from the horse hour of
the auspicious yellow horse day of the black tiger first spring month of the
water monkey year of the seventeenth 60-year cycle''. 
[noon 9/1 = 12 February 1992]
\cite[p.~240]{mongolian}
} 
in particular,  the Mongolian New Year (\emph{Tsagaan Sar}) 
is again celebrated as a major national holiday. See
\cite{mongolian}, \cite{MongoliansWelcome},
\cite[4.1.3]{legalinfo.mn}. 
(I guess that otherwise the Mongolian calendar
is mainly used for religious purposes and astrology.)
In 2012, a second public holiday was declared that is calculated by the
Mongolian calendar (the other holidays have fixed dates in the Gregorian
calendar \cite{legalinfo.mn}), \viz{} Genghis Khan's birthday, 
the first day in the first winter month (month 10)
\xfootnote{
His actual date of birth is not known, and there are different dates in
different sources; for example Terbish \cite{Terbish-GenghisKhan} believes in
day 16 in month 4, 1162. (However, this was 500 years before the present
vesion of the calendar was introduced, so the calculation of the
corresponding Julian date (1 May) in \cite{Terbish-GenghisKhan} seems
uncertain.)}
\cite{Q++-chinggis},
\cite{embassy-holidays},
\cite[4.1.8]{legalinfo.mn}.

According to  \cite{Berzin}, 
the Buryats and Tuvinians of
Siberia (in Russia) follow the New Genden version
(\ie, the version described here),
while
the Kalmyk Mongols in Russia follow
the \PH{} version.
(Inner Mongolia, in China,
uses instead the Chinese ``yellow system'', see \refS{Ayellow}.)

The information below is mainly based on 
\citet{Berzin}, \citet[Epoch data]{kalacakra} and \citet{Salmi}, but I have no
reliable Mongolian sources, and not even a printed calendar as an example.
An example of its use is in the daily horoscope at \cite{news.mn},
but as noted by \cite{Salmi}, this contains some obvious errors and is thus not
reliable. 
\xfootnote{In, for example, May, June and September 2012, some Mongolian
dates are out of order: others are missing.
Nevertheless, most of the skipped or repeated dates in \cite{news.mn} 
for March 2011 -- August 2013 agree with the rules below, as does the
the leap month 6 in 2011.
(Note that the skipped and repeated days are sensitive to also
  small changes in the constants. 
For example, they usually differ between the New Genden version and the \PH{}
version, see the examples in \refT{T4+-}.)
The discrepancies that exist may be due to further errors in their dates. It
is also possible, although less likely, that different versions of the
Mongolian calendar are used. 
}
A list of Gregorian dates of the  Mongolian New Years 
1896--2008 is given in \cite{olloo.mn}, and the same list extended to 2013
in
\cite{mongolnews}.
\xfootnote{These dates seem to be calculated by L. Terbish, who seems to be
the present expert on the calendar in Mongolia (see \eg{} \cite{mongolnews}).
The dates all agree with the rules below.
Note also that the list of New Years immediately yields the leap years. 
Since the leap months are regularly spaced, see \refS{Smonths},
the sequence of leap years for 65 consecutive years uniquely determines the
leap months for all years, and they agree with the leap year rule below.
}\,
\xfootnote{However, there are also other, contradictory lists, for example
  in \cite{mongolian} and on various unreliable web sites.
When this is written, Mongolian Wikipedia \cite[Tsagaan Sar]{mn.wikipedia}
gives two lists, one 1989--2013 (partial) and one 2000--2099; of the 13
common years, only 4 agree. (The first list is attributed to Terbish, and
agrees with the formulas here, except for a possible typo. The second list
agrees with neither the New Genden version,
the \PH{} version, nor the Chinese calendar.)
}
See \cite{Salmi} for many further references.

The Mongolian calendar follows,
see  \cite{Berzin},  \cite{Salmi}, \cite{MongoliansWelcome} and
\cite{mongolnews}, 
a version of the Tibetan calendar known as New Genden 
(Mongolian \emph{T\"ogs buyant}, ``Very virtuous'' \cite{Salmi})
which was created by Sumpa Khenpo \YP{}  
\xfootnote{
An prominent 18th century Tibetan monk of Mongolian origin.
\cite{treasury}
}
(\tibx{sum-pa  mkhan-po ye-shes dpal-'byor};  
Mongolian S\"umbe Khamba Ishbaljir) 
in 1786
\xfootnote{
According to  \cite{Berzin}, who comments that the starting point
is the 40th year of the 60 year cycle and 
claims that "Because of this
difference, the Mongolian calendar works out to be unique." 
However, as said in \refR{Repoch}, the choice of epoch in itself has no
importance. Furthermore,
\YP{} used the epoch 1747
(the beginning of the corresponding 60 year cycle)
 \cite[Epoch data]{kalacakra}.
}.

The constants $m_1$, $s_1$ and $a_1$ for mean motions are the same 
in the New Genden version
as in the
\PH{} version, see \refS{Sastro}, 
and thus so are $m_2$, $s_2$ and $a_2$, see \refR{R30},
while the epoch values $m_0$, $s_0$ and $a_0$ are different and given by,
see  \cite[Epoch data]{kalacakra}:
\begin{align}
\JD&=2359237 \quad \text{(Sunday, 9 April 1747 (Greg.))}
\notag
\\
m_0
&
=  1;55,13,3,31,394 \rr{60,60,6,67,707}
+ 2359236
\notag
\\& 
 =2359237+\frac{2603}{2828},
\label{m0genden}
\\
  s_0&
= 26,39,51,0,18 \rr{27,60,60,6,67} =\frac{397}{402}, \label{s0genden}
\\
a_0&
=  24,22\rr{28,126}  = \frac{1523}{1764}.
\label{a0genden}
\end{align}
I have here, as in \eqref{m0}, modified $m_0$ to yield the JD, see \refR{R7}.

\begin{remark}
In traditional calculations, one uses
a different sequence of radices 
than the standard one, but the mean values can be expressed
exactly also with the standard \PH{} radices;
the constant $m_0$ in \eqref{m0genden} above equals
$1;55,13,3,333 \rr{60,60,6,707}$.
Similarly, $m_1$ is given by \YP{} as
$1;31,50,0,45,345 \rr{60,60,6,67,707}$, which equals the standard value
$1;31,50,0,480 \rr{60,60,6,707}$; as always these traditional values are
interpreted modulo 7, see \refR{R7},
and we add 28 and use the value in \eqref{m1}.
\end{remark}

\begin{remark}
For comparison, 
  the constants for the \PH{} version for this epoch, which thus give the
  standard \PH{} calendar, can be calculated to be
\begin{align}
m_0  &
 = 
1;52,41,2,524  \rr{60,60,6,707} 
+ 2359236
 =2359237+\frac{4967}{5656},
\\
s_0 &
 =
26,9,37,3,45 \rr{27,60,60,6,67} =\frac{779}{804}, 
\\
a_0 &
 = 
24,19\rr{28,126}  =  \frac{3043}{3528}.
\end{align}
\end{remark}

To find the leap months, the true month is calculated by 
\eqref{tm0} with a constant $\gbx =10$
(the epoch value of the intercalation index). However, the Phugpa rule
\eqref{leaprule-P} is modified to
\xfootnote{The reason seems to be that counting backwards to the \Kc{}
  epoch \nag{} (\Caitra) 806 yields the intercalation index 46 and thus a leap
  month,
agreeing with the \KT{} (as in the \TS{} version but not in \PH).  
As a consequence, the Mongolian version has the same leap months as the
\TS{} version (although one may start a day before the other),
see \refT{T4Leap}.
}
\begin{equation}\label{leaprule-G}
\hskip-1em
\vbox {\narrower\narrower\narrower\noindent\em
A leap month is inserted when the intercalation index\\
$ix = 46\text{ or\/ }47$.
}
\hskip-3em
\end{equation}
(There is also a corresponding change of \eqref{mcrule-P}, with 48 replaced
by 46. Similarly, \refR{Rix} applies again, with 49 replaced by 47.)

By \eqref{jeppe}, 
\eqref{leaprule-G} agrees with
the general rule \eqref{leaprule-gb0} 
if $\gb$ is defined such that
\begin{equation}\label{gbx-G}
  \gb+  \gbx \equiv 52 \bmod{65}.
\end{equation}
For the epoch above, with $\gbx=10$,
we take $\gb=172\equiv42\pmod{65}$, see
\refApp{ASMdp}.

By \eqref{gam} and \eqref{gamx}, 
the rules
\eqref{lygam}, \eqref{lygamx} and the formula \eqref{lys} hold with
$\gamma =55$ and $\gamx =20$,
 \cf{} \cite{Salmi}.
(The same constants as for the \TS{} version in \refApp{ATS}, 
since as said above, these two versions have the same leap
months.)

The years in the Mongolian calendar
are named by Element + Animal as described in \refS{Syear}, but
often the element is replaced by the corresponding colour according to the
correspondence in \refT{T5} \cite{mongolian}.  
\xfootnote{
The Mongolians use the colours in
parenthesis in \refT{T5} \cite[p.~155]{mongolian}.
(I suspect that this is a difference in translation of the colours more
than an actual difference in colour.
}\,
\xfootnote{ 
According to \cite[p.~155]{mongolian}, 
it is common to use the element in male years and the colour in female
years,
which would give the cycle Wood, Blue, Fire, Red, Earth, Yellow, Iron,
White, Water, Black, see \refT{T5}.
However, I have not seen this confirmed in other sources.
For example, the lists (in Mongolian) of New Years in \cite{olloo.mn} 
and \cite{mongolnews} use
Element + Animal for all years.
}
The Tibetan 60-year \tibx{rab byung} cycle,
see \refS{Syear},
is recognized in Mongolia too
(called \tibx{jaran} in Mongolian), and the cycles are numbered as in Tibet
(with the first starting in 1027), 
see
\cite{MongoliansWelcome}, \cite{olloo.mn} and
\cite{mongolnews} (both with formulas similar to
\eqref{indian1}--\eqref{indian3}), 
and the example in \refF{fmongol}. 

The months are not numbered, but are named. 
\xfootnote{
Numbers are used for Gregorian months. 
See \cite[p.~106]{mongolian} and, for examples, \cite{legalinfo.mn}.
}
Two different systems are used.
(Both have been used earlier in Tibet, see \refS{SSmonths}.)
One method is by Animal, or 
Colour (or perhaps Element?) + Animal, 
by 
the same rules as described for the \TS{} version 
in \refApp{ASattributes-months} below
(which are the rules of the Chinese calendar); thus month 1 is Tiger, etc.
The months are also named as beginning, middle and end of each of the four
seasons spring, summer, autumn, winter, with month 1 beginning of spring
(as in the Chinese calendar, and in \refS{SSmonths}\ref{month-season-animal})
\cite[pp.~155, 241--242]{mongolian}.
See \refF{fmongol} for an example where both methods are used together.

The calendar days can also be named by Colour (or Element?) + Animal,
in a simple 60-day cycle, see \refApp{ASattributes-day} and the example
 in \refF{fmongol} \cite[p.~242]{mongolian}.

\newcommand\Dzo{Dzongkha}

\subsection{Bhutan}\label{ABhutan}
Bhutan uses a version of the Tibetan calendar as an official calendar;
see \eg{} the government's web page \cite{Bhutan} for an example, where 
a calendar with both Bhutanese and Gregorian dates is shown.
In official documents such as acts, both the Bhutanese and Gregorian dates
are given (in both the \Dzo{} and English versions), see many examples at
\cite[Acts]{BhutanNA}. Of the public holidays, some have fixed dates in
the Bhutanese calendar 
(New Year and some dates connected to Buddha and Buddhism);
some have fixed dates in the Gregorian calendar (\eg{} the National Day
and the King's Birthday); finally, the Winter Solstice is a holiday, but it is
calculated according to the Bhutanese calendar (see below); it usually
occurs at the Gregorian date 2 January. (See 
\cite[Open source Bhutanese calendar software]
{kalacakra} and \cite[Public holidays]{Bhutan} for a complete list.)

The version was described by Lhawang Lodro 
\tib{lha-dbang blo-gros}
in the 18th century
\cite[Open source Bhutanese calendar software]
{kalacakra}, but is said
to be older. 

An important difference from the other main versions of the Tibetan
calendar is that a leap month is given the number of the \emph{preceding} month
insted of the next month,
see \cite[Open source Bhutanese calendar software]{kalacakra}.
(This is the system in the Chinese calendar, but not in Indian calendars
\cite{CC}, see \refR{RChina}; it seems to have been the original \KT{}
system,
see \refApp{AKT}.)

A unique feature of the Bhutanese calendar is that its day of week differs by
one day from Tibet (and the rest of the world, see \refF{fplanets}).
The names of the days of week are the same as in Tibetan
\xfootnote{
\Dzo{} (the national language of Bhutan) is, depending on one's view,
a dialect of Tibetan or a language closely related to Tibetan;
it is written with Tibetan script.
},
but day 0 (Saturday) is \tibx{spen pa} in Tibet but \tibx{nyi ma} 
(really meaning sun) in Bhutan,
day 1 (Sunday) is \tibx{nyi ma} in Tibet but \tibx{zla ba} (really meaning
moon) in Bhutan,
etc., see \refT{T7}.
(This difference is lost when writing dates in English, since of
course the correct English day of week is chosen in both cases, for example
in the 
calendar at \cite{Bhutan}.)
See further \cite[Bhutan calendar problem]{kalacakra}.

The constants $m_1$, $s_1$ and $a_1$ for mean motions are the same 
in the Bhutanese version
as in the \PH{} version (and in the \TS{} and Mongolian versions), see
\refS{Sastro},  
and thus so are $m_2$, $s_2$ and $a_2$, see \refR{R30},
while the epoch values are different and given by,
see Henning \cite[Epoch data]{kalacakra}:
\begin{align}
\JD&=2361807 \quad \text{(Monday, 22 April 1754 (Greg.))}
\notag
\\
m_0
&
= 2;4,24,552 \rr{60,60,707}
+ 2361805
 =2361807+\frac{52}{707},
\label{m0bh}
\\
  s_0&
= 0,24,10,50 \rr{27,60,60,67} =\frac{1}{67}, \label{s0bh}
\\
a_0&
=  3,30\rr{28,126}  = \frac{17}{147}.
\label{a0bh}
\end{align}
Again,  $m_0$ is here modified to yield the JD, see \refR{R7}.

\begin{remark}
The traditional calculations use
different sequences of radices than the standard ones (with fewer radices),
see 
\eqref{m0bh}--\eqref{s0bh},  but the standard values can be expressed
exactly also with these radices.
The constant $m_1$ is in this tradition given by 
$1;31,50,80 \rr{60,60,707}$, which equals the standard value \eqref{m1}
after our usual addition of 28;
similarly, $s_1$ is given by 
$2,10,58,14 \rr{27,60,60,67}=65/804$ which equals \eqref{s1}.
\cite[Epoch data]{kalacakra}. 	
\end{remark}

To find the leap months, the true month is calculated by 
\eqref{tm0} with a constant $\gbx =2$
(the epoch value of the intercalation index)
\cite[Epoch data]{kalacakra}. 	
Furthermore,
 recalling that in Bhutan 
a leap month gets the number of the preceding month,
the Phugpa rule
\eqref{leaprule-P} is modified to 
\begin{equation}\label{leaprule-B2}
\hskip-1em
\vbox {\narrower\narrower\narrower\noindent\em
A month is a leap month if and only if
its intercalation index
$ix = 59\text{ or\/ }60$.
}
\hskip-3em
\end{equation}
Here have to use the correct (traditional) intercalation index, see
\refR{Rix},
which is \eqref{ix} increased by 2 for the leap month (and all later months
until $ix$ passes $65$).
A computationally simpler, but
more clumsy, formulation where our simplified
definition \eqref{ix} can be used directly is:
\begin{equation}\label{leaprule-B}
\hskip-1em
\vbox {\narrower\narrower\narrower\noindent\em
Regular month $M$ in year $Y$ is followed by a leap month if and only if
its intercalation index 
$ix = 57\text{ or\/ }58$.
}
\hskip-3em
\end{equation}

There is a corresponding change of \eqref{mcrule-P}, with 48 replaced
by 59, and the true month rounded up for a leap month. 

Comparing \eqref{leaprule-B} with
\eqref{jeppe} (and taking into account that now a leap month $M$ comes after
the regular month $M$), or better by \eqref{leaprule-B2} and \eqref{jeppe+}
below, 
\eqref{leaprule-B2} and
\eqref{leaprule-B} agree with
the general rule \eqref{leaprule-gb0} 
if $\gb$ is defined such that
\begin{equation}\label{gbx-B}
  \gb+  \gbx \equiv 63 \bmod{65}.
\end{equation}
For the epoch above, with $\gbx=2$,
we take $\gb=191\equiv61\pmod{65}$, see
\refApp{ASBdp}.

By \eqref{gam} and \eqref{gamx}, 
\eqref{lygam} and \eqref{lygamx} hold with
$\gamma =12$ and $\gamx =28$.

\smallskip
The winter solstice is, as said above, a holiday in Bhutan. 
It is defined as the
day the mean solar longitude \eqref{meansun} reaches 
$250\grad$, 
which at present occurs 
on 2 January,
see \cite[Open source Bhutanese calendar software]{kalacakra}.
\xfootnote{The correct astronomical definition is $270\grad$, for the
  (astronomical) true longitude; the correct Gregorian date is 
(at present)
21 or 22 December.
}
(The date drifts slowly, by almost 3 days per century, see \refS{Smean}.
It will be 3 January for the first time in 2020. It coincided with the
astronomical winter solstice about 400 years ago, in the early 17th century,
which was before the present version of the Bhutanese calendar was described.)

The months in Bhutan are numbered, as in Tibet.
(Names exist as in Tibet, see \refS{SSmonths}, but are usually not used.)
The years are 
named by Element + Gender + Animal, see \refS{Syear}. 
(See the examples in 
\cite[Acts]{BhutanNA}, where the Bhutanese dates are given in this way,
while the Gregorian dates are given with the number of the year.)


\subsection{\KT\  (\Karana)}\label{AKT}
The original \KT{} calculations are explained by 
\citet[Chapter 5]{Schuh} and 
\citet[Chapter V]{Henning}, \cite{Henning-Karana}; they are
(as far as I know) not used to produce calendars today; 
however, 
complete Tibetan almanacs usually give (by tradition, and for no other
obvious reason) 
some values calculated by this version in addition
to values calculated by one of the versions above, see \refS{Sfurther}.
The \KT{} version is called \emph{\karana} \tib{byed rtsis}, while the \PH{}
and other similar versions are called \emph{\siddhanta} \tib{grub rtsis}.

The \KT{} (\karana) 
version uses  different values of the mean motions $m_1$ and $s_1$
than the versions described above but the same $a_1$
(see also \cite[Epoch data]{kalacakra}):
\begin{align}
m_1  &
 = 
29;31,50  \rr{60,60} 
 =\frac{10631}{360}
=29+\frac{191}{360},
\label{m1KT}
\\
s_1 &
=
2,10,58,2,10 \rr{27,60,60,6,13}
=\frac{1277}{15795}, 
\label{s1KT}
\\
a_1 &
 = 
2,1\rr{28,126}  =  \frac{253}{3528};
\intertext{the epoch values are} 
\JD&=2015531 \quad \text{(Monday, 23 March 806 (Julian))}
\\
m_0  &
 = 
2;30,0  \rr{60,60} 
+ 2015529
 =2015531+\frac{1}{2},
\\
s_0 &
 =
26,58 \rr{27,60} =\frac{809}{810}, \label{s0K}
\\
a_0 &
 = 
5,112\rr{28,126}  =  \frac{53}{252}.
\intertext{%
Modern \karana{} calculations (\eg{} in the almanac \cite{tib2013})
use (at least in the cases I know) the natural values 
$m_1/30$ and $s_1/30$
for the daily increments $m_2$ and $s_2$, while
$a_2=1/28$ is given by \eqref{a2} as for the versions above.
However, the \KT{} gives instead the
following simplified (rounded) values for $m_2$ and $s_2$,
\cf{} \refR{R30}, which thus not are used today:}
m_2  &
 = 
0;59  \rr{60} 
 =\frac{59}{60},
\\
s_2 &
 =
0,4,20 \rr{27,60,60} = \frac{13}{4860}, 
\intertext{%
Furthermore, 
the \KT{} gives also 
another computational simplification (perhaps suggested as an alternative?),
also not used today 
(as far as I know, and certainly not in the almanac \cite{tib2013}),
see
\cite[p.~232--233]{Henning},
\cite[Chapter 2, v.~30]{Henning-Karana},
where
the values $m_1,s_1,a_1$ above are used only for the first month of a
calendar calculation (i.e, month 3, \tibx{nag pa}), and for each following
month one 
adds instead the simpler\footnotemark}
m_1'  &
 = 
29;32,0  \rr{60,60} 
 =\frac{443}{15}
=29+\frac{8}{15},
\label{m1'KT}
\\
s_1' &
=
2,11 \rr{27,60}
=\frac{131}{1620}, 
\label{s1'KT}
\\
a_1' &
 = 
2,0\rr{28,126}  =  \frac{1}{14}.
\end{align}
\footnotetext{It seems surprising that anyone should bother about these
  minor simplifications, 
once each for each month, which are very minor compared to the mass of
calculations for every day.}

The \KT{} calculates the true month by \eqref{tm0} with an epoch value
$\gbx=0$. The true month is simply rounded down, 
as in the later \TS{} version, see \refApp{ATS}. 
This would correspond to the
intercalation rule \eqref{leaprule-T}
if a leap month is given the same
number as the following month, as in Indian and later Tibetan calendars.
However, it seems that the \Kc{} system instead 
gave a leap month the same name as the month preceding it, 
as in the Bhutanese calendar in \refS{ABhutan} 
(and the Chinese calendar \cite{CC}),  
see \cite[Published calendar explanation]{kalacakra}.
In this case, the leap month is inserted before the month with
intercalation index $0$ or $1$; in the traditional formulation,
these indices are repeated (see \refR{Rix}) and thus the leap month
gets the same intercalation index, so
the rule \eqref{leaprule-T} nevertheless holds.
Cf.\ the rules \eqref{leaprule-B2} and \eqref{leaprule-B} for the  Bhutanese
version. 
\xfootnote{
The \Karana{} calculations in modern almanacs do not give any names or
numbers, or years, for the months, as far as I know,
so the question of leap months is irrelevant to them. 
However, in \eg{} \cite{tib2013},
the \tmc{} and intercalation index are given for each month, for both  \PH{}
and \karana{} calculations.
This almanac includes a month with \karana{}
intercalation index given as 65,
in accordance with the rules presented here; this month would thus be
a leap month (a second \tibx{M\=agha} = \tibx{mChu}, or month 1, 2014 if
we anachronistically use the same numbering as  for the other
versions, although this was introduced long after the \KT, 
see \refS{SSmonths}).
}

More precisely, this means that there is a leap month $(Y,M)$ when the
regular month $(Y,M)$ has intercalation index \eqref{ix} 63 or 64,
\ie, its \tm{} is $n+63/65$ or $n+64/65$ where $n$ is the  \tmc.
The month $(Y,M+1)$ then has a \tm{} that is
$\xfrac{67}{65}=1+2/65$ larger,  \ie{} $n+2$ or $n+2+1/65$; hence it has
\tmc{} that is $n+2$ and intercalation index $0$ or $1$.
There is thus a gap, which is filled by a leap month which is given the
missing \tmc{} $n+1$, and (conventionally) intercalation index $65$ or 66;
this month is thus leap month $(Y,M)$.
We thus have the rule, \cf{} \eqref{leaprule-B},
\begin{equation}\label{leaprule-K}
\hskip-1em
\vbox {\narrower\narrower\narrower\noindent\em
Regular month $M$ in year $Y$ is followed by a leap month if and only if
its intercalation index 
$ix = 63\text{ or\/ }64$.
}
\hskip-3em
\end{equation}
Furthermore, we can say, \cf{} \eqref{leaprule-B},
\begin{equation}\label{leaprule-K2}
\hskip-1em
\vbox {\narrower\narrower\narrower\noindent\em
A month is a leap month if and only if
its intercalation index
$ix = 65\text{ or\/ }66$.
}
\hskip-3em
\end{equation}
This, however, is a tautology since we use the exceptional values 65 and 66
only for leap months. We cannot replace them by 0 and 1 in \eqref{leaprule-K2},
since the following (regular) month has intercalation index 0 or 1.

Summarizing, 
in all cases, the  \tmc{} is given by \eqref{tm0} rounded down, and
increased by 1 in the case of a leap month; the intercalation index is given
by \eqref{ix}, increased by 2 (to 65 or 66) for a leap month.
(Cf.\ \refR{Rix}.)

For the epoch 806 above, with $\gbx=0$, it follows from \eqref{ix} and the
discussion above that \eqref{leaprule-K} is equivalent to
\eqref{leaprule-gb0} with $\gb\equiv 4\pmod{65}$, see also \eqref{jeppe+} below.

By \eqref{gam} and \eqref{gamx}, 
\eqref{lygam} and \eqref{lygamx} hold with
$\gamma =28$ and $\gamx =22$.

\smallskip
The value \eqref{m1KT} of $m_1$ (the mean length of the month) is
$\approx 29.530556$ days, about 2.7 seconds shorter than the mean
length \eqref{meanmonth} for the standard versions. The \Kc{} value 
differs from the modern astronomical value by about $2.9$ seconds, and is thus
less exact than the value \eqref{meanmonth}, but the difference amounts to
less than 1 hour every century.

The mean length of the year is, \cf{} \eqref{meanyear},
\begin{equation}
\label{meanyearK}
\frac{m_1}{s_1}
=\frac{3731481}{10216}
\approx 365.258516
\text{ days},
\end{equation}
which is better than the value in \eqref{meanyear} for the later versions;
it is $0.01633$ days longer than 
the astronomical value of the tropical year ($365.24219$ days)
and only $0.002156$ days longer than the 
sidereal year ($365.25636$ days),
see \refS{Smean} and \cite[Table 15.3]{AA}.
(Nevertheless, the year is intended to be a tropical year
\cite[p.~298]{Henning}.)

\subsection{Further historical versions}
Several further historical versions of the Tibetan calendar are described by
\citet{Schuh}. (I do not know to which extent these were actually used.)

A proposed (but never adopted) version by Zhonnu Pal \tib{gzhun-nu dpal}
from 1443  
is described by 
\citet[pp.~307--321]{Henning} and
\cite[Error correction system]{kalacakra}.

\subsection{Sherab Ling}\label{ASL}
A modern attempt at a reformed Tibetan calendar developed by Kojo
Tsewang Mangyal (Tsenam) at the Sherab Ling monastery in Bir, India,
is described in \cite[pp.~342--345]{Henning}, 
see also \cite[A reformed Tibetan calendar]{kalacakra}
and \cite[Epoch data]{kalacakra}.
It uses the value of $m_1$ in \refS{Sastro}, but uses 
\begin{equation}
  s_1=2,10,58,2,564,5546
\rr{27,60,60,6,707,6811}
=\frac{3114525}{38523016}.
\end{equation}
This deviation from the traditional value, 
and from the equivalent fundamental relation \eqref{6765},
means that leap months no
longer will be regularly spaced in the traditional pattern with 2 leap
months every 65 regular months.
(However, the difference is small; the average number of leap months for 65
regular months is $1.9978$, about $0.1\%$ less than in the standard
calendars, which means that on the average a leap month is delayed a month 
after about 75 years.)

\subsection{Sarnath}
A version of the Sherab Ling calendar (see \ref{ASL})
is published by the Jyotish
Department of    	
the Central University of Tibetan Studies 
at Sarnath, India,
see 
\cite[p.~346]{Henning},
\cite[A reformed Tibetan calendar]{kalacakra}.
This version differs from all other Tibetan calendars that I have
heard of in that the months begin at full moon (as in Indian calendars
in northern India \cite{CC}).

\subsection{Sherpa}
The Sherpas (a minority in Nepal) have a Tibetan calendar
called the Khunu almanac. 
(The official Nepalese calendar is of the Indian type.)
The Sherpa calendar for 2002 I found on the web \cite{sherpa} 
is identical to the \PH{} version, and this seems to be the general rule.

\subsection{Drikung Kagyu}
According to \cite{Berzin}, the Drikung Kagyu tradition follows a system
that combines the \TS{} and \PH{} traditions. 

\subsection{Yellow calculations}\label{Ayellow} (Chinese-style)
According to \cite{Berzin},
the yellow system is like the Chinese in that it has no repeated or
skipped days. Months have 29 or 30 days, numbered consecutively "and
determined according to several traditions of calculation". The way it
adds leap months is similar to, but not equivalent to, the
Chinese. Unlike the Chinese calendar, it uses the basic calculations
from \KT. 
This seems to be a version of the Chinese calendar 
(or astrology) rather than the Tibetan; it might also be the Chinese
calendar (without modification), but then some of the statements just made
are incorrect. 
See also \cite{MongoliansWelcome} and \cite{whenMongolia}.
Inner Mongolia (in China) follows the yellow system \cite{Berzin}.

\subsection{An astronomical version}
Henning \cite{kalacakra} has constructed 
a reformed Tibetan calendar combining the principles of the \KT{}
with modern astronomical calculations.


\subsection{A comparison} \label{Acomp}
We give some examples of the (small) differences between the 
four different versions of the Tibetan calendar in
Appendices \ref{APH}--\ref{ABhutan}.

The epoch values given above for the different versions use different
years. To enable an easy comparison, \refT{T4epok} gives the epoch values
for the four versions calculated for the same epoch, here chosen as
the \KT{} epoch at the beginning of \nag{} (\Caitra) 806.
(In all four versions, the epoch, \ie, the mean new moon,
is at JD 2015531, Monday 23 March 806 (Julian), \cf{}
\refR{Repoch2}.) 
\xfootnote{
This happens to be the last day of the preceding month in all four versions;
this is day 30 in month 2 for the \PH{} and Bhutanese versions, and day 30
in leap month 3 in the \TS{} and Mongolian versions.
}
Since all four versions have the same mean motions $m_1,s_1,a_1$, the 
differences between them for
mean dates, solar longitudes and lunar anomalies will be constant. For
example, the \TS{} mean solar longitude is always
$0.018261-0.004975=0.013286=4.78\grad$ larger than the \PH{} value,
and the mean dates differ by 0.046, about $1/20$.
The differences in true date and true solar longitude will be varying because of
the corrections in \eqref{truedate} and \eqref{truesun}, but the average
difference is the same as for the mean values. 
Thus the
\PH{} and \TS{} dates will differ for, on the average, about one day in 20,
\ie, typically one or two days in a month. (See \refT{T4+-} for some
examples, with 0--4 days differing each month.)

\begin{table}[!htpb]
\begin{tabular}{lllll}
 & \PH & \TS & Mongolia & Bhutan\\
\hline 
$m_0$ & $2; 22, 34, 2, 518$ & $2; 25, 20, 2, 352$ & $2; 25, 6, 3, 327$ 
  & $2; 24, 37, 5, 431$ \\
 & 2.376238 & 2.422338 & 2.418494 & 2.410537 \\
\hline 
$s_0$ & $0, 8, 3, 3, 33$ & $0, 29, 34, 5, 37$ & $0, 38, 17, 0, 6$
 & $0, 28, 12, 3, 15$ \\
 & 0.004975 & 0.018261 & 0.023632 & 0.017413 \\
\hline 
$a_0$ & $5, 98$ & $5, 112$ & $5, 101$ & $6, 22$ \\
 & 0.206349 & 0.210317 & 0.207200 & 0.220522 \\
\hline 
  \end{tabular}
\caption{Epoch data for four versions of the Tibetan calendar
for JD 2015531 (23 March 806), given both in
Tibetan form with radices
$\xpar{60,60,6,707}$,
$\xpar{27, 60, 60, 6, 67}$,
$\xpar{28,126}$
and with 6 decimals.
We have here given $m_0$ in the traditional form modulo 7; to get the result
in JD as in this paper, add 2015529.
}
\label{T4epok}
\end{table}

\begin{remark}
The differences in the correction terms can be estimated as follows:  
Consider two versions whose epoch values $m_0,s_0,a_0$ differ by 
$\gD m_0$, $\gD s_0$, $\gD a_0$. By \eqref{moonequ}, the arguments used
in the table $moon\_tab$ differs by $28\gD a_0$; since the maximum 
derivative (slope) in the table \eqref{moontab}
is 5, the difference in $moon\_equ$ 
satisfies $|\gD moon\_equ|\le 140|\gD a_0|$; this correction is divided by
60 in \eqref{truedate},
so the corrections to the true date differ by at most $\frac73|\gD a_0|$.
Similarly, by \eqref{anosun}, \eqref{sunequ} and \eqref{suntab},
the difference in $sun\_equ$ satisfies
$|\gD sun\_equ|\le 6\cdot12|\gD s_0|$; this correction is divided by
60 in \eqref{truedate} and by $27\cdot60$ in \eqref{truesun},
so the corrections to the true date and true solar longitude differ by at
most $1.2|\gD s_0|$ and $0.045|\gD s_0|$. 
In particular, the difference between the true solar longitudes differ
always by between $0.95$ and $1.05$ times the mean difference $\gD s_0$; for
\PH{} and \TS{} the difference is thus between $4.57\grad$ and $5.00\grad$.

For the date, the mean difference $\gD m_0$
between \PH{} and \TS{} is as said above
$0.046$, and the maximum differences in the correction terms
$\frac1{60}\gD moon\_equ$ and $\frac1{60}\gD sun\_equ$ are
$\frac73\gD a_0=0.009$ and $1.2\gD s_0=0.016$, respectively; 
hence the difference of the true
dates is $0.046\pm0.009\pm0.016$, \ie, between $0.021$ and $0.071$,
with \TS{} always larger. 
\xfootnote{
The two correction terms are periodic, with periods the anomalistic month
\eqref{meanano} and  the year \eqref{meanyear}; 
the maximum values will thus coincide several times each year.
}
As a consequence, the \PH{} version is slightly
behind the \TS, so when the \PH{} and \TS{} dates
of a calendar day differ, the \PH{} date is always larger (by 1).
(See the example of differences in \refT{T4+-}.)

Similarly, the difference in mean date
between the  Mongolian and \TS{} versions is 
only $\gD m_0=0.0038$; the differences in solar longitude and lunar anomaly
are
$\gD s_0=-0.0054$ and
$\gD a_0=0.0031$; hence the true date differs by
$0.0038\pm0.0073\pm0.0064$,
\ie, between $-0.010$ and $0.018$. (So we expect a difference of the
calendar date only a few
days each year, \cf{} \refT{T4+-}.)
\end{remark}


\refT{T4Leap} shows leap months for a range of (Gregorian) years.
Since all four versions have a simple periodic pattern with alternating 32
or 33 regular months between the leap months, the same pattern repeats for
ever.
Note that the \TS{} and Mongolian versions have the same leap months, as said
before. 
We see also that leap months in Bhutan come 2 or 3 months after the \PH{}
leap months (3 or 4 months later if we take into account the different
numbering of Bhutanese leap months), and the \TS{} and Mongolian leap months
come an additional 4 (really 3) months later.

\smallskip

\refT{T4+-} shows all repeated and skipped days during 2012 for the four
versions. (This year is chosen since none of the versions has a leap month.
The versions also all start this year on the same day.)
It is seen that four versions are very similar; often the same day is
repeated or skipped, but it also frequently happens that different
versions differ by a day (meaning that some days get different dates). 
We can see that this year each month has the same
length in all four versions.

\smallskip

\refT{T4Losar} shows the Gregorian dates of New Year 
for a range of (Gregorian) years.
We see that most years, all four calendars coincide. However, some times the
\PH{} and Bhutanese versions, or just the \PH{} version, is a month later
(due to the difference in leap months); sometimes there is also a difference
of a day between two versions.
\xfootnote{
Bhutan actually had New Year 3/3 2003, with day 1 of month 1 repeated
3/3 and 4/3, according to the government web site calendar \cite{Bhutan}
[Henning, personal communication];  
the calculations described here yield 
(as do the ones by Henning 
\cite[Open source Bhutanese calendar software]{kalacakra})
3/3 instead as a repeated day 30 in
the last (leap) month of 2002. I have no explanation for this discrepancy.
}

In particular, we do not see any differences between the \TS{} and Mongolian
New Year in \refT{T4Losar}; a computer search reveals that the last time
these versions had different New Year was 1900 (31/1 vs 1/2) and the next
will be 2161 (26/2 vs 25/2), and only four more differences were found for
the present millenium.


\begin{table}[!htpb]
\begin{tabular}{r r r r r}
 & \PH & \TS  & Mongolia & Bhutan\\
\hline
2000  & 1  & 8  & 8  & 4  \\ 
2001  &    &    &    &    \\ 
2002  & 10  &    &    & 12  \\ 
2003  &    & 4  & 4  &    \\ 
2004  &    &    &    &    \\ 
2005  & 6  &    &    & 9  \\ 
2006  &    & 1  & 1  &    \\ 
2007  &    &    &    &    \\ 
2008  & 3  & 9  & 9  & 5  \\ 
2009  &    &    &    &    \\ 
2010  & 11  &    &    &    \\ 
2011  &    & 6  & 6  & 2  \\ 
2012  &    &    &    &    \\ 
2013  & 8  &    &    & 10  \\ 
2014  &    & 2  & 2  &    \\ 
2015  &    &    &    &    \\ 
2016  & 4  & 11  & 11  & 7  \\ 
2017  &    &    &    &    \\ 
2018  &    &    &    &    \\ 
2019  & 1  & 7  & 7  & 3  \\ 
2020  &    &    &    &    \\ 
  \end{tabular}
\caption{Leap months for four versions of the Tibetan calendar.}
\label{T4Leap}
\end{table}

\begin{table}[!htpb]
\begin{tabular}{r| l|l|l|l}
month & \PH & \TS  & Mongolia & Bhutan\\
\hline
1 & 5, -19 & 4, -20 & 4, -20 & 4, -19 \\
2 & 9, -12, -25, 27 & 8, -13 & 8, -13 & 8, -13 \\
3 & -17 & -17 & -17 & -17 \\
4 & 3, -10 & 2, -11 & 2, -11 & 2, -10 \\
5 & -13, 29 & -14, 28 & -14, 28 & -13, 28 \\
6 & -6 & -6 & -6 & -6 \\
7 & -9, 25 & -9, 25 & -9, 25 & -9, 24 \\
8 & -1 & -2 & -2 & -1 \\
9 & -5, 20, -29 & -6, 19, -29 & -6, 20, -29 & -5, 19, -29 \\
10 &  &  &  &  \\
11 & -3, 13, -27 & -3, 12, -28 & -4, 12, -28 & -3, 12, -27 \\
12 & 17, -21 & 15, -22 & 15, -22 & 15, -21 \\
  \end{tabular}
\caption{Repeated and skipped (marked with -) days in each month
2012 for four versions of the Tibetan
  calendar.} 
\label{T4+-}
\end{table}


\begin{table}[!htpb]
\newcommand\xx{\textbf}
\begin{tabular}{r r r r r}
 & \PH & \TS & Mongolia & Bhutan\\
\hline
2000 & 6/2 &  6/2 &  6/2 &  6/2 \\ 
2001 & 24/2 &  24/2 &  24/2 &  24/2 \\ 
2002 & 13/2 &  13/2 &  13/2 &  13/2 \\ 
2003 & \xx{3/3} &  2/2 &  2/2 &  \xx{4/3} \\ 
2004 & 21/2 &  21/2 &  21/2 &  21/2 \\ 
2005 & 9/2 &  9/2 &  9/2 &  9/2 \\ 
2006 & \xx{28/2} &  \xx{30/1} &  \xx{30/1} &  \xx{28/2} \\ 
2007 & 18/2 &  18/2 &  18/2 &  18/2 \\ 
2008 & \xx{7/2} &  8/2 &  8/2 &  8/2 \\ 
2009 & 25/2 &  25/2 &  25/2 &  25/2 \\ 
2010 & 14/2 &  14/2 &  14/2 &  14/2 \\ 
2011 & \xx{5/3} &  3/2 &  3/2 &  3/2 \\ 
2012 & 22/2 &  22/2 &  22/2 &  22/2 \\ 
2013 & 11/2 &  11/2 &  11/2 &  11/2 \\ 
2014 & \xx{2/3} &  \xx{31/1} &  \xx{31/1} &  \xx{2/3} \\ 
2015 & 19/2 &  19/2 &  19/2 &  19/2 \\ 
2016 & 9/2 &  9/2 &  9/2 &  9/2 \\ 
2017 & 27/2 &  27/2 &  27/2 &  27/2 \\ 
2018 & 16/2 &  16/2 &  16/2 &  16/2 \\ 
2019 & 5/2 &  5/2 &  5/2 &  5/2 \\ 
2020 & 24/2 &  24/2 &  24/2 &  24/2 \\ 
2021 & 12/2 &  12/2 &  12/2 &  12/2 \\ 
2022 & \xx{3/3} &  \xx{2/2} &  \xx{2/2} &  \xx{3/3} \\ 
2023 & 21/2 &  21/2 &  21/2 &  21/2 \\ 
2024 & 10/2 &  10/2 &  10/2 &  10/2 \\ 
2025 & \xx{28/2} &  \xx{1/3} &  \xx{1/3} &  \xx{28/2} \\ 
2026 & 18/2 &  18/2 &  18/2 &  18/2 \\ 
2027 & 7/2 &  7/2 &  7/2 &  7/2 \\ 
2028 & 26/2 &  26/2 &  26/2 &  26/2 \\ 
2029 & 14/2 &  14/2 &  14/2 &  14/2 \\ 
2030 & \xx{5/3} &  3/2 &  3/2 &  3/2 \\ 
  \end{tabular}
\caption{%
Gregorian dates for 
New Year (\tibx{Losar}) for four versions of the Tibetan calendar.
Dates differing from the majority in boldface.}
\label{T4Losar}
\end{table}

\clearpage

\section{The 60 year cycle}\label{A60}

The Indian 60 year cycle is derived from an Indian cycle of 60 names for the
``Jovian years''. 
\xfootnote{
It takes Jupiter almost 12 years to orbit the sun.
A 1/12 of the orbital period can be called a Jovian year, 
and 
the traditional Indian Jovian cycle gives each
Jovian year a name (\emph{samvatsara}) 
from a list of 60 names, so the names repeat with a cycle of 5
revolutions. The solar years are given the same names,
based on the calculated position of
Jupiter at the beginning of the year; hence sometimes (every 85 or 86
years) a name is skipped (expunged) from the list. In southern India
this has from the 9th century been simplified to a simple 60 year
cycle of names, 
and the same is done in Tibet. \cite{CC}, \cite[p.~143f]{Henning}.
}
\refT{T60} 
(taken from \cite{Henning})
gives the full list of names, in Tibetan and Sanskrit,
together with the corresponding 
year in the Chinese 60 year cycle and the
Gregorian years in the last and current cycle.
(Somewhat different transliterations of the Sanskrit names are given in
\cite{CC}.)  

See also \refT{Tlosar}, which gives the Gregorian dates of  New Year for
the years in the last and current cycles.


\begin{longtable}
{r r l l l r r }
\rlap{year}\phantom{0} && 
 element--animal & Tibetan & Sanskrit & & \\
\hline
1 & 4 & Fire--Rabbit & rab byung & prabhava & 1927 & 1987\\ 
2 & 5 & Earth--Dragon & rnam byung & vibhava & 1928 & 1988\\ 
3 & 6 & Earth--Snake & dkar po & suklata & 1929 & 1989\\ 
4 & 7 & Iron--Horse & rab myos & pramadi & 1930 & 1990\\ 
5 & 8 & Iron--Sheep & skyes bdag & prajapati & 1931 & 1991\\ 
6 & 9 & Water--Monkey & anggi ra & ankira & 1932 & 1992\\ 
7 & 10 & Water--Bird & dpal gdong & srimukha & 1933 & 1993\\ 
8 & 11 & Wood--Dog & dngos po & bhava & 1934 & 1994\\ 
9 & 12 & Wood--Pig & na tshod ldan & yuvika & 1935 & 1995\\ 
10 & 13 & Fire--Mouse & 'dzin byed & dhritu & 1936 & 1996\\ 
11 & 14 & Fire--Ox & dbang phyug & isvara & 1937 & 1997\\ 
12 & 15 & Earth--Tiger & 'bru mang po & vahudhvanya & 1938 & 1998\\ 
13 & 16 & Earth--Rabbit & myos ldan & pramadi & 1939 & 1999\\ 
14 & 17 & Iron--Dragon & rnam gnon & vikrama & 1940 & 2000\\ 
15 & 18 & Iron--Snake & khyu mchog & brisabha & 1941 & 2001\\ 
16 & 19 & Water--Horse & sna tshogs & citra & 1942 & 2002\\ 
17 & 20 & Water--Sheep & nyi ma & bhanu & 1943 & 2003\\ 
18 & 21 & Wood--Monkey & nyi sgrol byed & bhanutara & 1944 & 2004\\ 
19 & 22 & Wood--Bird & sa skyong & virthapa & 1945 & 2005\\ 
20 & 23 & Fire--Dog & mi zad & aksaya & 1946 & 2006\\ 
21 & 24 & Fire--Pig & thams cad 'dul & sarvajit & 1947 & 2007\\ 
22 & 25 & Earth--Mouse & kun 'dzin & sarvadhari & 1948 & 2008\\ 
23 & 26 & Earth--Ox & 'gal ba & virodhi & 1949 & 2009\\ 
24 & 27 & Iron--Tiger & rnam 'gyur & vikrita & 1950 & 2010\\ 
25 & 28 & Iron--Rabbit & bong bu & khara & 1951 & 2011\\ 
26 & 29 & Water--Dragon & dga' ba & nanda & 1952 & 2012\\ 
27 & 30 & Water--Snake & rnam rgyal & vijaya & 1953 & 2013\\ 
28 & 31 & Wood--Horse & rgyal ba & jaya & 1954 & 2014\\ 
29 & 32 & Wood--Sheep & myos byed & mada & 1955 & 2015\\ 
30 & 33 & Fire--Monkey & gdong ngan & durmukha & 1956 & 2016\\ 
31 & 34 & Fire--Bird & gser 'phyang & hemalambha & 1957 & 2017\\ 
32 & 35 & Earth--Dog & rnam 'phyang & vilambhi & 1958 & 2018\\ 
33 & 36 & Earth--Pig & sgyur byed & vikari & 1959 & 2019\\ 
34 & 37 & Iron--Mouse & kun ldan & sarvavati & 1960 & 2020\\ 
35 & 38 & Iron--Ox & 'phar ba & slava & 1961 & 2021\\ 
36 & 39 & Water--Tiger & dge byed & subhakrita & 1962 & 2022\\ 
37 & 40 & Water--Rabbit & mdzes byed & sobhana & 1963 & 2023\\ 
38 & 41 & Wood--Dragon & khro mo & krodhi & 1964 & 2024\\ 
39 & 42 & Wood--Snake & sna tshogs dbyig & visvabandhu & 1965 & 2025\\ 
40 & 43 & Fire--Horse & zil gnon & parabhava & 1966 & 2026\\ 
41 & 44 & Fire--Sheep & spre'u & pravamga & 1967 & 2027\\ 
42 & 45 & Earth--Monkey & phur bu & kilaka & 1968 & 2028\\ 
43 & 46 & Earth--Bird & zhi ba & saumya & 1969 & 2029\\ 
44 & 47 & Iron--Dog & thun mong & sadharana & 1970 & 2030\\ 
45 & 48 & Iron--Pig & 'gal byed & virobhakrita & 1971 & 2031\\ 
46 & 49 & Water--Mouse & yongs 'dzin & paradhari & 1972 & 2032\\ 
47 & 50 & Water--Ox & bag med & pramadi & 1973 & 2033\\ 
48 & 51 & Wood--Tiger & kun dga' & ananda & 1974 & 2034\\ 
49 & 52 & Wood--Rabbit & srin bu & raksasa & 1975 & 2035\\ 
50 & 53 & Fire--Dragon & me & anala & 1976 & 2036\\ 
51 & 54 & Fire--Snake & dmar ser can & vingala & 1977 & 2037\\ 
52 & 55 & Earth--Horse & dus kyi pho nya & kaladuti & 1978 & 2038\\ 
53 & 56 & Earth--Sheep & don grub & siddhartha & 1979 & 2039\\ 
54 & 57 & Iron--Monkey & drag po & rudra & 1980 & 2040\\ 
55 & 58 & Iron--Bird & blo ngan & durmati & 1981 & 2041\\ 
56 & 59 & Water--Dog & rnga chen & dundubhi & 1982 & 2042\\ 
57 & 60 & Water--Pig & khrag skyug & rudhirura & 1983 & 2043\\ 
58 & 1 & Wood--Mouse & mig dmar & raktaksi & 1984 & 2044\\ 
59 & 2 & Wood--Ox & khro bo & krodhana & 1985 & 2045\\ 
60 & 3 & Fire--Tiger & zad pa & ksayaka & 1986 & 2046\\ 
\caption{The Chinese and Indian 60 year cycles of names, with the
  names from the Indian cycle both in Tibetan and in Sanskrit.
The first number on each line shows the number in the
Prabhava cycle; the second shows the number in  the
Chinese cycle. The last two numbers show the Gregorian
years in the last and current cycles.}
\label{T60}  
\end{longtable}

\clearpage

\section{Leap months and the mean sun}\label{Aleap}

We give here an explanation of the leap month rules in \refS{Smonths},
based on Tibetan astronomy. 
(This appendix is based on the description in \citet{Henning};
see also \cite[On intercalary months]{kalacakra}.
For a detailed historical discussion of different leap month rules,
see \citet[pp.~107--117]{Schuh};
see further 
\citet{Yamaguchi} (which also contains tables with actually observed
leap months from historical data) and 
\citet[Early epochs]{kalacakra}.
\xfootnote{
In particular, according to 
Schuh \cite[Kalenderrechnung]{tibetenc},
see also \cite{Schuh},
the \dpp{s} described here were introduced (for the \PH{} version)
in 1696. 
However, some similar method to align the year seems to have been used earlier; 
the principle to use the position of the mean sun goes back to early
Indian calendars and
the calculations described 
in \refS{Smonths} are only minor modifications of the ones in
the \KT{} which thus seem to be based on considerations related
to the ones
given here, although possibly in different formulations.
See also the discussion in \refS{ASKdp}.
})

\subsection{General theory}\label{ASdp}
The key is the position of the mean sun, which is a fictitious version of the
sun that travels along the ecliptic with uniform speed.
Its longitude (the mean solar longitude)
at the beginning of
\tm{} $n$ 
was denoted by $\MSL(0,n)$ in \refS{Sastro};
we use here the simplified notation $\MSL(n)$.
It is by \eqref{meansun} given
by the linear formula 
\begin{equation}
  \label{msl}
\MSL(n)=s_1n+s_0,
\end{equation}
where $s_0=\MSL(0)$ is the mean solar longitude at the epoch.
The Tibetan constant $s_1$, the mean motion of the sun per (lunar) month,
 is given by
 \begin{equation}\label{s1q}
s_1=  \frac{65}{804}=\frac{65}{67\cdot12}, 
 \end{equation}
see \eqref{s1}.
(This can also be expressed as 65/67 signs, or  $29\frac{7}{67}\grad$.)

Note that \eqref{s1q} says that the mean sun goes 
around the
ecliptic exactly 65 times in 804 lunar months, \ie,
\begin{equation}
804 \text{ lunar months } = 65 \text{ years} .
\end{equation}
Since $804=67\cdot12$, this is equivalent to the relation \eqref{6765},
and it explains the leap year cycle of 65 years, see \refR{R6567}.

Moreover, the zodiac contains 12 evenly spaced \emph{\dpp s} 
\xfootnote{We use this term from \citet{Henning}.} 
\tib{sgang}, 
one in each sign. Let us denote these (and their
longitudes) by $p_1,\dots,p_{12}$.
The rule for naming months is, see  \cite{Henning}:
\begin{equation}\label{rule-dp}
\hskip-1em
\vbox {\narrower\narrower\narrower\noindent\em
A month where the mean sun passes a
\dpp{} $p_M$ is given number $M$.  
}
\hskip-3em
\end{equation}
Since $s_1<\frac1{12}$, the
spacing of the \dpp{s}, the mean sun can never pass two \dpp{s} in one
month, but sometimes it does not pass any of them; in that case the
month is designated as a leap month, and is given the number of the
next month. In both cases, the number of the month is thus given by
the first \dpp{} $p_M$ that comes after $\MSL(n)$.
(We do not have to worry about the exact definition in the ambiguous
case when $\MSL(n)$ exactly equals some $p_M$; the constants in the
\PH{} system are such that this never will happen, see 
\refR{Rnotequal}
below. I have therefore
just chosen one version in the formulas below.)
The leap month rule is thus:
\begin{equation}\label{leaprule-dp}
\hskip-1em
\vbox {\narrower\narrower\narrower\noindent\em
A month where the mean sun does not pass any definition point
is a leap month.
}
\hskip-3em
\end{equation}

\begin{remark}\label{RChina}
The rule above for numbering the months is 
the same as in many Indian calendars \cite[Chapters 9 and 17]{CC}, 
and almost the same as the rule in the Chinese  calendar 
\cite[Chapters 16]{CC};
however, in these calendars, the
\dpp{s} are beginnings of the zodiacal signs, \ie, multiples of $30\grad=1/12$
while in the Tibetan calendar, the \dpp{s} are shifted, see \refS{SSdp}.
(In the Chinese calendar, the \dpp{s} are
called \emph{(major) solar terms}. 
In most (but not all) Indian calendars,
month 1 (New Year) is defined by the vernal equinox, \ie, $p_1=0$ in
our notation;
the exact rule in the Chinese calendar is that the
winter solstice ($270\grad=3/4$) occurs in month 11, which corresponds to
$p_{11}=3/4$, but there are some complications for other months.) 
Note however, that the (present)
Chinese and Indian calendars use the 
true motion of sun and moon, while the Tibetan uses the mean motion,
leading to regularly spaced leap months in the Tibetan calendar, but
not in the Chinese and Indian ones. (It also leads to skipped months
sometimes in the Indian calendars.)
Note further that the numbering of leap months differs in the Chinese
calendar, where a leap month is given the number of the preceding month.
\end{remark}

The longitude in \eqref{msl} is naturally taken modulo 1, \ie,
considering only the fractional part. But a moment's consideration
shows that the integer part shows the number of elapsed full circles
of the sun, \ie, the number of years; in this appendix we thus use $\MSL(n)$
for the 
real number defined by \eqref{msl}.
Let, as in \refS{Smonths},
$Y$ and $M$ be the year and the number of the month with \tmc{} $n$, 
and let the Boolean variable $\ell$ indicate whether the month is a leap month.
The rule \eqref{rule-dp} then yields the following relations determining
$(Y,M,\ell)$,
with $p_0=p_{12}-1$,
\begin{gather}
  Y-Y_0+p_{M-1}
<
\MSL(n)
\le
  Y-Y_0+p_{M},
\label{xa1}
\\
\ell =[\MSL(n+1)\le  Y-Y_0+p_{M}].
\label{xa2}
\end{gather}
Furthermore, the points $p_M$ are evenly spaced, so $p_M=p_0+M/12$.

\begin{remark}\label{Rs0}
The initial longitude $s_0$ is, as any longitude, really defined modulo 1,
\ie, only the fractional part matters. However, when we regard $\MSL(n)$ as a
real number, we have to make the right choice of integer part of $s_0$.
  Since the epoch is assumed to be in  year $Y_0$, with a month $M$
  satisfying $1\le M\le12$, taking $n=0$ in \eqref{xa1} shows that 
we must have
  \begin{equation}\label{s0krav}
	p_0<s_0\le p_{12}=p_0+1.
  \end{equation}
Equivalently, $\ga$ defined in \eqref{ga} below must satisfy
\begin{equation}\label{gakrav}
  0<\ga\le12.
\end{equation}
(This follows also by taking $n=0$ and $Y=Y_0$ in \eqref{b3}--\eqref{b5}
below, see \refR{Repoch-ga}.)
Consequently,
for the formulas for month numbers and leap months below,
we have to assume that the integer part of $s_0$ is chosen
such that \eqref{s0krav}--\eqref{gakrav} hold; 
this sometimes means adding 1 to the traditional value
(which does not affect the solar longitude seen as an angle, \ie{}
modulo 1).
\xfootnote{
For calculations one can use any $s_0$ and instead normalize $\ga$ in
\eqref{ga} modulo 12 so that \eqref{gakrav} holds.
}
\end{remark}

Let us for simplicity write $Y'=Y-Y_0$. Then \eqref{xa1} can be
rewritten
\begin{equation}
  \label{b0}
Y'+\frac{M-1}{12}+p_0
<
\MSL(n)
\le
Y'+\frac{M}{12}+p_0
\end{equation}
or
\begin{equation}  \label{b0x}
12Y'+M-1
<
12(\MSL(n)-p_0)
\le
12Y'+M.
\end{equation}
Hence, if we use \eqref{msl} and further define
\begin{equation}
  \label{ga}
\ga=12(s_0-p_0),
\end{equation}
we have
\begin{equation}
  \label{b1}
12Y'+M=\ceil{12(\MSL(n)-p_0)}
=\ceil{12s_1n+\ga}.
\end{equation}
Consequently, we can calculate $(Y,M)$ from $n$ by
\begin{align}
  x&=\ceil{12s_1n+\ga},
\label{b3}\\
M&= x \amod 12,
\label{b4}\\
Y&=\frac{x-M}{12}+Y_0.
\label{b5}
\end{align}
To complete the calculations of $(Y,M,\ell)$ from $n$, we 
find similarly from \eqref{xa2}, or simpler by \eqref{b1}
because a month is leap if and only if it gets
the same number as the following one,
\begin{equation}
  \label{b2}
\ell=\bigbool{\ceil{12s_1(n+1)+\ga}=\ceil{12s_1n+\ga}}.
\end{equation}

For the Tibetan value of $s_1$ in \eqref{s1q}, we have
$12s_1=\frac{65}{67}$,
and thus \eqref{b3} can be written
\begin{equation}\label{bxx3}
  x = \Ceil{\frac{65}{67}n+\ga}
=\Ceil{\frac{65 n + 67\ga}{67}}
=\Ceil{\frac{65 n + \gb}{67}},
\end{equation}
where we define the integer $\gb$ to be $67\ga$ rounded up, \ie,
\begin{equation}\label{gb}
  \gb=\ceil{67\ga}.
\end{equation}
Note that \eqref{bxx3} and \eqref{b4}--\eqref{b5} are the same as
\eqref{bx3}--\eqref{bx5}, which shows that the algorithmic
calculations in \refS{Smonths}
yield the same correspondence between $(Y,M)$ and \tmc{} $n$ as 
the rule \eqref{rule-dp} used in this appendix, provided the value of
$\gb$ is the same; in particular, the two methods yield the same leap months.

\begin{remark}\label{Repoch-ga}
By \eqref{b3}--\eqref{b5}, the epoch month with \tmc{} $n=0$ is month
$M=\ceil{\ga}$ in the epoch year. 
(Recall that we have $0<\ga\le12$ by \eqref{gakrav}, so this yields a value
$1\le M\le 12$. Conversely, we see again that \eqref{gakrav} is necessary
for $Y_0$ to be the epoch year.)
By \refR{Repoch}, the traditional choices of epoch always yield $M=2$ or 3
for $n=0$. There are thus only two possibilities:
either $1<\ga\le2$ and the epoch month with \tmc{} 0 is month 2
(so the nominal epoch month 3 has \tmc{} 1), or $2<\ga\le3$ and 
the \tmc{} is 0 for month 3. 
\end{remark}

\begin{remark}\label{Rnotequal}
To verify the assertion above that $\MSL(n)$ never equals some \dpp, note
that the calculations above show that this would happen if and only
if  $12s_1n+\ga$ would be an integer. Since $12s_1=\frac{65}{67}$,
this can happen only if $67\ga=804(s_0-p_0)$ is an integer (and in that case it
would happen for some $n$); we will see in \eqref{ga1S} below that
for the \PH{} version,
this is not the case ($804\,s_0$ is an integer by \eqref{s0}
but $804\, p_0$ is not).
Similarly, it will never happen for the other versions of the Tibetan
calendar with the \dpp{s} defined for them below.
\end{remark}

Conversely, 
given $(Y,M,\ell)$, we can find 
the \tmc{} $n$ by
the theory in this appendix
from \eqref{b0},
noting that if there are two possible values of $n$, then we should
choose the smaller one if $\ell=\true$ (a leap month)
and the larger one if
$\ell=\false$ (a regular month).
If $\ell=\false$, then \eqref{b0} and \eqref{msl} thus show that $n$ is the
largest integer such that 
\begin{equation}\label{qu6}
  s_1 n+s_0-p_0 \le Y'+\frac{M}{12}
\end{equation}
or, recalling \eqref{ga}, 
\begin{equation}\label{qu7}
  s_1 n \le Y'+\frac{M-\ga}{12}
=\frac{12(Y-Y_0)+M-\ga}{12};
\end{equation}
if $\ell=\true$, this value of $n$ should be decreased by 1.
Hence, in all cases, $n$ can be computed from $(Y,M,\ell)$ by
\begin{equation}\label{b6}
  n =\Floor{\frac{12(Y-Y_0)+M-\ga}{12s_1}}-\boolx{\ell}.
\end{equation}
(For the \PH{} version, it is easily verified directly
that this is equivalent to \eqref{c2}, using \eqref{MM} 
and \eqref{gbxga} below.)

An alternative formula, which has the advantage that if there is no
leap month $M$ in year $Y$, then $(Y,M,\true)$ gives the same result
as $(Y,M,\false)$, is given by
\begin{equation}\label{b6x}
  n =\Floor{\frac{12(Y-Y_0)+M-\ga-(1-12s_1)\boolx{\ell}}{12s_1}};
\end{equation}
to see this, note that if $\ell=\false$, then \eqref{b6x} and \eqref{b6}
give the same result, while if $\ell=\true$, then \eqref{b6x} gives 1
more than \eqref{b6} applied to $(Y,M-1,\false)$, which is the month
preceding leap month $M$ (also if $M=1$, when this really is
$(Y-1,12,\false)$).

Let us write $M'=12(Y-Y_0)+M$; this can be interpreted as a solar 
month count from the beginning of year $Y_0$
(or we can interpret month $M$ year $Y$ as month $M'$ year
$Y_0$); note that $\MM$ in \refS{Smonths} is 
given by, by \eqref{MM},
\begin{equation}
  \label{MMM}
\MM=M'-M_0=M'-3.
\end{equation}
Since $12s_1=\xfrac{65}{67}$, we can write \eqref{b6} as
\begin{equation}\label{c0}
  n 
=\Floor{\frac{67}{65}(M'-\ga)}-\boolx{\ell}
=\Floor{\frac{67M'-67\ga}{65}}-\boolx{\ell},
\end{equation}
and since ${67}M'$ is an integer, 
this value is not affected if $67\ga$ is replaced by the integer
$\gb=\ceil{67\ga}$, see \eqref{gb}.
We thus have
\begin{equation}\label{c1}
  n 
=\Floor{\frac{67M'-\gb}{65}}-\boolx{\ell}
=M'+\Floor{\frac{2M'-\gb}{65}}-\boolx{\ell}.
\end{equation}
There is a leap month $M$ in year $Y$ if and only if
the \tmc{} for $(Y,M,\false)$ jumps by 2 from the preceding regular month
$(Y,M-1,\false)$. By \eqref{c1}, this happens exactly when $2M'-\gb$
just has passed a multiple of 65, \ie, when 
$2M'-\gb\equiv \text{0 or 1}\pmod{65}$.
Thus, the general leap month rule \eqref{leaprule-dp} is equivalent to:
\begin{Rule}
There is a leap month $M$ in year $Y$ if and only if
\begin{equation}\label{leaprule-gb}
2M'\equiv \gb\text{ or\/ }\gb+1\pmod{65},
\end{equation}
where $M'=12(Y-Y_0)+M$.
\end{Rule}
Since this is the same as \eqref{leaprule-gb0}, we see again that
\eqref{leaprule-dp} leads to the same results as the rules in \refS{Smonths}.

\subsection{\PH{} definition points}\label{SSdp}
In the \PH{} version, 
the first \dpp{} $p_1$ is 
$23;6 \rr{60}$ mansions, \ie{} 
\begin{equation}
  \label{p1}
p_1=23,6 \rr{27,60}
=\frac{77}{90}
\end{equation}
and thus
\begin{equation}
  \label{p0}
p_0=p_1-\frac1{12}
=\frac{139}{180};
\end{equation}
in degrees, this is $p_0=278\grad$, $p_1=308\grad$, and so on, with
intervals of $30\grad$, \ie{} 8 degrees after the beginning of each sign, see
\cite{Henning}, 
\cite{Schuh-review}. 
\xfootnote{\label{fPHpoints}%
For other purposes, Phugpa astronomy
regards the winter solstice to be at solar longitude
$18,31,30\rr{27,60,60}=247/360=247\grad$, one degree earlier than
  $p_{11}=248\grad$, and similarly for the summer solstice,
see 
\cite[p.~322--328]{Henning},
\cite[On intercalary months]{kalacakra},
\cite[pp.~114--115]{Schuh} and 
\cite[pp.~223, 225]{Schuh-review}.
I have no explanation for this
difference, but it does not affect the calendar which uses the values above.
}

Recall that in the formulas above, the integer part of $s_0$ has to be
chosen such that \eqref{s0krav} holds, \ie, $139/180 < s_0 \le 319/180$.
This is satisfied by the values \eqref{s0} and \eqref{s0E1927} for the
epochs \eS{} and \eH, but for \eX{} the value $s_0=0$ in \eqref{s0E1987} has
to be replaced by 1. (Recall that this does not matter in \refS{Sastro}.)
The constant $\ga$ defined by \eqref{ga} equals thus
\begin{align}
  \label{ga1S}
\ga&=12(s_0-p_0)
=\frac{1832}{1005}
=1+\frac{827}{1005}
\qquad(\eS);
\intertext{or}
\ga&=12(s_0-p_0)
=\frac{1922}{1005}
=1+\frac{917}{1005}
\qquad(\eH);
\\
\ga&=12(s_0-p_0)
=\frac{41}{15}
\phantom{00}
=2+\frac{11}{15}
\phantom{00}
\qquad(\eX).   \label{ga1X}
\end{align}
Hence,  \eqref{gb} yields
\begin{align}
  \label{gb1}
\gb&=123  
\qquad(\eS),
\\
\gb&=129  
\qquad(\eH),
  \label{gb1H}
\\
\gb&=184  
\qquad(\eX)
  \label{gb1X},
\end{align}
in agreement with
\eqref{gb1S0}--\eqref{gb1X0}. 
Consequently, we have verified that the 
arithmetic calculations in \refS{Smonths} and the astronomical theory in
this appendix yield the same result.

Recall that in \refS{Smonths}, $\gb$ was defined by \eqref{gbgbx}, so we see,
using \eqref{ga} and \eqref{gb},
that the initial value $\gbx$ in the \tm{}
calculation \eqref{tm0} is related to the definition points by 
\begin{equation}\label{gbxga}
  \gbx
 = 184 - \gb 
=184 - \ceil{67\ga}
=184 + \floor{804(p_0-s_0)}.
\end{equation}
This yields directly
\begin{align}
  \label{gbxS=}
\gbx&=61
\qquad(\eS),
\\
\gbx&=55
\qquad(\eH),
  \label{gbxH=}
\\
\gbx&=0\phantom0
\qquad(\eX),
  \label{gbxX=}
\end{align}
in agreement with
\eqref{gbxS}--\eqref{gbxX}.

Note also that by \refR{Repoch-ga},
the values  of $\ga$  in \eqref{ga1S}--\eqref{ga1X}
show again that
for \eS{} and \eH{}, the epoch month with \tmc{} 0 is month 2,
while for \eX, it is month 3.

\subsection{\TS{} \dpp{s}}\label{ASTSdp}

The values for $\gb$ given in 
\refApp{ATS} are consistent with, \cf{} \eqref{p1} for \PH,
\begin{equation}
  \label{p1t}
p_1=23,1,30 \rr{27,60,60}
=\frac{307}{360}
=307\grad,
\end{equation}
which gives 
\begin{equation}\label{p0t}
p_0=\frac{277}{360}=277\grad.  
\end{equation}
To see this, note first that
the values of $s_0$ in \eqref{s0T1732} and \eqref{s0T1852} both have to be
increased by 1 in order to satisfy \eqref{s0krav}, see \refR{Rs0}.
Thus we now take, for E1732 and E1852, respectively,
\begin{align}
  s_0&
=1-\frac{5983}{108540}
=\frac{102557}{108540}, 
\\
  s_0&=
0,1,22,2,4,18 \rr{27,60,60,6,13,67}+1
=1+\frac{23}{27135},
\end{align}
which together with \eqref{p0t} yield
\begin{align}
\ga&=12(s_0-p_0)
=2+\frac{1903}{18090}
\qquad\text{(E1732)},
\\
\ga&=12(s_0-p_0)
=2+\frac{14053}{18090}
\qquad\text{(E1852)}.
\end{align}
By \eqref{gb} this yields 
the values  $\gb=142$ (E1732) and 187 (E1852), 
as said in \refApp{ATS}.

Note that $67\ga$ is not an integer, so (as for \PH{})
 $\MSL(n)$ never equals some \dpp, see \refR{Rnotequal}.

\begin{remark} 
  Any $p_1$ satisfying
  \begin{equation}
\frac{92432}{27\cdot4020}
<
p_1
<
\frac{92567}{27\cdot4020}
  \end{equation}
would give the same $\gb$ and thus by \eqref{leaprule-gb} the same calendar.
The value \eqref{p1t} is from
Henning [personal communication]; 
it is determined from calendars and not explicitly taken from any text.
\end{remark}

\subsection{Mongolian (New Genden) \dpp{s}}\label{ASMdp}

As far as I know, \YP{} does not discuss \dpp{s} explicitly for the New
Genden version, but the value
of $\gbx$ is consistent with 
increasing the \dpp{} $p_1$ from \eqref{p1} to
\begin{equation}
  \label{p1mongo}
p_1=23,9 \rr{27,60}
=\frac{463}{540}
=308\tfrac23\grad,
\end{equation}
and thus
\begin{equation}
  \label{p0mongo}
p_0
=\frac{209}{270}
=278\tfrac23\grad,
\end{equation}
which together with \eqref{s0genden}
yields, by \eqref{ga} and \eqref{gb},
\begin{equation}
  \ga =\frac{7724}{3015}
=2+\frac{1694}{3015}
\end{equation}
and $\beta=172\equiv 42\pmod{65}$, 
as said in \refApp{AMongo}.
(Again, note that $67\ga$ is not an integer, \cf{} \refR{Rnotequal}.)

In fact, any value of $p_1$ with
\begin{equation}
  \frac{689}{804} < p_1 < \frac{690}{804}
\end{equation}
gives the same $\gb$, and thus the same calendar.
The value \eqref{p1mongo} is just our choice within this range.

We do not know whether such definition points are used at all in the 
Mongolian calendar, or whether the rule \eqref{leaprule-G} is used as it is
without further justification.

\subsection{When the leap month is the second}\label{ASdp+}

In the Bhutanese version of the calendar (\refApp{ABhutan}, see also
\refApp{AKT}), 
a leap month takes the number of the \emph{preceding} month, so the leap
month is the second of the two months with the same number.
For such versions, the theory above has to be modified. 

The basic rules \eqref{rule-dp} and 
\eqref{leaprule-dp} remain unchanged, but they now imply
that the number of a month is the number of
the last \dpp{} before the end of the month.
Hence \eqref{xa1}--\eqref{xa2} are replaced by
\begin{gather}
  Y-Y_0+p_{M}
<
\MSL(n+1)
\le
  Y-Y_0+p_{M+1},
\label{xa1b}
\\
\ell =[\MSL(n)> Y-Y_0+p_{M}].
\label{xa2b}
\end{gather}
Note that \eqref{xa1b} is the same as \eqref{xa1} with $M$ and $n$ replaced
by $M+1$ and $n+1$;
\xfootnote{
This implies that \eqref{s0krav} and \eqref{gakrav} have to be slightly
modified to
  \begin{equation}
	p_1<s_0+s_1\le p_{13}=p_1+1
  \end{equation}
and, using \eqref{gabh}, 
\begin{equation}
  0<  
\gax=
\ga-\frac{2}{67}\le12;
\end{equation}
this is no difference in practice since the epochs traditionally are chosen at
the beginning of 
month 3 (\nag), so $\ga$ (or rather $\ga-2/67$ in this case)
is between 1 and 3, see \refR{Repoch}.
}
thus \eqref{b0}--\eqref{b1} hold with the same modification,
which yields,
recalling $12s_1=\frac{65}{67}$,
\begin{equation}
  \label{b1b}
  \begin{split}
12Y'+M &
=\Ceil{12s_1(n+1)+\ga}-1
=\Ceil{12s_1n+\ga-\tfrac2{67}}
.  
\end{split}
\end{equation}
Hence we now define, in addition to  \eqref{ga},
\begin{equation}
  \label{gabh}
\gax  =\ga-\tfrac2{67} =12(s_0-p_0)-\tfrac2{67}
\end{equation}
and have, instead of \eqref{b1},
\begin{equation}
  \label{b1bh}
12Y'+M
=\ceil{12s_1n+\gax}.
\end{equation}
Consequently, \eqref{b3}--\eqref{b2} hold with $\ga$ replaced by $\gax$.

We now change the definition \eqref{gb} to 
\begin{equation}\label{gbb}
  \gb=\ceil{67\gax}=\ceil{67\ga}-2;
\end{equation}
then \eqref{b1bh} 
can be written, \cf{} \eqref{bxx3},
\begin{equation}\label{bxx3bh}
12Y'+M = \Ceil{\frac{65}{67}n+\gax}
=\Ceil{\frac{65 n + 67\gax}{67}}
=\Ceil{\frac{65 n + \gb}{67}}.
\end{equation}
Hence, with our new definition of $\gb$, \eqref{bx3}--\eqref{bx5} are valid
and yield the same correspondence between $(Y,M)$ and true month $n$ as
before (but we have to remember that the leap month now is the later of two
months with the same number).

For the converse problem, to find the true month $n$ given $(Y,M, \ell)$,
note first that if the month is regular ($\ell=\false$), then the month
before has number $M-1$ (interpreted modulo 12) and \tmc{} $n-1$, and thus by
\eqref{xa1b} 
\begin{equation}
  \MSL(n)\le Y-Y_0+p_M<\MSL(n+1).
\end{equation}
Hence, $n$ is the largest integer such that \eqref{qu6}--\eqref{qu7} hold, 
just as in \refApp{Aleap}. 
For a leap month ($\ell=\true$), this value of $n$ now should be increased by 1.
Hence, in all cases, \eqref{b6} holds if $-[\ell]$ is changed to $+[\ell]$.
Similarly,
\eqref{b6x} is modified to 
\begin{equation}\label{b6xb}
  n =\Floor{\frac{12(Y-Y_0)+M-\ga+(1-12s_1)\boolx{\ell}}{12s_1}};
\end{equation}
as before this gives the same result for
$(Y,M,\true)$ and $(Y,M,\false)$ if there is no leap month $M$.
(If $\ell=\false$, \eqref{b6xb}  agrees with \eqref{b6}; if $\ell=\true$ it
gives 1 less than \eqref{b6} applied to $(Y,M+1,\false)$.)

Also \eqref{c0} holds with if $-[\ell]$ is changed to $+[\ell]$.
There is now a leap month $M$ in year $Y$ if and only if
the \tmc{} jumps by 2 from the regular month $(Y,M,\false)$ to the next
$(Y,M+1,\false)$. It follows from (the modified)
\eqref{c0}, \cf{} \eqref{c1} and the argument after it, that
this happens exactly
when $2(M'+1)-\ceil{67\ga}\equiv \text{0 or 1}\pmod{65}$.
With our new definition \eqref{gbb} of $\gb$,
this means that the rules \eqref{leaprule-gb} and  the equivalent
\eqref{leaprule-gb0} still hold. 

If there is a leap month $(Y,M)$, then \eqref{ixq} yields the intercalation
index for the preceding regular month, so for the leap month we have to
increase the result by $2$, yielding
$(2M'+\gbx-4)\bmod{65}$.
By \eqref{leaprule-gb}, this shows that \eqref{jeppe} is modifed when the
leap month comes after the 
regular month with the same number; 
a leap month now has intercalation index
\begin{equation}\label{jeppe+}
  (\gb+\gbx-4)\bmod65 
\qquad \text{or}\qquad
  (\gb+\gbx-3)\bmod65 .
\end{equation}

Finally, \eqref{cu}--\eqref{lys} were derived from \eqref{leaprule-gb0}, and
thus hold in the present case too.

\subsection{Bhutanese \dpp{s}}\label{ASBdp}

The value  $\gb=191\equiv 61\pmod{65}$ for the Bhutanese version of the calendar
in \refApp{ABhutan} (and the epoch used there) is
consistent
with the \dpp{s}  defined by, for example,
\xfootnote{
This is our own choice. Any $p_1$ with 
$\frac{690}{804} < p_1 <\frac{691}{804}$ yields the same result.
} 
\begin{equation}
  \label{p1B}
p_1=23,10,30 \rr{27,60,60}
=\frac{103}{120}
=309\grad
\end{equation}
and thus
\begin{equation}
  \label{p0B}
p_0=p_1-\frac{1}{12} = 20,55,30 \rr{27,60,60}
=\frac{31}{40}
=279\grad.
\end{equation}
Again, in accordance with \refR{Rs0} and \eqref{s0krav},
we have to add 1 to the value of $s_0$ given in \eqref{s0bh};
we thus now use
\begin{equation}
    s_0
= 0,24,10,50 \rr{27,60,60,67}+1 =1+\frac{1}{67}, \label{s0bh+}
\end{equation}
which together with \eqref{p0B}
yields, by \eqref{ga},
\begin{equation}
  \ga =\frac{1929}{670}
=2+\frac{589}{670}
\end{equation}
and 
thus by \eqref{gbb}
$\beta=193-2=191$ as asserted above.
By \eqref{gbx-B}, this is consistent with the value $\gbx=2$ used in
Bhutan (for the epoch above).

However, this is \emph{not} consistent with the original text by
Lhawang Lodro and published calendars, which give the definition points
as [Henning, personal communication] 
\begin{equation}
  \label{p1Bx}
p_1=23,15 \rr{27,60}
=\frac{31}{36}
=310\grad
\end{equation}
and thus
\begin{equation}
  \label{p0Bx}
p_0=p_1-\frac{1}{12} = 21,0 \rr{27,60}
=\frac{7}{9}
=280\grad,
\end{equation}
which together with \eqref{s0bh}
would yield, by \eqref{ga} and \eqref{gbb}, 
\begin{equation}
  \ga =\frac{572}{201}
=2+\frac{170}{201}
\end{equation}
and $\beta=189\equiv 59\pmod{65}$, which does not agree with \eqref{gbx-B},
nor with actual leap months.
\xfootnote{For example,
month 5 was a leap month in 2008. 
(This is shown by the Election act 
\cite[National Assembly, Acts]{BhutanNA}   
who is enacted the
``26th Day of the Second 5th Month of the Earth Male Rat Year
corresponding to the 28th Day of the 7th Month of the Year 2008''.)
A simple calculation shows
that at the beginning of the second (leap) month 5, the mean solar longitude 
\eqref{meansun} was 
$5,14,19,47 \rr{27,60,60,67}=13/67=69.851\grad$, 
so if the definition point 
$p_5=5,15\rr{27,60}=70\grad$, as implied by \eqref{p1Bx},
then the mean sun reaches $p_5$ just after the beginning of the month, which
therefore would not be a leap month. 
However, with \eqref{p1B} we have $p_5=69\grad$, passed at the end of the
preceding month, and the mean sun does not pass any \dpp{} during the month;
hence the month is (correctly) a leap month.
(The next \dpp{} is $p_6=99\grad$ and the mean sun has longitude $98.955\grad$
at the end of the month.)
}
I have no explanation for this. One possibility is that at some time, either
the definition points or the epoch value for true month and intercalation
index was adjusted, ignoring the fact that they in theory are connected and
that one cannot be changed without changing the other.
(Cf.\ the similar problem in \refF{fPHpoints}.)

Note also that, as said in \refApp{ABhutan}, the winter solstice 
is defined as when the mean solar longitude is
$250\grad$; this is the definition point $p_{11}$ if we use
\eqref{p1Bx}, but not if we use \eqref{p1B}.

\subsection{\Karana{} \dpp{s}}\label{ASKdp}

The leap month rule \eqref{leaprule-K} may seem to be in accordance 
with the definition points given by
\begin{equation}
  p_1=22,30 \rr{27,60} = \frac56
=300\grad.
\end{equation}
Thus the \dpp{s} are at the beginning of zodiacal signs. 
In particular,
the \dpp{} for \nag{} (\Caitra) is
$p_3=0\grad$, 
the first point of Aries.
\xfootnote{
This agrees with Indian calendars, see \refR{RChina}, and seems to be the most
natural choice.}
This would give
\begin{equation}
  p_0= \frac9{12}=\frac34
=270\grad
\end{equation}
and \eqref{s0K}, \eqref{ga} and \eqref{gbb} would yield
\begin{align}
  \ga =\frac{403}{135}=2+\frac{133}{135}
\end{align}
and
$\gb=199\equiv4\pmod{65}$, in accordance with \refApp{AKT}, see also
\eqref{jeppe+}. 

However, the theory of \dpp{s} above
is based on the standard value \eqref{s1} for the mean
solar motion $s_1$, and is thus logically inconsistent with the 
slightly larger value
\eqref{s1KT} used in the \KT; the mean sun will move a little faster, and if
the leap year rule \eqref{leaprule-K}
is used, the mean sun will advance relative
to the position indicated by the number of the month. 
(But the difference is
small (0.003\%), and amounts to about 1 sign (= 1 month) in 2600 years.) 
\xfootnote{
The relation \eqref{6765} is thus only an approximation in
the original \KT{} version, but it has been treated as exact in later
versions of the calendar.}
Conversely, if the leap month rule \eqref{leaprule-dp} is used with any fixed
\dpp{s} and $s_1$ is the \karana{} value \eqref{s1KT}, then the leap months
will not follow a strict 65 year cycle.
See further \cite[On intercalary months]{kalacakra}. 
In order to use the formulas in Section \ref{ASdp}, 
or
rather their modifications described in \refS{ASdp+}
with leap months numbered by the preceding month, 
we thus have to use the \siddhanta{} value
$s_1=65/804$ from \eqref{s1} and not the \karana{} value \eqref{s1KT},
which seems illogical.

As far as I know, \dpp{s} are not mentioned explicitly in the \KT.
Nevertheless, it is obvious that the idea was to insert leap months so
that the year on the average agrees with the tropical solar year, and thus
the sun has more or less the same longitude for the, say, first month any year.
The simple relation \eqref{6765} was at some stage chosen as an
approximation, perhaps believed to be exact, and the mean solar motion $s_1$
was at some stage chosen to be the \siddhanta{} value $65/804$; these
choices are logically connected, but I do not know whether that was
realized when these values were introduced in the calendar, and whether
these choices were made together or at different times (as the \karana{}
system suggests).

\clearpage

\section{Planets}
\label{Aplanet}

The positions of the planets are, 
for astrological purposes,
calculated by the following procedure, also based on the \KT.
The calculations yield the longitudes at the end of the calendar
(=solar) day.

In modern terms, the (geocentric) longitude of a planet is found by
first calculating the heliocentric longitude and the longitude of the
sun, and then combining them (using trigonometry).
The Tibetan calculations do effectively this, although the theory
behind (which is not explicitly mentioned) rather is an old Indian
system using epicycles similar to the Ptolemaic system \cite[p.~57]{Henning}.
As in the calculations of $true\_date$ in \refS{Sastro}, the
calculations use special tables as substitutes for trigonometrical
calculations. 
For a more detailed description of the traditional calculation,
including conversions between different radices, 
se \citet[Chapter II]{Henning}, which also includes discussions of
the astronomical background and geometry, and of the accuracy of the formulas.

The constants given below are 
(as the main text of the present paper)
for the \PH{} tradition; 
the \TS{} and other versions use their own
slightly different epoch values, see \cite[p.~340]{Henning} and 
\cite[Epoch data]{kalacakra}.


\subsection{General day}

The general day \tib{spyi zhag}
is a count of (solar) days since the epoch. It is thus simply given by
\begin{equation}
  general\_day=\JD-epoch.
\end{equation}
For \eH{} (as in \cite{Henning}), 
the epoch is JD 2424972 (Friday, 1 April 1927; RD 703547)
and thus
\begin{equation}
  general\_day=\JD-2424972.
\end{equation}

\begin{remark}
  \label{Rgeneral}
The traditional method to calculate the general day is to first find
the true month $n$ as
in Section \ref{Smonths}, 
see \eqref{tm0} and \eqref{mcrule-P} or \eqref{c2}, and
then use it to find the 
number of elapsed lunar days
since the epoch, $ld$ say,  by
$ld=30n+D$, where $D$ is the date.

Next, $ld$ is multiplied by the ratio $11135/11312$ in
\eqref{meanlunarday}
(the mean length of a lunar day). A small constant $gd_0$ is
subtracted; for \eH, $gd_0=2,178 \rr{64,707}=199/5656$.
This constant is the fraction of the solar day remaining at the end of
the mean lunar day at the epoch.
The difference $\frac{11135}{11312}ld-gd_0$ is rounded up to the
nearest integer, and one sets provisionally  
\begin{equation}
general\_day=\Ceil{\frac{11135}{11312}ld-gd_0}.  
\end{equation}
This gives an approximation of the general day, but may be off by a
day because only the mean motion of the moon is considered (this is
equivalent to ignoring the corrections in \eqref{truedate}). 
Hence the day of week is calculated by the simple formula
\begin{equation}
\label{wd}
  \bigpar{general\_day+wd_0} \mod 7
\end{equation}
where the constant $wd_0$ is the day of week at the epoch; for \eH, $wd_0=6$.
If the value in \eqref{wd} differs from the correct day of week, then
$general\_day$ is adjusted by $\pm1$ (the error cannot be larger) so
that \eqref{wd} becomes correct.
(Since this final check and correction has to be done, one could as
well ignore the subtraction of $gd_0$ above, but it is traditionally done.)
\end{remark}
 
                                   
\subsection{Mean heliocentric motion}

The mean heliocentric position of a planet is represented by an
integer for each plane called \emph{particular day} \tib{sgos zhag},
calculated as
\begin{equation}
  particular\_day
=
\begin{cases}
  (100\cdot general\_day+pd_0) \bmod R, & \text{Mercury}, \\
  (10\cdot general\_day+pd_0) \bmod R, & \text{Venus}, \\
  (general\_day+pd_0) \bmod R, & \hskip-12pt\text{Mars, Jupiter, Saturn},
\end{cases}
\end{equation}
where the modulus $R$ and the epoch value 
$pd_0$ are given in \refT{Tplanet}. The periods of the planets are
thus, exactly,
87.97, 224.7, 687, 4332 and 10766 days, respectively.
(Modern astronomical values are 87.9684, 224.695, 686.93, 4330.6,
and 10746.9 \cite[Table 15.6]{AA}.)

The particular day is 0 at the first point of Aries (\ie, when the
longitude is 0), and thus the mean heliocentric longitude is
\begin{equation}
  mean\_helio\_long
= \frac{particular\_day}{R} \pmod 1.
\end{equation}
This is traditionally expressed in the radices $\rr{27,60,60,6,r}$
with the final radix $r$ depending on the planet and given in
\refT{Tplanet}. 
(Note that $r$ always is a divisor of $R$.)

\begin{table}[!htpb]
\begin{tabular}
{l l l l l l}
& Mercury & Venus & Mars & Jupiter & Saturn\\
\hline
$R$ & 8797 & 2247 & 687 & 4332 & 10766 \\
$pd_0$ (\eH) & 4639 & 301 & 157 & 3964 & 6286 \\
$r$ & 8797 & 749 & 229 & 361 & 5383 \\
$birth\_sign$ & 11/18 & 2/9 &19/54 & 4/9 & 2/3 \\
trad.\ \rr{27,60}& 16,30 & 6,0 & 9,30 & 12,0 & 18,0 \\
\end{tabular}  
\caption{Constants for planets.}
\label{Tplanet}
\end{table}

\subsection{Mean longitude of the sun}
The mean longitude of the sun (at the end of the calendar day) is 
calculated from scratch and not using the
calculation of $mean\_sun$ (at the end of the lunar day) in
\refS{Sastro}, but the results are consistent
\cite[pp.~87--88]{Henning}.
(The \TS{} tradition uses the calculation of $mean\_sun$,
and in one version $true\_sun$  
\cite[p.~341]{Henning}.) 
The formula used is
\begin{equation}
  mean\_solar\_long=s_1'\cdot general\_day+s_0',
\end{equation}
where
\begin{equation}
  s_1'=0,4,26,0,93156 \rr{27,60,60,6,149209}
=\frac{18382}{6714405}
\quad\Bigpar{=\frac{11312}{11135}s_2}
\end{equation}
and the epoch value is, for \eH,
\begin{equation}
  s_0'=25,9,20,0,97440 \rr{27,60,60,6,149209}
=1-\frac{458772}{6714405}
.
\end{equation}

\begin{table}[!htpb]
\begin{tabular}
{l l l l l l}
\qquad\qquad& Mercury & Venus & Mars & Jupiter & Saturn\\
\hline
0 & 0 & 0 & 0 & 0 & 0 \\ 
1 & 10 & 5 & 25 & 11 & 22 \\ 
2 & 17 & 9 & 43 & 20 & 37 \\ 
3 & 20 & 10 & 50 & 23 & 43 \\ 
\end{tabular}  
\caption{Equation for planets.}
\label{Tplanet-equ}
\end{table}

\begin{table}[!htpb]
\begin{tabular}
{l l l l l l}
\qquad\qquad& Mercury & Venus & Mars & Jupiter & Saturn\\
\hline
0 & 0 & 0 & 0 & 0 & 0 \\ 
1 & 16 & 25 & 24 & 10 & 6 \\ 
2 & 32 & 50 & 47 & 20 & 11 \\ 
3 & 47 & 75 & 70 & 29 & 16 \\ 
4 & 61 & 99 & 93 & 37 & 20 \\ 
5 & 74 & 123 & 114 & 43 & 24 \\ 
6 & 85 & 145 & 135 & 49 & 26 \\ 
7 & 92 & 167 & 153 & 51 & 28 \\ 
8 & 97 & 185 & 168 & 52 & 28 \\ 
9 & 97 & 200 & 179 & 49 & 26 \\ 
10 & 93 & 208 & 182 & 43 & 22 \\ 
11 & 82 & 202 & 171 & 34 & 17 \\ 
12 & 62 & 172 & 133 & 23 & 11 \\ 
13 & 34 & 83 & 53 & 7 & 3 \\ 
\end{tabular}  
\caption{Final correction for planets.}
\label{Tplanet-corr}
\end{table}


\subsection{Slow longitude and step index}
The remaining calculations are based on the mean heliocentric
longitude and the mean solar longitude, but these quantities are
treated differently for the inner (or ``peaceful'') 
planets Mercury and Venus 
and for the outer (or ``wrathful'') planets Mars, Jupiter,
Saturn.
The reason is that the mean motion of an outer planet is given by the
mean heliocentric longitude, while the mean motion of an inner planet is
given by the mean longitude of the sun.
In both cases, this main term is called the \emph{mean slow longitude}
\tib{dal ba},
while the other quantity is called the \emph{step index} \tib{rkang 'dzin}.
In other words, for the inner planets 
\begin{align*}
  mean\_slow\_long&=mean\_solar\_long, \\
  step\_index&=mean\_helio\_long
\end{align*}
and for the outer planets
\begin{align*}
  mean\_slow\_long&=mean\_helio\_long, \\
  step\_index&=mean\_solar\_long.
\end{align*}

Next, \cf{} the calculations for the moon and sun in \refS{Sastro},
the anomaly is calculated by
\begin{equation}
  anomaly=mean\_slow\_long - birth\_sign \pmod1,
\end{equation}
where the ``birth-sign'' \tib{skyes khyim}
is given in \refT{Tplanet}, both as a
rational number and in the traditional form in mansions \rr{27,60}.
The anomaly is used to find the equation from
\begin{equation}
  equ=planet\_equ\_tab(12\cdot anomaly),
\end{equation}
where $planet\_equ\_tab(i)$ is given in \refT{Tplanet-equ} for $i=0,\dots,3$,
which extends by the symmetry rules
$planet\_equ\_tab(6-i)=planet\_equ\_tab(i)$ and
$planet\_equ\_tab(6+i)=-planet\_equ\_tab(i)$;
linear interpolation is used beween integer arguments.
Finally, the true slow longitude \tib{dal dag} is given by 
\begin{equation}
  true\_slow\_long=mean\_slow\_long - equ/(27\cdot 60).
\end{equation}

\subsection{Geocentric longitude}

The final step is to combine the true slow longitude and the step
index. First, the difference of these is found:
\begin{equation}
  \mathit{diff}=step\_index - true\_slow\_long.
\end{equation}
This is used to find a correction by another table look-up:
\begin{equation}
  corr=planet\_corr\_tab(27\cdot \mathit{diff}) \pmod1,
\end{equation}
where $planet\_corr\_tab(i)$ is given in \refT{Tplanet-corr} for
$i=0,\dots,13$, 
which extends by the symmetry rules
$planet\_corr\_tab(27-i)=-planet\_corr\_tab(i)$ and
$planet\_corr\_tab(27+i)=planet\_corr\_tab(i)$;
as always, linear interpolation is used beween integer arguments.

Finally the (geocentric, or \emph{fast}) longitude 
\tib{myur ba} is given by
\begin{equation}
  fast\_long=true\_slow\_long + corr/(27\cdot 60).
\end{equation}

\subsection{Rahu}
\emph{Rahu} is the name of the nodes of the lunar orbit, \ie, the
intersections of the orbit and the ecliptic. More precisely, the
ascending node is called the \emph{Head of Rahu} and the descending node is
called the \emph{Tail of Rahu}. In Tibetan (as in Indian) astrology, Rahu is
treated as a planet, or perhaps two planets.
Rahu is further essential for prediction of eclipses 
\cite[Chapter III]{Henning}. 

Rahu has a slow motion that is retrograde (\ie, to the west, with
decreasing longitude, unlike the real planets). In the Tibetan system,
the period is exactly
230 lunar months = 6900 lunar days.

\begin{remark}
  In calendar (solar) days, this is, \cf{} \eqref{meanlunarday},
\begin{equation}
\frac{11135}{11312} \cdot 6900 = \frac{19207875}{2828}
\approx 6792.035
\text{ days}
\end{equation}
or, \cf{} \eqref{meanyear},
\begin{equation}
\frac{11135}{11312} \cdot 
6900 s_2 =6900\cdot\frac{s_1}{30}=\frac{7475}{402}
\approx 18.5945
\text{ Tibetan years}.
\end{equation}
The modern astronomical value is 
6798 days = 18.61 Gregorian years
\cite[Table 15.4]{AA}.
\end{remark}

To find the position of Rahu at day $D$, month $M$, year $Y$, one
first calculates the  \tmc{}
$n$, for example by \eqref{MM} and \eqref{c2}. Next, the number $x$
of elapsed lunar days since the Head of Rahu had longitude 0 is
calculated by, \cf{} the calculation of $ld$ in \refR{Rgeneral}, 
\begin{equation}
  x=30(n+rd_0)+D = 30n+D+30rd_0 = ld+30rd_0,
\end{equation}
where $rd_0$ is an epoch value; $rd_0=187$ for \eH.
(For \eX, $rd_0=10$.)
Finally, the longitudes of the head and tail of Rahu are given by
\begin{align}
  rahu\_head\_long&=-\frac{x}{6900} \pmod 1, \\
  rahu\_tail\_long&=rahu\_head\_long +\frac12 \pmod 1.
\end{align}

\begin{remark}
  Since Rahu has a retrograde motion, these are decreasing. In
  traditional calculations, 
$-rahu\_head\_long=\xfrac{x}{6900} \pmod 1$ is called the longitude
  of the \emph{Source of Rahu}. 
The longitudes are traditionally expressed in $\rr{27,60,60,6,23}$, and
  $1/6900$ is then written as $0,0,14,0,12$.
\end{remark}

\clearpage

\section{Further astrological calculations}\label{Aastro}
As explained in \refS{Syear}, each year has a name consisting of an
element and an animal.
For astrological 
\xfootnote{See Footnote \ref{f:astro}.}
purposes (in the Chinese or elemental astrological
system), 
there are many further associations, assigning a year, month, lunar
day or solar day one of the 12 animals in \refT{T12},
one of the 5 elements in \refT{T5},
one of the 8 trigrams in \refT{T8},
or one of the 9 numbers $1,\dots,9$ in \refT{T9};
as shown in Tables \ref{T12}, \ref{T5}, \ref{T8}, \ref{T9}, these have further
associations to, for example, numbers, colours and directions.
There are also further attributes in the Indian system.
This appendix is only a brief introduction, and
only some of the attributes and their connections are mentioned here. 
See \citet{Tseng1,Tseng2} and
 \citet{Henning} for further
details.


\begin{table}[!p]
\begin{tabular}{c l l}
number & element & colour\\
\hline
1 & wood & green (blue) \\
2 & fire & red \\
3 & earth & yellow \\
4 & iron & white \\
5 & water & dark blue (black) \\
\end{tabular}  
\caption{The 5 elements and the associated numbers and colours.}
\label{T5}
\end{table}


\newlength{\yinwd}
\newlength{\yinht}
\newlength{\yindp}
\setlength{\yinwd}{6pt}
\setlength{\yinht}{0pt}
\setlength{\yindp}{1pt}
\newcommand\yina{\vrule height \yinht width \yinwd	depth\yindp}
\newcommand\yin{\hbox{\yina\hskip\yinwd\yina}}
\newcommand\yang{\hbox{\vrule height \yinht width 3\yinwd depth\yindp}}
\newcommand\yy[1]{\vskip2pt\ifodd#1\yang \else \yin \fi}
\newcommand\tri[3]{\vbox{\yy{#3}\yy{#2}\yy{#1}}}
\begin{table}[!htpb]
\begin{tabular}{l c l l l l l }
& binary & trigram & Tibetan & Chinese & direction & element\\
\hline
1 & 5 & \tri101 & li & l\'i & S & fire \\
2 & 0 & \tri000 & khon & k\=un & SW & earth \\
3 & 6 & \tri110 & dwa & du\`i & W & iron \\
4 & 7 & \tri111 & khen & qi\'an & NW & sky \\
5 & 2 & \tri010 & kham & k\v{a}n & N & water \\
6 & 1 & \tri001 & gin & g\`en & NE & mountain \\
7 & 4 & \tri100 & zin & zh\`en & E & wood \\
8 & 3 & \tri011 & zon &  x\`un & SE & wind
\end{tabular}  
\caption{The 8 trigrams with some attributes.
The ordering is the standard (King Wen, Later Heaven) order; the
numbering is perhaps not traditional, and the binary coding (reading
the trigrams bottom-up) is mathematically natural but not traditionally used.}
\label{T8}
\end{table}


\begin{table}[!htpb]
\begin{tabular}{ c l l l }
\qquad\quad& colour & element & direction\\
\hline
1 & white & iron & N \\
2 & black & water & SW  \\
3 & blue & water & E \\
4 & green & wood & SE \\
5 & yellow & earth & Centre \\
6 & white & iron & NW \\
7 & red & fire & W \\
8 & white & iron & NE \\
9 & red & fire & S 
\end{tabular}  
\caption{The 9 numbers and their attributes.}
\label{T9}
\end{table}

\begin{remark}
The order of the directions in \refT{T9} may seem jumbled,
yet there is method in it.
The 9 numbers are often arranged in a $3\times3$ square according to
the directions as in \refT{T3x3} (upside down the standard Western
orientation), and then the numbers form a magic square with all rows,
columns and diagonals summing to 15.
\end{remark}

\begin{table}[!htpb]
  \begin{tabular}{l | l | l}
	4 & 9 & 2 \\
\hline
	3 & 5 & 7 \\
\hline
	8 & 1 & 6 \\
  \end{tabular}
\qquad
  \begin{tabular}{l | l | l}
	SE & S & SW\\
\hline
	E & C & W\\
\hline
	NE & N & NW\\
  \end{tabular}
\caption{A magic square of numbers and their directions}
\label{T3x3}
\end{table}

In formulas below, $Y$ is the Gregorian number of the year, see
\refS{Syear}. By \refS{Syear}, year $Y$ has in the Chinese 60 year
cycle number
\begin{equation}
  \label{yz}
 (Y-3) \amod 60,
\end{equation}
and hence numbers $(Y-3) \amod 10$ and $(Y-3) \amod 12$ in the Chinese
10 and 12 year cycles.


\subsection{Attributes for years}
\subsubsection{Elements}
Each year 
is given 4 or 5 elements: 
the power element \tib{dbang thang}, 
life element \tib{srog},
body element \tib{lus},
fortune element \tib{klung rta},
and sometimes also the
spirit element \tib{bla}.

The \emph{power element} is the element associated to the celestial
stem, see \refT{T10}; it is thus repeated in a cycle of 10 years, with
each element repeated 2 consecutive years, in the standard
order wood, fire, earth, iron, water
(the order in \refT{T5}).
As said above, year $Y$ is $(Y-3) \amod 10$ in the Chinese 10 year
cycle, and by Tables \refand{T10}{T5}, its power element has number
\begin{equation}\label{power}
  \Ceil{\frac{(Y-3)\amod10}{2}}=\Ceil{\frac{Y-3}2} \amod 5.
\end{equation}

The \emph{life element} is repeated in a cycle of 12 years, and is
thus determined by the animal name of the year. The list is given in
\refT{Tlife}. Note that each third year is earth (the years $\equiv
2 \bmod 3$); the remaining four elements come repeated 2 years each, in the
same (cyclic) order wood, fire, iron, water as the power element.


\begin{table}[!tpb]
\begin{tabular}{r r l l l}
\rlap{year}\phantom{0} && animal & life element & spirit element\\
\hline
1 &10&  Mouse & water & iron\\
2 &11&  Ox & earth & fire\\
3 &12&  Tiger & wood & water\\
4 &1&  Rabbit & wood & water\\
5 &2&  Dragon & earth & fire\\
6 &3& Snake & fire & wood \\
7 &4& Horse & fire & wood\\
8 &5& Sheep & earth & fire\\
9 &6& Monkey & iron & earth\\
10 &7& Bird & iron & earth\\
11 &8& Dog & earth & fire\\
12 &9& Pig & water & iron\\
\end{tabular}  
\caption{The 12 year cycle of life elements.
The first number on each line shows the year mod 12 counted from the
  start of a Chinese cycle; the second shows the year mod 12 counted from the
start of a Prabhava cycle.}
\label{Tlife}
\end{table}

The \emph{fortune element} is repeated in a cycle of 4 years, 
in the order wood, water, iron, fire.
(Earth is not used. Note that
the order of the four used elements is different from the order used
for the power and life elements.)
Since 4 is a divisor of 12, the fortune element is determined by the
animal name of the year.
See \refT{Tfortune}. 


\begin{table}[!htpb]
\begin{tabular}{r r l l }
\rlap{year}\phantom{0} && animals & fortune element\\
\hline
1 &2&  Mouse, Dragon, Monkey & wood \\
2 &3&  Ox, Snake, Bird & water\\
3 &4&  Tiger, Horse, Dog & iron \\
4 &1&  Rabbit, Sheep, Pig & fire\\
\end{tabular}  
\caption{The 4 year cycle of fortune elements.
The first number on each line shows the year mod 4 counted from the
  start of a Chinese cycle; the second shows the year mod 4 counted from the
start of a Prabhava cycle.}
\label{Tfortune}
\end{table}

The \emph{body element} is calculated in two steps.
(See \citet{Henning} for traditional ways of doing the calculations.)
First, an element is determined by the animal name; I do not know any
name for this intermediate element so let us call it $x$. The element
$x$ is repeated in a cycle of 6 years, with the 3 elements wood,
water, iron repeated 2 years each, see \refT{T6}. 
Then count the number of steps from $x$ to the power element $y$, or
equivalently, replacing the elements by the corresponding numbers in
\refT{T5}, calculate $y-x\bmod 5$.
Finally, this difference determines the body element by
\refT{Tbody}.
(Note that the order in this table is not the standard one.)


\begin{table}[!htpb]
\begin{tabular}{r r l l c}
\rlap{year}\phantom{0} && animals & element & number\\
\hline
1 &4&  Mouse, Horse & wood & 1\\
2 &5&  Ox, Sheep & wood & 1\\
3 &6&  Tiger, Monkey & water & 5 \\
4 &1&  Rabbit, Bird & water & 5 \\
5 &2&  Dragon, Dog & iron & 4 \\
6 &3&  Snake, Pig & iron & 4 \\
\end{tabular}  
\caption{The 6 year cycle of the element $x$ 
and its number
in the body element calculation.
The first number on each line shows the year mod 6 counted from the
  start of a Chinese cycle; the second shows the year mod 6 counted from the
start of a Prabhava cycle.}
\label{T6}
\end{table}


\begin{table}[!htpb]
\begin{tabular}{ c l }
$y-x \mod 5$ & body element\\
\hline
0 & iron \\
1 & water \\
2 & fire \\
3 & earth \\
4 & wood \\
\end{tabular}  
\caption{The final step in the body element calculation.}
\label{Tbody}
\end{table}

Since both $x$ and the power element $y$ are repeated 2 consecutive years each,
the same is true for the body element. If we consider only the even
years, say, then $x$ is by \refT{T6} given by 
1, 0, $-1$, 1, 0, $-1$, $\dots\pmod5$, 
while $y$ by \refT{T10} is given by 1, 2, 3, 4, 5, 1, 2, 3, \dots.
Hence, $y-x\mod 5$ repeats in a cycle of length 15:
0, 2, 4, 3, 0, 2, 1, 3, 0, 4, 1, 3, 2, 4, 1
(with differences $+2,+2,-1,+2,+2,-1,\dots$).

Consequently, the body element repeats in a cycle of 30 years, with
each element repeated for 2 consecutive years.
The full cycle is given in \refT{T30}.

The \emph{spirit element} is the element preceding the life element in
the standard order, see \refT{Tlife}.
Thus it too is repeated in a cycle of 12 years, and is
determined by the animal name of the year. Each third year is fire 
(the years $\equiv 2 \bmod 3$) and the remaining four elements come
repeated 2 years each, in the standard (cyclic) order.

Since all elements for the year
have periods dividing 60, they all repeat in the
same order in each 60 year cycle as is shown in Table \ref{T60e}.


\subsubsection{Numbers}
Each year 
is also associated with a set of three numbers, 
in the range
$1,\dots,9$, the \emph{central number},
the \emph{life number} and the
\emph{power number}.
(These numbers are used for persons born that year. The central number is also
called the \emph{body number} or 
\emph{birth number}.)
The numbers are associated with elements and directions according to
\refT{T9}. 
The numbers decrease by 1 $\pmod9$ for each new year, and thus repeat
in a cycle of 9 years; they are given simply by, for Gregorian year
$Y$:
\begin{alignat}2
 & \text{central number}&& = (2-Y) \amod 9,
\\
 & \text{life number}&& 
= 
(\text{central number}-3) \amod 9
=(8-Y) \amod 9,
\\
 & \text{power number}&& 
= 
(\text{central number}+3) \amod 9
=(5-Y) \amod 9.
\end{alignat}

Since 9 does not divide 60, these numbers do not follow the 60 year
cycle. The period for repeating all elements and numbers is 180 years,
\ie, 3 cycles of 60 years.


\begin{table}[!phbt]
\begin{tabular}{r r l l}
\rlap{year}\phantom{0} && power & body element \\
\hline
1 & 28 & wood & iron\\ 
2 & 29 & wood & iron\\ 
3 & 30 & fire & fire\\ 
4 & 1 & fire & fire\\ 
5 & 2 & earth & wood\\ 
6 & 3 & earth & wood\\ 
7 & 4 & iron & earth\\ 
8 & 5 & iron & earth\\ 
9 & 6 & water & iron\\ 
10 & 7 & water & iron\\ 
11 & 8 & wood & fire\\ 
12 & 9 & wood & fire\\ 
13 & 10 & fire & water\\ 
14 & 11 & fire & water\\ 
15 & 12 & earth & earth\\ 
16 & 13 & earth & earth\\ 
17 & 14 & iron & iron\\ 
18 & 15 & iron & iron\\ 
19 & 16 & water & wood\\ 
20 & 17 & water & wood\\ 
21 & 18 & wood & water\\ 
22 & 19 & wood & water\\ 
23 & 20 & fire & earth\\ 
24 & 21 & fire & earth\\ 
25 & 22 & earth & fire\\ 
26 & 23 & earth & fire\\ 
27 & 24 & iron & wood\\ 
28 & 25 & iron & wood\\ 
29 & 26 & water & water\\ 
30 & 27 & water & water\\ 
\end{tabular}  
\caption{The 30 year cycle of the body element.
The first number on each line shows the year mod 30 counted from the
  start of a Chinese cycle; the second shows the year mod 30 counted from the
start of a Prabhava cycle.}
\label{T30}
\end{table}


\begin{longtable}
{r r l l l l l}
\rlap{year}\phantom{0} && name (power) &  life & body & fortune & spirit\\
\hline
1 & 58 & Wood--Mouse & water & iron & wood & iron\\ 
2 & 59 & Wood--Ox & earth & iron & water & fire\\ 
3 & 60 & Fire--Tiger & wood & fire & iron & water\\ 
4 & 1 & Fire--Rabbit & wood & fire & fire & water\\ 
5 & 2 & Earth--Dragon & earth & wood & wood & fire\\ 
6 & 3 & Earth--Snake & fire & wood & water & wood\\ 
7 & 4 & Iron--Horse & fire & earth & iron & wood\\ 
8 & 5 & Iron--Sheep & earth & earth & fire & fire\\ 
9 & 6 & Water--Monkey & iron & iron & wood & earth\\ 
10 & 7 & Water--Bird & iron & iron & water & earth\\ 
11 & 8 & Wood--Dog & earth & fire & iron & fire\\ 
12 & 9 & Wood--Pig & water & fire & fire & iron\\ 
13 & 10 & Fire--Mouse & water & water & wood & iron\\ 
14 & 11 & Fire--Ox & earth & water & water & fire\\ 
15 & 12 & Earth--Tiger & wood & earth & iron & water\\ 
16 & 13 & Earth--Rabbit & wood & earth & fire & water\\ 
17 & 14 & Iron--Dragon & earth & iron & wood & fire\\ 
18 & 15 & Iron--Snake & fire & iron & water & wood\\ 
19 & 16 & Water--Horse & fire & wood & iron & wood\\ 
20 & 17 & Water--Sheep & earth & wood & fire & fire\\ 
21 & 18 & Wood--Monkey & iron & water & wood & earth\\ 
22 & 19 & Wood--Bird & iron & water & water & earth\\ 
23 & 20 & Fire--Dog & earth & earth & iron & fire\\ 
24 & 21 & Fire--Pig & water & earth & fire & iron\\ 
25 & 22 & Earth--Mouse & water & fire & wood & iron\\ 
26 & 23 & Earth--Ox & earth & fire & water & fire\\ 
27 & 24 & Iron--Tiger & wood & wood & iron & water\\ 
28 & 25 & Iron--Rabbit & wood & wood & fire & water\\ 
29 & 26 & Water--Dragon & earth & water & wood & fire\\ 
30 & 27 & Water--Snake & fire & water & water & wood\\ 
31 & 28 & Wood--Horse & fire & iron & iron & wood\\ 
32 & 29 & Wood--Sheep & earth & iron & fire & fire\\ 
33 & 30 & Fire--Monkey & iron & fire & wood & earth\\ 
34 & 31 & Fire--Bird & iron & fire & water & earth\\ 
35 & 32 & Earth--Dog & earth & wood & iron & fire\\ 
36 & 33 & Earth--Pig & water & wood & fire & iron\\ 
37 & 34 & Iron--Mouse & water & earth & wood & iron\\ 
38 & 35 & Iron--Ox & earth & earth & water & fire\\ 
39 & 36 & Water--Tiger & wood & iron & iron & water\\ 
40 & 37 & Water--Rabbit & wood & iron & fire & water\\ 
41 & 38 & Wood--Dragon & earth & fire & wood & fire\\ 
42 & 39 & Wood--Snake & fire & fire & water & wood\\ 
43 & 40 & Fire--Horse & fire & water & iron & wood\\ 
44 & 41 & Fire--Sheep & earth & water & fire & fire\\ 
45 & 42 & Earth--Monkey & iron & earth & wood & earth\\ 
46 & 43 & Earth--Bird & iron & earth & water & earth\\ 
47 & 44 & Iron--Dog & earth & iron & iron & fire\\ 
48 & 45 & Iron--Pig & water & iron & fire & iron\\ 
49 & 46 & Water--Mouse & water & wood & wood & iron\\ 
50 & 47 & Water--Ox & earth & wood & water & fire\\ 
51 & 48 & Wood--Tiger & wood & water & iron & water\\ 
52 & 49 & Wood--Rabbit & wood & water & fire & water\\ 
53 & 50 & Fire--Dragon & earth & earth & wood & fire\\ 
54 & 51 & Fire--Snake & fire & earth & water & wood\\ 
55 & 52 & Earth--Horse & fire & fire & iron & wood\\ 
56 & 53 & Earth--Sheep & earth & fire & fire & fire\\ 
57 & 54 & Iron--Monkey & iron & wood & wood & earth\\ 
58 & 55 & Iron--Bird & iron & wood & water & earth\\ 
59 & 56 & Water--Dog & earth & water & iron & fire\\ 
60 & 57 & Water--Pig & water & water & fire & iron\\ 
\caption{The 60 year cycle of combinations of different elements.
The first number on each line shows the year mod 60 counted from the
  start of a Chinese cycle; the second shows the year mod 60 counted from the
start of a Prabhava cycle. The power element is the first part of the name.}
\label{T60e}
\end{longtable}


\subsection{Attributes for months}\label{ASattributes-months}
Each regular calendar month is given attributes as follows.
A leap month is given the same attributes as the regular
month with the same number.

\subsubsection{Animals}
The 12 months are assigned one each of the 12 animals, in standard
order but with different starting points in the \PH{} and \TS{}
traditions,
see \refS{Smonths}\ref{month-animal}.
The full lists are given in \refT{T12months}. 
(The \TS{} system is the same as the Chinese, see \eg{} 
\cite[\S 126]{Ginzel} and \cite{10000}.
\xfootnote{\label{f10000}%
The Chinese astrological system in \cite{10000} assigns animals,
  gender and elements to solar months, defined by the minor solar terms, and
  not to the (lunar) calendar months. However, there are different
  astrological traditions in China. 
})

\subsubsection{Gender}
Each month is given the gender (male or female) associated to its
animal. As Tables \refand{T12}{T12months} show, this simply means
(both in the \PH{} and \TS{} traditions) that odd-numbered months are
male and even-numbered female.
(This is in accordance with the general Chinese principle that odd
numbers are male (\emph{yang}) and even numbers female (\emph{yin}).)

\subsubsection{Elements}
Each month is assigned one of the 5 elements in \refT{T5}.

In the \TS{} version, the months cycle continuously through 
the cycle in \refT{T10}, year after year.
This combines with the 12 month cycle for animals
to a 60 month cycle, exactly as for years, see \refS{Syear} and 
Tables \ref{Tlosar} and \ref{T60}.
Since each element is repeated 2 months in the 10 month cycle, the element
of the first  month advances one step in the list \refT{T5} each year.
More precisely, 
the first  month (Tiger) Gregorian year $Y$ has 
month element number $(Y-2) \amod 5$.
(This is exactly as in the Chinese calendar, see \eg{} 
\cite[\S 126]{Ginzel} and \cite{10000}.)

In the \PH{} version, the Tiger month 
(the first in the Chinese calendar),
which is month 11 \emph{the preceding year}
(see \refT{T12months} and
\refS{Smonths}\ref{month-animal}),
gets the element following the
element of the year given in \refT{T10}, using the standard element
order in \refT{T5}. (The element following another, $x$, in this cyclic
order is called the \emph{son} of $x$.) 
By \eqref{power},
this is element $\ceil{(Y-1)/2} \amod 5$.
Having determined the element of the Tiger month, the elements for the
12 month period starting with it are assigned in the same pattern as
for years in \refT{T10}: each element is repeated 2 months, the
first male and the second female, and then followed by the next
element. (This is the same sequence as for \TS, but only for 12 month periods.)

We thus obtain the following formulas for the number of the element
for month $M$ year $Y$. 
\begin{description}
  \item[\PH] 
\begin{equation}
  \begin{cases}
 \lrpar{ \Ceil{\frac{Y-1}2}+\Floor{\frac{M+1}2}} \amod 5, & 1\le M\le 10,\\
\\
 \lrpar{ \Ceil{\frac{Y}2}+\Floor{\frac{M-11}2}} \amod 5, & 11\le M\le 12 .	
  \end{cases}
\end{equation}
\item[\TS]
\begin{equation}
\lrpar{ Y-2+\Floor{\frac{M-1}2}} \amod 5.
\end{equation}
\end{description}

\begin{table}[!htpb]
\begin{tabular}{r l l l}
month & \PH & \TS & gender\\
\hline
1 & Dragon & Tiger & male \\
2 & Snake & Rabbit & female \\
3 & Horse & Dragon & male \\
4 & Sheep & Snake & female \\
5 & Monkey & Horse & male \\
6 & Bird & Sheep & female \\
7 & Dog & Monkey & male \\
8 & Pig & Bird & female \\
9 & Mouse & Dog & male\\
10 & Ox & Pig & female\\
11 & Tiger & Mouse & male\\
12 & Rabbit & Ox & female
\end{tabular}  
\caption{The animals for the months.}
\label{T12months}
\end{table}

\subsubsection{Numbers}
Only \TS{} calendars give one of the 9 numbers to each month.
The number decreases by $1 \pmod 9$ for each month (except leap months),
and thus by $12\equiv 3\pmod 9$ for each year.
Month $M$ year $Y$ has number
\begin{equation}
 \bigpar{3-(12Y+M)}\amod 9.
\end{equation}

The \TS{} rules for animal, gender and element agree with the rules in the
Chinese calendar, see  \cite{10000}. Mongolia uses also the same rules, see
\refApp{AMongo}. 

\goodbreak 

\subsection{Attributes for lunar days}
\subsubsection{Animal}
Each lunar day has an animal. These repeat in the usual cycle of 12,
see \refT{T12}, with 
each odd-numbered month starting with Tiger (number 3 in \refT{T12}) and
each even-numbered month starting with Monkey (number 9 in \refT{T12});
since there are exactly $30\equiv6\pmod{12}$ lunar days in each month,
the animals thus repeat in a continuous cycle broken only at leap
months, where the animals are repeated in the same order in the leap
month and the following regular month and there is a discontinuity
between the two months.
The number (in \refT{T12}) of the animal for lunar day $D$ in month
$M$ is thus
\begin{equation}
  (D+30M+8) \amod 12
=
  (D+6M+8) \amod 12.
\end{equation}

\subsubsection{Element}
Each lunar day has an element; these cycle inside each month in a
cycle of 5 in the usual order given in \refT{T5} (without repetitions
as for years and months), beginning with the element following the
element of the month calculated above.
If the month has element $x$, then lunar day $D$ of the month thus has
element $(x+D)\amod5$.
(There is thus often a jump in the sequence at the beginning of a new month.)

\subsubsection{Trigram}
The trigrams for lunar days cycle in the usual cycle of 8 in
\refT{T8};
as for the animals, this is continuous across months except for leap
months.
A Tiger month begins with trigram number 1, Li.
Thus the trigram for lunar day $D$ in a month with animal number $A$
(in \refT{T12}) has number
\begin{equation}
  (D+30(A-3))\amod 8
=
(D-2A-2)\amod 8
=
(D+6A+6)\amod 8.
\end{equation}
Hence months Tiger, Horse, Dog begin with Li;
Rabbit, Sheep, Pig begin with Zin;
Mouse, Dragon, Monkey begin with Kham;
Ox, Snake, Bird begin with Dwa.

\subsubsection{Number}
The nine numbers (with their colours) for lunar days cycle forward
in the usual cycle of 9 in \refT{T9};
as for the animals, this is continuous across months except for leap
months.
A Tiger month begins with number 1 (white).
Thus the number for lunar day $D$ in a month with animal number $A$
(in \refT{T12}) has number
\begin{equation}
  (D+30(A-3))\amod 9
=
(D+3A)\amod 9.
\end{equation}
Hence, months Tiger, Snake, Monkey, Pig begin with 1 (white);
Mouse, Rabbit, Horse, Bird begin with 4 (green);
Ox, Dragon, Sheep, Dog begin with 7 (red).

\subsection{Attributes for calendar days}\label{ASattributes-day} 

As explained in Sections \ref{Sdays}, \refand{Scal}{Sweek}, each
calendar (solar) day has a number (the date) and a day of week.

Each (calendar) day is also given an element (and its colour according to
\refT{T5}, gender, animal, trigram and number
from the  
Chinese system; these are simple cyclic with periods 10, 2, 12, 8, 9,
respectively, with the elements repeated twice each as for years in
\refT{T10}.
At least the element, gender and animal are the same as in the Chinese
calendar for the same day, see \cite{10000}.  

\subsubsection{Element}
The element corresponds to number $\JD \amod 10$ in the (Chinese) cycle in
\refT{T10}.
The number of the element in \refT{T5} is thus
\begin{equation}\label{dayelement}
  \Ceil{\frac{\JD \amod 10}{2}}
=
  \Ceil{\frac{\JD}{2}} \amod 5.
\end{equation}

\subsubsection{Gender}
The gender is male when $ \JD$ is odd, and female when $\JD$ is even.
\xfootnote{
This is in accordance with the Chinese rule mentioned above that odd numbers
are male and even numbers female, but this is just a coincidence since
Julian Day Numbers were not invented as part of Chinese astrology.
}

\subsubsection{Animal}
The animal has number $(\JD+2) \amod12$ in the (Chinese) cycle in
\refT{T12}.

\subsubsection{Trigram}
The trigram has number $(\JD+2) \amod8$ in \refT{T8}.

\subsubsection{Number}
The number in \refT{T9} is $(-\JD) \amod9$. 

\begin{remark}
  \citet[pp.~208--209]{Henning} describe these using different,
and presumably traditional,
  numberings for the 10 day cycle of elements and the 12 day cycle of
  animals; his numbers (for the same element and animal as given
  above)
are $(\JD-2) \amod 10$ and $(\JD-2)\amod12$.
He further calculates the number as $10-((\JD+1)\amod9)$.
\end{remark}

The element calculated in \eqref{dayelement} is the power element of the day. 
Exactly as for years, see above, further elements (life, fortune,
body, spirit) can be calculated from the animal--element pair; these
elements thus follow a cycle of 60 days, which is equal to the cycle
in \refT{T60e}. A day has the elements given in \refT{T60e} on line
$(\JD -10) \amod 60$.

\begin{remark}
  \TS{} calendars use a different method to assign the nine numbers;
  the first Wood--Mouse day after the (true astronomical) winter
  solstice is 5 (yellow); then the numbers increase by $1 \pmod 9$ each
  day, until 
  the first Wood--Mouse day after the (true astronomical) summer
  solstice, which is 4 (green), and then the numbers decrease by $1 \pmod
  9$ each day for half a year until the first
Wood--Mouse day after the next winter solstice.
This requires accurate astronomical calculations of the solstices,
  which is a central part of the Chinese calendar system, but foreign
  to the Tibetan calendar calculations.
\end{remark}

\subsubsection{Elemental yoga}
In the Indian system, each day of week has an associated element from
the set \{earth, fire, water, wind\},
see \refT{T7}.
Further, each lunar mansion is also associated to an element from the
same set, see \citet[Appendix I]{Henning} for a list. Each calendar
day is thus given a combination of two elements, for the day of week
and for the lunar mansion (calculated in \eqref{weekday} and \eqref{mansion});
the order of these two elements does not matter and the combination is
regarded as an unordered pair. There are thus 10 possible different
combinations (\emph{yogas}), each having a name, see \cite[p.~204]{Henning}.

\clearpage

\newcommand\jour{\emph}
\newcommand\book{\emph}
\newcommand\vol{\textbf}
\newcommand\no{\unskip:}  
\def\no#1#2,{\unskip#2, no. #1,} 


\end{document}